\PassOptionsToPackage{monochrome}{xcolor}
\documentclass[12pt]{article}
\usepackage{amssymb}
\usepackage{mathrsfs}
\usepackage{amsfonts}
\usepackage{graphicx}
\usepackage{colortbl,dcolumn}
\usepackage{latexsym}
\usepackage{amsbsy,amsthm}
\usepackage{amsmath}
\usepackage{float}
\usepackage{bbm}
\usepackage[cmyk]{xcolor}
\usepackage{psfrag}
\usepackage{float}
\usepackage{subfigure}
\usepackage{textcase}
\usepackage{booktabs}
\usepackage{verbatim}
\usepackage{enumerate}
\usepackage{cite}
\usepackage{amsthm}
\usepackage[top=0.8in, bottom=0.8in, left=0.8in, right=0.8in]{geometry}
\usepackage{hyperref}
\numberwithin{equation}{section}

\newtheorem{assumption}{Assumption}[section]
\newtheorem{remark}[assumption]{Remark}
\newtheorem{lemma}[assumption]{Lemma}
\newtheorem{theorem}[assumption]{Theorem}   
\newtheorem{proposition}[assumption]{Proposition}

\hypersetup{colorlinks=true,                         
linkcolor=orange,
anchorcolor=orange,
citecolor=orange}

\allowdisplaybreaks[4] 

\begin{document}
\title
{
    {\color{blue}
    Perturbation estimates for order-one strong approximations of SDEs without 
    globally monotone coefficients}\footnote
    {
    This work is supported by Natural Science Foundation of China (12471394, 12071488, 12371417), 
    Postgraduate Scientific Research Innovation Project of Hunan Province (No: CX20230349) and Fundamental Research Funds for the Central Universities of Central South University (No: 2021zzts0480). The authors would like to thank Prof. Arnulf Jentzen for his helpful comments to improve the manuscript, when he was invited by Xiaojie Wang to visit Central South University in June of 2024.
    }
}

\author
{
    Lei Dai, Xiaojie Wang\footnote
    {
    Corresponding author.
    \newline
    \hspace*{0.4cm} Emails: \{dailei, x.j.wang7\}@csu.edu.cn; x.j.wang7@gmail.com.
    }
    \\
    \footnotesize  School of Mathematics and Statistics, HNP-LAMA, Central South University, Changsha, China
    \\
}

\maketitle
\begin{abstract}
{
    \rm\small
    To obtain strong convergence rates of numerical schemes, an overwhelming majority of existing works impose a global monotonicity condition on coefficients of SDEs. 
    {\color{black}
    Nevertheless, there are still many SDEs from applications that do not have globally monotone coefficients.}
    As a recent breakthrough, 
    the authors of [Hutzenthaler, Jentzen, Ann. Probab., 2020] originally presented a perturbation theory
    for stochastic differential equations (SDEs), 
    which is crucial to recovering strong convergence rates 
    of numerical schemes in a non-globally monotone setting.
    However, only a convergence rate of order $1/2$
    was obtained there for time-stepping schemes such as 
    a stopped increment-tamed Euler--Maruyama (SITEM) method.
    {\color{black}
    An interesting question arises, 
    also raised by the aforementioned work,
    as to whether a higher 
    convergence rate than $1/2$ can be obtained
    when higher order schemes are used.}
    The present work attempts to give a positive answer to this question. To this end, we develop 
    some new perturbation  estimates that are
    able to reveal the order-one 
    strong convergence of numerical methods.
    As the first application of the newly developed estimates,
    we identify the expected order-one pathwise uniformly
    strong convergence
    of the SITEM method for additive noise driven SDEs 
    and multiplicative noise driven second order SDEs
    with non-globally monotone coefficients. 
    As the other application, we propose and analyze a
    positivity preserving explicit Milstein-type method for 
    Lotka--Volterra competition model driven by multi-dimensional
    noise, with a pathwise uniformly strong convergence rate of 
    order one recovered under mild assumptions. 
    These obtained results are completely new
    and significantly improve the existing theory.
    Numerical experiments are also provided to confirm the theoretical findings. 
} 

\end{abstract}

\section{Introduction}	
    In order to describe the time evolution of many dynamical processes under random environmental effects, stochastic differential equations (SDEs) 
    \begin{equation} \label{eq:DaiLei-intro-SODE}
X_t = X_0 + \int_0^t f(X_s){\rm d}s +\int_0^t g(X_s){\rm d}W_s,\ t\in [0,T]
\end{equation}
    are widely used in various science and engineering fields such as finance, chemistry, physics and biology. In practice, the closed-form solutions of non-linear SDEs are rarely available and one usually falls back on their numerical approximations. 
    For SDEs possessing globally Lipschitz coefficients, the monographs  \cite{Milstein2004,1992Numerical} established a fundamental framework to analyze a batch of numerical schemes including typical  methods such as the explicit Euler--Maruyama (EM) method and explicit Milstein method. 
    Nevertheless, a majority of SDEs arising from applications have superlinearly growing coefficients and the globally Lipschitz condition is violated. A natural question thus arises as to whether the traditional numerical methods designed in the globally Lipschitz setting are still able to perform well when  used to solve SDEs with superlinearly growing coefficients.
%
    Unfortunately, the authors of \cite{hutzenthaler2011strong} gave a negative answer, by showing that the popularly used EM method is divergent in the sense of both strong and weak convergence, when used to
    solve a large class of SDEs with superlinearly growing coefficients.  
    Therefore special care must be taken to construct and analyze convergent numerical schemes in the absence of the Lipschitz regularity of coefficients.
    In recent years, a prospering growth of relevant works is devoted to numerical approximations 
    of SDEs with non-globally Lipschitz coefficients 
    (see, e.g., \cite{hutzenthaler2015numerical} and references therein).
    Roughly speaking, people either rely on implicit Euler/Milstein schemes \cite{alfonsi2013strong,higham2002strong,mao2013strong,andersson2017mean,higham2000mean,wang2020mean,wang2023mean,neuenkirch2014first,zong2018convergence} or 
    some explicit schemes based on modifications of the traditional explicit EM/Milstein methods \cite{Fang2020,Kelly2023,hutzenthaler2012strong,sabanis2016euler,wang2013tamed,kumar2019milstein,liu2013strong,mao2015truncated,beyn2017stochastic,Tretyakov2013fundamental,Li2019explicit,brehier2023approximation}.
%
%
To get the desired convergence rates of numerical schemes, a frequently used argument is based on Gronwall's lemma together with the popular global monotonicity condition, {\color{blue} for all $x,y \in \mathbb{R}^d$},
    \begin{equation}\label{eq:oneside_coef_condi}
        \begin{aligned}
            \langle x - y,f(x) - f(y)\rangle  + \tfrac{q}{2}\|g(x) - g(y)\|^2 \le {K}|x - y{|^2},
        \end{aligned}
    \end{equation}
    where $f$ and $g$ are the drift and diffusion coefficients of SDEs 
    \eqref{eq:DaiLei-intro-SODE}, respectively and 
    $K$ is a positive constant independent of $x,y$.
    Indeed, an overwhelming majority of existing works on convergence rates 
    carry out the error analysis under the
    global monotonicity condition \eqref{eq:oneside_coef_condi}.

    However, such a condition is still  
    restrictive and many momentous SDEs from applications fail to obey  \eqref{eq:oneside_coef_condi}. Examples include stochastic van der Pol oscillator, stochastic Lorenz equation, stochastic Langevin dynamics and stochastic Lotka--Volterra (LV) competition  model (see, e.g., \cite{hutzenthaler2015numerical,mao2007stochastic}).
    What if we did not have the condition \eqref{eq:oneside_coef_condi} available?
    In fact, the analysis of the convergence rates 
    of numerical schemes without the global monotonicity condition turns out to be highly non-trivial (see \cite{hutzenthaler2015numerical,Hutzenthaler2020}).
    As a recent breakthrough, 
    Hutzenthaler and Jentzen \cite{Hutzenthaler2020} made significant progress in this direction and originally
    developed a framework known as perturbation theory for SDEs beyond the global monotonicity assumption \eqref{eq:oneside_coef_condi}. 
    This theory, combined with exponential integrability properties of both numerical solutions and exact solutions (see \cite{cox2024local,hutzenthaler2018exponential}) enables one to reveal strong convergence rates of numerical schemes in a non-globally monotone setting.
    %
    Following this argument, the authors of \cite{Hutzenthaler2020} analyzed 
    the pointwise strong error
    \begin{equation}
        \sup_{t\in[0,T]}
        \big\| X_t - Y_t \big\|_{L^{r}(\Omega;\mathbb{R}^d)}
    \end{equation}
     of an explicit stopped increment-tamed EM (SITEM) method $\{Y_t\}_{t \in [0, T]}$ proposed by  \cite{hutzenthaler2018exponential} (cf. \eqref{eq:stop_tamed_method} therein),
     which was shown there to inherit exponential integrability properties of SDEs.
    Successfully, the authors identified the pointwise strong convergence rate of 
    order $\tfrac12$ for the SITEM method.
    {\color{black}
    An interesting question arises as to whether a higher 
    convergence rate than order $\tfrac12$ can be obtained
    when the considered SDEs are driven by additive noise
    or when high-order (Milstein-type) schemes 
    are used, which is also expected by \cite{Hutzenthaler2020}
    (see Remark 3.1 therein).}
    Unfortunately, following \cite[Theorem 1.2]{Hutzenthaler2020}, the convergence rates
    of any schemes would not exceed order $\tfrac12$,
which is nothing but the order of the H\"{o}lder regularity of the approximation process.
    %
    %

In the present article, we attempt to present some new perturbation  estimates that can be used to reveal the order-one pathwise uniformly strong convergence of numerical methods for several SDEs with non-globally monotone coefficients (see Theorem \ref{thm:main_thm}).
{\color{black}
Different from \cite{Hutzenthaler2020}, we use 
the $\rm It\hat{o}$ formula to expand the difference for the drift term, i.e., $f(Y_s) - a (s)$ in  Lemma \ref{lem:pre_estimate_theory} and rely on Burkholder-Davis-Gundy type inequalities to carefully treat the related terms (see estimates of $\mathbb{S}_2$ and $\mathbb{T}_2$ in Theorem \ref{thm:main_thm}). This approach essentially enables one to attain order one strong convergence for the error analysis of numerical schemes.
} 

As the first application of the newly developed estimates,
    we identify the expected order-one 
    pathwise uniformly strong convergence
    of the SITEM method for 
    some SDEs with non-globally monotone coefficients
    (see Theorem \ref{thm:stopped_tamed_method}
    and subsequent example models), including the additive noise driven SDEs (e.g., the stochastic Lorenz equation with additive noise, Brownian dynamics and Langevin dynamics)
    and
    {\color{blue}
    multiplicative noise driven second order SDEs (i.e., second order ordinary differential equations perturbed by multiplicative white noise) such as the stochastic van der Pol oscillator and stochastic Duffing-van der Pol oscillator:
    }
\begin{equation}
    \label{eq:intro_stopped_tamed_method-converg}
            \Big\| \sup_{t\in[0,T]}|X_t-Y_t| \Big\|_{L^{r}(\Omega;\mathbb{R})}
            \leq 
            C h.
        \end{equation}
Here $\{ Y_t \}_{ t \in [0, T]}$ is produced 
by the SITEM method \eqref{eq:stop_tamed_method},  {\color{blue} $r>0$ is an arbitrary constant and $h>0$ is the uniform stepsize.}
{\color{black}
These findings thus fill the gap left by \cite{Hutzenthaler2020}
and also significantly improve relevant convergence results in \cite{Hutzenthaler2020},
where the pointwise strong convergence rate of 
only order $\tfrac12$ was obtained for the SITEM method applied to these models.}

    As the other application, we propose and analyze a
    positivity preserving explicit Milstein-type method \eqref{eq:linear_mil_method} for stochastic
    LV competition model driven by multi-dimensional
    noises, with a pathwise uniformly strong convergence 
    rate of order one recovered (Theorem \ref{thm:LV_convergence_rate}).
    To the best of our knowledge, 
    this is the first paper to obtain 
    the order-one pathwise uniformly strong convergence of an explicit positivity preserving scheme
    for the general LV competition model.

    The paper is structured as follows. In the next section, 
    we introduce some notations and inequalities that may be used later. In Section \ref{sect:perturbation-estimates}, we present new perturbation estimates for SDEs beyond the global monotonicity assumption. Equipped with
    these estimates, we derive the order-one 
    strong convergence of the SITEM method for 
    some additive noise driven or second-order SDE models. 
     In Section \ref{sect:application3}, we propose and analyze an explicit Milstein method for the LV competition model with
    multi-dimensional noise. Some numerical experiments are provided to confirm the theoretical findings and a short conclusion is made in Section \ref{section:conclusion}.

\section{Preliminaries}
\label{sect:preliminary}

\subsection{Notations}
    Throughout this paper, unless otherwise specified, the following notations are used.
    Let $(\Omega,\mathcal{F},$\ $\{\mathcal{F}_t\}_{t\geq 0},\mathbb{P})$ be a complete probability space with a filtration $\{\mathcal{F}_t\}_{t\geq 0}$ fulfilling the usual conditions, that is,  the filtration is right continuous and increasing, and $\mathcal{F}_0$ contains all $\mathbb{P}$-null sets.
    Let $\{W_t\}_{t\geq0}$ be an $m$-dimensional standard Brownian motion defined on $(\Omega,\mathcal{F},\{\mathcal{F}_t\}_{t\geq 0},\mathbb{P})$. 
    For $a \in \mathbb{R}$, we define $\tfrac{a}{\infty}:=0$ and for $a \in \mathbb{R}\backslash \{ 0\},\tfrac{a}{0}:=\infty.$
    For a fixed integer number $d \geq 1$ and a vector $x \in \mathbb{R}^d,x^{(i)},i=1,...,d$ denotes the $i$-th component of $x$ and $|x|$ denotes the Euclidean norm induced by the vector inner product $\langle\cdot,\cdot\rangle$.
    For a matrix  $A \in \mathbb{R}^{d\times m},d,m \in \mathbb{N}$, $A^{(i)},i=1,...,m$ denotes the $i$-th column of $A$ and $A^{(ij)},i=1,...,d,j=1,...,m$ represents the element at $i$-th row and $j$-th column of $A$. 
    Let $A^*$ be the transpose of $A$ and $\|A\|:= \sqrt{\operatorname{trace}(A^{*}A)}$ be the Hilbert-Schmidt norm 
    induced by Hilbert-Schmidt inner product $\langle\cdot,\cdot\rangle_{HS}$. 
    For a random variable $\xi: \Omega \rightarrow \mathbb{H}$, where $\mathbb{H}$ is a separable Banach space endowed with norm $\|\cdot\|_\mathbb{H}$, $\mathbb{E}[\xi]$ denotes its expectation and for any $r>0$, $\|\xi\|_{L^r(\Omega;\mathbb{H})}:= (\mathbb{E}[\|\xi\|_{\mathbb{H}}^r])^{1/r}$. 
    For {\color{red} $f=(f^{(1)},...,f^{(d)})^*\in C^2(\mathbb{R}^d, \mathbb{R}^d)$}, we use $f^\prime(x)$ to denote the Jacobian matrix of $f(x)$, in which the $i$-th row is $f^{(i)\prime}(x):= \big(\nabla f^{(i)}(x)\big)^*:\mathbb{R}^d \rightarrow \mathbb{R}^{1 \times d}$. The notation $\operatorname{Hess}_x\big(f(x)\big)$ is a generalized Hessian matrix including $d$ components where the $i$-th is the Hessian matrix $\operatorname{Hess}_{x}\big(f^{(i)}(x)\big): \mathbb{R}^d \rightarrow \mathbb{R}^{d \times d}$ of $f^{(i)}(x)$. Further, for  $f: \mathbb{R}^d \rightarrow\mathbb{R}^d,g:\mathbb{R}^d \rightarrow\mathbb{R}^{d\times m}$ and $U\in C^2(\mathbb{R}^d,\mathbb{R})$, we denote 
    \begin{equation}
        (\mathcal{A}_{f,g} U )(x):=U'(x)f(x)+\tfrac{1}{2}\operatorname{trace}\left(g(x)g(x)^{*}\operatorname{Hess}_x\big(U(x)\big)\right).
    \end{equation}
    Let $T \in (0,\infty),N \in \mathbb{Z}^{+}$, let a uniform mesh 
    \begin{equation}\label{eq:uniform_mesh}
        0 = t_{0} < t_{1}< \cdots <t_{N}=T
    \end{equation}
    be constructed with the time step $h=T/N$, and for $s\in[0,T]$, {\color{black} define 
    \[
    \lfloor s \rfloor_N:=\sup_{n=0,...,N}\{t_n:t_n\leq s\}.
    \]}
    Moreover, we introduce a class of functions $\mathcal{C}^3_{\mathcal{D}}(\mathbb{R}^d,\mathbb{R})$
    as follows:
    \begin{equation}\label{eq:C_D_3_functions}
        \mathcal{C}^3_{\mathcal{D}}(\mathbb{R}^d,\mathbb{R}):=
        \left\{
        \Lambda \in C^2\big(\mathbb{R}^d,\mathbb{R}\big):
      \begin{array}{llll}
          &\text{Every element of}\ \operatorname{Hess}_x\big( \Lambda(x)\big)\ \text{is locally Lipschitz } 
          \\
          &\text{continuous and for}\  i \in \{1,2,3\}, a.s.\ (\text{Lebesgue }
          \\
          &\text{measure}) \ x\in{\mathbb{R}^d}, \text{there exist}\ p, c\geq3 \ \text{such that}
          \\
          &\|\Lambda^{[i]}(x)\|_{L^{[i]}(\mathbb{R}^d,\mathbb{R})}\leq c\big(1+| \Lambda ( x )|\big)^{1-i/p}.
      \end{array}\right\}.
    \end{equation}
    Here for $i=1,2,3$, we denote
    \begin{equation}
        \| \Lambda^{[i]}(x)\|_{L^{[i]}(\mathbb{R}^d,\mathbb{R})}:=\sup_{v_1,...,v_i\in{\mathbb{R}^d}\backslash \{0\}}\frac{\vert \Lambda ^{[i]}(x)(v_1,...,v_i)\vert}{|v_1|\cdots|v_i|},
    \end{equation}
    where
    \begin{equation} \label{eq:F-derivatives}
        \Lambda^{[i]}(x)(v_{1}, \ldots, v_{i})=\sum_{l_{1}, \ldots, l_{i}=1}^{d}\Big(\tfrac{\partial^{i} \Lambda }{\partial x_{l_{1}} \ldots \partial x_{l_{i}}} \Big)(x) \cdot v_{1}^{(l_{1})} \cdot v_{2}^{(l_{2})} \cdot \ldots \cdot v_{i}^{(l_{i})}.
    \end{equation}
    Note that $\mathcal{C}^3_{\mathcal{D}}(\mathbb{R}^d,\mathbb{R})$ forms a linear space containing a batch of functions such as
    $$
        \Lambda(x)=\Big(\sum_{i=1}^{d} x^{2c_i}_i \Big)^{r},\ c_i\geq 1,\ r \geq 1.
    $$
    {\color{black}
    For a metric space $(E,\rho)$, we say $\Lambda \in \mathcal{C}^1_{\mathcal{P}}(\mathbb{R}^d,E)$ with constants $K_{\Lambda}, c_{\Lambda}$,  
    {\color{red} if  $\Lambda \in  C\big(\mathbb{R}^d,E\big)$ and there exist constants $K_{\Lambda}, c_{\Lambda}\geq 0$} such that 
    \begin{equation}
        \rho( \Lambda(x), \Lambda(y))\leq K_{\Lambda}(1+|x|+|y|)^{c_{\Lambda}}|x-y|
    \end{equation}
    {\color{red}holds for all $x,y\in \mathbb{R}^d$.}
%
    One can easily see that if $\Lambda: \mathbb{R}^d \rightarrow \mathbb{R}$ is differentiable  and $\Lambda \in \mathcal{C}^1_{\mathcal{P}}(\mathbb{R}^d,\mathbb{R})$ with constants $K_{\Lambda}, c_{\Lambda}$, then there exists some constant $K_1$ such that {\color{blue} for all $x \in \mathbb{R}^d$},
    $$
        |\Lambda ( x ) | \leq K_1(1+|x|)^{c_ \Lambda + 1 }, |\Lambda'(x)|\leq K_1(1 + |x| )^{c_\Lambda}.
    $$
    {\color{red}
    Finally, we use $C$ (resp. $C$ with some subscripts)  to denote a generic positive constant independent of the time step (resp. dependent on the subscripts), which may differ from one place to another.
    }
    
\subsection{Burkholder-Davis-Gundy type inequalities}
In what follows we recall two Burkholder-Davis-Gundy type inequalities, which are frequently used in the subsequent analysis.
    \begin{lemma}\label{lem:BDG_continuous_inequality}
        Let $S:[0,T] \times \Omega \rightarrow  \mathbb{R}^{d \times m}$ be a predictable stochastic process satisfying $\mathbb{P}(\int_0^T\|S_t\|^2$ ${\rm d}t<\infty)=1$ and let $\{W_t\}_{t\geq0}$ be a $m$-dimensional standard Brownian motion. Then for any $p\geq2$,
        \begin{equation}
            \Big\| \sup_{t\in[0,T]}\Big|\int_0^t S_r{\rm d}W_r\Big|\Big\|_{L^p(\Omega;\mathbb{R})}\leq p\left(\int_0^T \sum\limits_{i=1}^{m}\|S^{(i)}_r\|^2_{L^p(\Omega;\mathbb{R}^d)}{\rm d}r\right)^{1/2}.
        \end{equation}
    \end{lemma}
    \begin{lemma}\label{lem:BDG_discrete_inequality}
        Let $M\in \mathbb{N}$ and $ S_1,...,S_M:\Omega \rightarrow \mathbb{R}$ be random variables satisfying $\sup_{i\in \{1,...,M\}}\|S_i$
        $\|_{L^2(\Omega;\mathbb{R})}$ $<\infty$ and for any $i\in \{1,...,M-1\},\ \mathbb{E}[S_{i+1}|S_1,...,S_i]=0$. Then for any $p\geq 2$, there exists some positive constant $C_p$  such that 
        \begin{equation}
            \|S_1+\cdots+S_M\|_{L^p(\Omega;\mathbb{R})}\leq C_p\left(\|S_1\|^2_{L^p(\Omega;\mathbb{R})}+\cdots+\|S_M\|^2_{L^p(\Omega;\mathbb{R})}\right)^{1/2}.
        \end{equation}
    \end{lemma}
The first lemma can be found in \cite[Lemma 2.7]{wang2013tamed} and 
the other one is quoted from \cite[Lemma 4.1]{hutzenthaler2011convergence}.

\section{New perturbation estimates for SDEs}
\label{sect:perturbation-estimates}
    In this section, let us focus on the following SDEs of It\^{o} type:
    \begin{equation} \label{eq:typical_sde}
      \left\{ 
        \begin{array}{l}
        X_t-X_0=\int_0^t f(X_s){\rm d}s +\int_0^t g(X_s){\rm d}W_s,\ t\in [0,T],\\
        X_0 = \xi_{X},
        \end{array} \right.
    \end{equation}
    where $f: \mathbb{R}^d \rightarrow\mathbb{R}^d$ stands for the drift coefficient, $g:\mathbb{R}^d \rightarrow\mathbb{R}^{d\times m}$ the diffusion coefficient and $\xi_{X}:\Omega \rightarrow  \mathbb{R}^d$ the initial data. 
    Also, we consider an approximation process given by
    \begin{equation} \label{eq:approxi_sde}
      \left\{ 
        \begin{array}{l}
        Y_t-Y_0=\int_0^t a(s){\rm d}s +\int_0^t b(s){\rm d}W_s,\ t\in [0,T],\\
        Y_0 = \xi_{Y},
        \end{array} \right.
    \end{equation}
    {\color{red} where $a:\Omega \times [0,T]\rightarrow \mathbb{R}^d$ and $b:\Omega \times [0,T]\rightarrow \mathbb{R}^{d\times m}$ are two stochastic processes that are integrable in the sense of Lebesgue integral and It\^{o} stochastic integral, respectively.}
    This can be regarded as a perturbation of the solution process of the original SDE  \eqref{eq:typical_sde}.
    For example, when the Euler-Maruyama method
    is used to  approximate \eqref{eq:typical_sde}
    on a uniform grid $\{t_n = n h\}_{0 \leq n \leq N}$ with stepsize $h = \tfrac{T}{N}$,
    one can get a continuous version of the approximation as
    \begin{equation} \label{eq:euler_approxi}
        Y_t-Y_0=\int_0^t f(Y_{\lfloor s \rfloor_N}){\rm d}s +\int_0^t g(Y_{\lfloor s \rfloor_N}){\rm d}W_s,\ t\in [0,T],\ Y_0 = \xi_{X},
    \end{equation}
    where in the notation of \eqref{eq:approxi_sde}
    we have $a(s)=f(Y_{\lfloor s \rfloor_N}),b(s)=g(Y_{\lfloor s \rfloor_N})$ and $\xi_{Y}=\xi_{X}$. 
%
The following lemma provides two estimates,
which will be essentially used later to obtain 
the desired perturbation estimates. 
{\color{black} The first assertion can be regarded as a modification of \cite[Proposition  2.9]{Hutzenthaler2020} and the second one is new.}
\begin{lemma}\label{lem:pre_estimate_theory}
    Let 
    $f:\mathbb R^d\rightarrow \mathbb R^d,g:\mathbb R^d\rightarrow \mathbb R^{d \times m} $ 
    be measurable functions. 
    {\color{black}
    Let $a: [0,T] \times \Omega \rightarrow \mathbb R^d,b: [0,T] \times \Omega \rightarrow \mathbb R^{d\times m}$ 
    be predictable stochastic processes and let $\tau:\Omega\rightarrow [0,T]$
    be a stopping time.
    }
    Let $\{ X_s \}_{s\in[0,T]}$ and 
    $\{ Y_s \}_{s\in[0,T]}$  be defined by \eqref{eq:typical_sde} and \eqref{eq:approxi_sde} with continuous sample paths, respectively.
    Assume that $\int_0^T |a(s)|+\|b(s)\|^2+|f(X_s)|+\|g(X_s)\|^2+|f(Y_s)|+\|g(Y_s)\|^2 {\rm{d}}s< \infty \ \mathbb{P}$-a.s. and for $\varepsilon \in (0,\infty),p\geq 2$ with $\mathbb{P}$-a.s.
        \begin{equation}\label{eq:assum_lem_inte_proper}
            \int_{0}^{\tau}\left[\tfrac{\left\langle X_{s}-Y_{s}, f\left(X_{s}\right)-f\left(Y_{s}\right)\right\rangle+\frac{(1+\varepsilon)(p-1)}{2}\left\|g\left(X_{s}\right)-g\left(Y_{s}\right)\right\|^{2}}{\left|X_{s}-Y_{s}\right|^{2}}\right]^{+} {\rm d} s<\infty.
        \end{equation}
    Then for any $u\in[0,T]$ it holds
    \begin{equation}\label{eq:pre_estimate_no_sup}
        \begin{aligned}
            &\sup_{t\in[0,u]}\Big\|\tfrac{{|{X_{t \wedge \tau}} - {Y_{t \wedge \tau}}{|}}}{{\exp (\int_0^{t \wedge \tau} {\frac{1}{p}{\eta_{p,r}}} {\rm d}r)}} \Big\|_{L^{p}(\Omega;\mathbb{R})}
            \\
            &\leq \sup_{t\in[0,u]} \Bigg[
            \|{\xi_X-\xi_Y}{\|^p_{L^{p}(\Omega;\mathbb{R}^d)}}
            +
            \underbrace{
            \mathbb{E}\bigg[\int_0^{t \wedge \tau} {\tfrac{{p|{X_s} - {Y_s}{|^{p - 2}}}}{{\exp (\int_0^s {{\eta_{p,r}}} {\rm d}r)}}\langle } {X_s} - {Y_s},\left(g(X_s)-b(s)\right){\rm d}{W_s}\rangle \bigg]
            }_{=:\color{red}\mathbb{S}_1}
            \\
            &\ \ \ +
            \underbrace{
            \mathbb{E}\bigg[\int_0^{t \wedge \tau} {\tfrac{p|{X_s} - {Y_s}{|^{p - 2}}{\langle {X_s} - {Y_s},f(Y_s)-a(s)\rangle
            }}{{\exp (\int_0^s {{\eta_{p,r}}} {\rm d}r)}}} {\rm d}s\bigg]
            }_{=:\color{red}\mathbb{S}_2}
            \\
            &\ \ \ +
            \underbrace{
            \mathbb{E}\bigg[\int_0^{t \wedge \tau} {\tfrac{p|{X_s} - {Y_s}{|^{p - 2}}{\frac{(p-1)(1+1/\varepsilon)}{2}\|g(Y_s)-b(s){\|^2}
            }}{{\exp (\int_0^s {{\eta_{p,r}}} {\rm d}r)}}} {\rm d}s\bigg]
            }_{=:\color{red}\mathbb{S}_3}\Bigg]^{1/p}.
    \end{aligned}
\end{equation}
Furthermore, we have
\begin{align*}
            &
            \Big\|\sup_{t\in[0,u]} \tfrac{|X_{t \wedge \tau}-Y_{t \wedge \tau}|}{{\exp (\int_0^{t \wedge \tau} {\frac{1}{2}{\eta_{2,r}}} {\rm{d}}r)}}\Big\|_{{L^{p}(\Omega;\mathbb{R})}}
            \\
            &\leq \Bigg[
            \|{\xi_X-\xi_Y}{\|^2_{L^{p}(\Omega;\mathbb{R}^d)}}+ 
            \underbrace{
            \Big\|\sup_{t\in[0,u]}\int_0^{t \wedge \tau} {\tfrac{{2\langle {X_s} - {Y_s},\left( {g({X_s}) -b(s)} \right){\rm{d}}W_s\rangle }}{{\exp (\int_0^s {{\eta_{2,r}}} {\rm{d}}r)}}}\Big\|_{{L^{p/2}(\Omega;\mathbb{R})}} 
            }_{=:\color{red} \mathbb{T}_1}\quad \quad \quad
            \\
            &\ \ \ +\underbrace{
            \Big\|\sup_{t\in[0,u]}\int_0^{t \wedge \tau} {\tfrac{{2\langle {X_s} - {Y_s},{f({Y_s}) - a(s)} \rangle }}{{\exp (\int_0^s {{\eta_{2,r}}} {\rm{d}}r)}}} {\rm{d}}s\Big\|_{{L^{p/2}(\Omega;\mathbb{R})}}
            }_{=: \color{red} \mathbb{T}_2}
            \\
            &\ \ \ +\underbrace{
            \Big\|\sup_{t\in[0,u]}\int_0^{t \wedge \tau} {\tfrac{{(1+\frac{1}{\varepsilon})\|g({Y_s}) - b(s){\|^2}}}{{\exp (\int_0^s {{\eta_{2,r}}} {\rm{d}}r)}}} {\rm{d}}s\Big\|_{{L^{p/2}(\Omega;\mathbb{R})}}
            }_{=:\color{red} \mathbb{T}_3}
            \Bigg]^{1/2}.
            \stepcounter{equation}
            \tag{\theequation}
            \label{eq:pre_estimate_with_sup}
\end{align*}
    Here for $z \geq2 $, we denote
    \begin{equation}\label{eq:pre_estimate_eta}
        \eta_{z,r}:= z \mathbbm{1}_{r\leq \tau} (\omega)
        \left[\tfrac{\langle {X_r} - {Y_r},f(X_r)-f(Y_r)\rangle+{\frac{(z-1)(1+\varepsilon)}{2}\|g(X_r)-g(Y_r){\|^2}
        }}{|{X_r} - {Y_r}|^2}\right]^{+}.
    \end{equation}  
\end{lemma}
\textbf{Proof}. 
{\color{blue} 
For fixed $p \geq 2$, it is easy to validate that $\eta_{p,r}: [0,T] \times \Omega$ defined by \eqref{eq:pre_estimate_eta} is well-defined due to \eqref{eq:assum_lem_inte_proper}. 
The $\rm It\hat{o}$ formula, the $\rm It\hat{o}$ product rule and the inequality $(a+b)^2 \leq (1+\varepsilon)a^2+(1+\tfrac{1}{\varepsilon})b^2$ yield 
\begin{align*}
        &\tfrac{{|{X_{t \wedge \tau}} - {Y_{t \wedge \tau}}{|^{p}}}}{{\exp (\int_0^{t \wedge \tau}{{\eta_{p,r}}} {\rm d}r)}} 
        =
        |{\xi_X} - {\xi_Y}{|^{p}}+ \int_0^{t \wedge \tau} {\tfrac{{p|{X_s} - {Y_s}{|^{p - 2}}}}{{\exp (\int_0^s {{\eta_{p,r}}} {\rm d}r)}}\langle } {X_s} - {Y_s},\left(g(X_s)-b(s)\right){\rm d}{W_s}\rangle  
        \\
        &\ \ \ +\int_0^{t \wedge \tau} {\tfrac{{p|{X_s} - {Y_s}{|^{p - 2}}\langle {X_s} - {Y_s},f(X_s)-a(s)\rangle  -|{X_s} - {Y_s}|^{p}\eta_{p,s}}}{{\exp (\int_0^s {{\eta_{p,r}}} {\rm d}r)}}} {\rm d}s
        \\
        &\ \ \ +\int_0^{t \wedge \tau} {\tfrac{{\frac{p(p - 2)}{2}|{X_s} - {Y_s}{|^{{p} - 4}}|({X_s} - {Y_s})^*(g(X_s)-b(s)){|^2} + \frac{p}{2}|{X_s} - {Y_s}{|^{p - 2}}\|g(X_s)-b(s){\|^2}}}{{\exp (\int_0^s {{\eta_{p,r}}} {\rm d}r)}}} {\rm d}s
        \\
        &\leq |{\xi_X} - {\xi_Y}{|^{p}} + \int_0^{t \wedge \tau} {\tfrac{{p|{X_s} - {Y_s}{|^{p - 2}}}}{{\exp (\int_0^s {\eta_{p,r}} {\rm d}r)}}\langle } {X_s} - {Y_s},\left(g(X_s)-b(s)\right){\rm d}{W_s}\rangle
        \\
        &\ \ \ +\int_0^{t \wedge \tau} {\tfrac{p|{X_s} - {Y_s}{|^{p - 2}}\left({\frac{(p-1)(1+\varepsilon)}{2}\|g(X_s)-g(Y_s){\|^2}+\langle {X_s} - {Y_s},f(X_s)-f(Y_s)\rangle
        }\right)-|{X_s} - {Y_s}|^{p}\eta_{p,s}}{{\exp (\int_0^s {{\eta_{p,r}}} {\rm d}r)}}} {\rm d}s
        \\
        &\ \ \ +\int_0^{t \wedge \tau} {\tfrac{p|{X_s} - {Y_s}{|^{p - 2}}\left({\frac{(p-1)(1+1/\varepsilon)}{2}\|g(Y_s)-b(s){\|^2}+\langle {X_s} - {Y_s},f(Y_s)-a(s)\rangle
        }\right)}{{\exp (\int_0^s {\eta_{p,r}} {\rm d}r)}}} {\rm d}s.
        \stepcounter{equation}\tag{\theequation}
        \label{eq:pre_estimate_ito_formula}
      \end{align*}}
{\color{black}
Then we arrive at  \eqref{eq:pre_estimate_no_sup} by
taking expectation of both sides of \eqref{eq:pre_estimate_ito_formula}.}
Similarly, to show \eqref{eq:pre_estimate_with_sup}, letting $p=2$ in \eqref{eq:pre_estimate_ito_formula} deduces that
    \begin{equation}\label{eq:pre_estimate_temp_with_sup}
        \begin{aligned}
            \tfrac{{|{X_{t \wedge \tau}} - {Y_{t \wedge \tau}}{|^2}}}{{\exp (\int_0^{t \wedge \tau} {\eta_{2,r}} {\rm{d}}r)}} &\leq |{\xi_X-\xi_Y}{|^2} + \int_0^{t \wedge \tau} {\tfrac{{2\langle {X_s} - {Y_s},\left( {g({X_s})-b(s)} \right){\rm{d}}W_s\rangle }}{{\exp (\int_0^s {{\eta_{2,r}}} {\rm{d}}r)}}}
            \\
            &\ \ \ +\int_0^{t \wedge \tau} {\tfrac{{2\langle {X_s} - {Y_s},{f({Y_s}) - a(s)}\rangle +(1+\frac{1}{\varepsilon})\|g({Y_s}) - b(s){\|^2}}}{{\exp (\int_0^s {{\eta_{2,r}}} {\rm{d}}r)}}} {\rm{d}}s.
        \end{aligned}
    \end{equation}
This clearly implies  \eqref{eq:pre_estimate_with_sup}, after
     taking supremum and $\|\cdot\|_{L^{p/2}(\Omega;\mathbb{R}^d)}$-norm of both sides of \eqref{eq:pre_estimate_temp_with_sup}.
\qed

As a consequence of Lemma \ref{lem:pre_estimate_theory}, we state the main
results of this section.

\begin{theorem}\label{thm:main_thm}
    Let $f:\mathbb R^d\rightarrow \mathbb R^d,g:\mathbb R^d\rightarrow \mathbb R^{d \times m}$ be measurable functions with $f\in C^2{(\mathbb R^d,\mathbb R^d})$, let $f\in \mathcal{C}^1_{\mathcal{P}}(\mathbb{R}^d,\mathbb{R}^{d})$ {\color{black} with constants $K_f,c_f$  and $g\in \mathcal{C}^1_{\mathcal{P}}(\mathbb{R}^d,\mathbb{R}^{d \times m})$ with constants $K_g, c_g$ }and
     let $a: [0,T] \times \Omega \rightarrow \mathbb R^d, b: [0,T] \times \Omega \rightarrow \mathbb R^{d\times m}$ be predictable stochastic processes.
     {\color{black}
     Let $\tau_N:\Omega\rightarrow [0,T]$ be a stopping time which may depend on $N$,} let $\{ X_s \}_{s\in[0,T]}$ and
    $\{ Y_s \}_{s\in[0,T]}$  be defined by \eqref{eq:typical_sde} and \eqref{eq:approxi_sde} with continuous sample paths, respectively and let the uniform mesh be constructed by \eqref{eq:uniform_mesh}. Assume that $\int_0^T |a(s)|+\|b(s)\|^2+|f(X_s)|+\|g(X_s)\|^2+|f(Y_s)|+\|g(Y_s)\|^2 {\rm{d}}s< \infty\ \mathbb{P}$-a.s. and for $\varepsilon \in (0,\infty)$ with $\mathbb{P}$-a.s.
        \begin{equation}
            \int_{0}^{\tau_N}\left[\tfrac{\left\langle X_{s}-Y_{s}, f\left(X_{s}\right)-f\left(Y_{s}\right)\right\rangle+\frac{(1+\varepsilon)(p-1)}{2}\left\|g\left(X_{s}\right)-g\left(Y_{s}\right)\right\|^{2}}{\left|X_{s}-Y_{s}\right|^{2}}\right]^{+} {\rm d} s<\infty.
        \end{equation}
    {\color{black}
    Moreover, let $K_{sup} >0$ be some constant that is independent of $N$, 
    let $\xi_0:=(\xi_X-\xi_Y) \in L^{p}(\Omega;\mathbb{R}^{d}), p\geq 4$, and suppose that }
    \begin{enumerate}[{\rm(}a{\rm)}.]
        \item 
            $\{s \leq \tau_N\} \in \mathcal{F}_{\lfloor s \rfloor_N};$
        \item 
            $\sup_{s\in[0,T]}\bigg\|\mathbbm{1}_{s\leq \tau_N}\left[\tfrac{\left\langle X_{s}-Y_{s}, f\left(X_{s}\right)-f\left(Y_{s}\right)\right\rangle+\frac{(1+\varepsilon)(p-1)}{2}\left\|g\left(X_{s}\right)-g\left(Y_{s}\right)\right\|^{2}}{\left|X_{s}-Y_{s}\right|^{2}}\right]^{+}\bigg\|_{L^{3p}(\Omega;\mathbb{R})}\leq K_{sup};$
        \item
        for any ${\color{red} i=1,...,d},\sup_{s\in[0,T]}\|\operatorname{Hess}_x(f^{(i)}(Y_s))\|_{L^{3p}(\Omega;\mathbb{R}^{d\times d})}\leq K_{sup}$ and
            \begin{equation}
                \begin{aligned}
                    &\sup_{s\in[0,T]}\|X_s\|_{L^{6pc_g \vee3pc_f\vee3p}(\Omega;\mathbb{R}^{d})} {\textstyle \bigvee}
                    \sup_{s\in[0,T]}\|Y_s\|_{L^{{6pc_g \vee3pc_f}}(\Omega;\mathbb{R}^{d})} {\textstyle \bigvee}
                    \\
                    &\sup_{s\in[0,T]}\|a(s)\|_{L^{3p}(\Omega;\mathbb{R}^{d})}{\textstyle \bigvee} \sup_{s\in[0,T]}\|b(s)\|_{L^{3p}(\Omega;\mathbb{R}^{d\times m})}
                    \leq K_{sup}.
                \end{aligned}
            \end{equation}
    \end{enumerate}
    \begin{enumerate}[\rm (1).]
    \item  {\color{black}Then for  any $u\in[0,T],v \in (0,\infty),q\in(0,\infty]$ with $\tfrac{1}{p}+\tfrac{1}{q}=\tfrac{1}{v}$, there exists some positive constant $C$ that might depend on $p,\varepsilon,d,m,T,K_f,c_f,K_g,c_g,K_{sup}$, but do not depend on $h$, such that
    }
        \begin{equation}\label{eq:perturbation_theory_no_sup}
            \begin{aligned}
            &\sup_{t\in[0,u]} \|X_{t\wedge\tau_N}-Y_{t \wedge \tau_N}\|_{L^{v}(\Omega;\mathbb{R}^d)} \leq \Bigg[ \|\xi_0{\|^p_{L^{p}(\Omega;\mathbb{R}^d)}}+C\bigg( h^p
            \\
            &\quad+\mathbb{E}\bigg[\int_0^{u} \mathbbm{1}_{s\leq \tau_N}\|g(Y_s)-b(s)\|^p {\rm d}s\bigg]
            +
             h^{\frac{p}{2}-1}\int_0^{u} \int_{\lfloor s \rfloor_N}^s \Big(\mathbb{E}\Big[ \mathbbm{1}_{r\leq \tau_N}\big\|g(Y_r)-b(r)\big\|^p\Big]\Big)^{\frac{1}{2}}{\rm d}r {\rm d}s
            \\
            &\quad +
            \mathbb{E}\bigg[\int_0^{u} \mathbbm{1}_{s\leq \tau_N}|f(Y_{\lfloor s \rfloor_N})-a(s)|^p {\rm d}s\bigg]
            +
            h^{\frac{p}{2}-1}\mathbb{E}\bigg[\int_0^{u} \int_{\lfloor s \rfloor_N}^s\mathbbm{1}_{r\leq \tau_N}|f(Y_{\lfloor r \rfloor_N})-a(r)|^p{\rm d}r {\rm d}s\bigg]
           \bigg) \Bigg]^{\frac{1}{p}}
            \\
            &\ \times\bigg\|\exp
            \bigg(\int_{0}^{\tau_N}\Big[\tfrac{\left\langle X_{s}-Y_{s}, f\left(X_{s}\right)-f\left(Y_{s}\right)\right\rangle+\frac{(1+\varepsilon)(p-1)}{2}\left\|g\left(X_{s}\right)-g\left(Y_{s}\right)\right\|^{2}}{\left|X_{s}-Y_{s}\right|^{2}}\Big]^{+} {\rm d} s\bigg)
            \bigg\|_{{L^{q}(\Omega;\mathbb{R})}}.
        \end{aligned}
        \end{equation}
    \item 
    In addition to the same settings as (1), if $g$ is Lipschitz, then {\color{black} it holds
    }
    \begin{equation}\label{eq:perturbation_theory_with_sup}
            \begin{aligned}
                &\big\|\sup_{t\in[0,u]} |X_{t \wedge \tau_N}-Y_{t \wedge \tau_N}|\big\|_{{L^{v}(\Omega;\mathbb{R})}} \leq \Bigg[\|\xi_0{\|^2_{L^{p}(\Omega;\mathbb{R}^d)}}+C\bigg(h^2
                \\
                &\quad+ 
                \int_0^u\Big\|\mathbbm{1}_{s\leq \tau_N}\|{g({Y_s})-b(s)}\| \Big\|^2_{{L^{p}(\Omega;\mathbb{R})}}{\rm d}s
                +
                h^{\frac{1}{2}}\int_0^{u} \Big( {{\int_{\lfloor s \rfloor_N}^s\Big\|\mathbbm{1}_{r\leq \tau_N}\|{g({Y_r})-b(r)}\| \Big\|^2_{{L^{p}(\Omega;\mathbb{R})}}{\rm d}r }}\Big)^{\frac{1}{2}}{\rm{d}}s
                \\
                &\quad+\int_0^u \Big\|\mathbbm{1}_{s\leq \tau_N}\big|{f(Y_{\lfloor s \rfloor_N}) - a(s)}\big|\Big\|^2_{{L^{p}(\Omega;\mathbb{R})}} {\rm{d}}s +
                h^{\frac{1}{2}}\int_0^{u} {{\int_{\lfloor s \rfloor_N}^s \Big\|\mathbbm{1}_{r\leq \tau_N}|f({Y_{\lfloor r \rfloor_N}}) - a(r)|\Big\|_{{L^{p}(\Omega;\mathbb{R})}} {\rm d}r}}{\rm{d}}s
                \bigg)\Bigg]^{\frac{1}{2}}
                \\
                &\ \times\bigg\|\exp
                \bigg(\int_{0}^{\tau_N}\left[\tfrac{\left\langle X_{s}-Y_{s}, f\left(X_{s}\right)-f\left(Y_{s}\right)\right\rangle+\frac{1+\varepsilon}{2}\left\|g\left(X_{s}\right)-g\left(Y_{s}\right)\right\|^{2}}{\left|X_{s}-Y_{s}\right|^{2}}\right]^{+} {\rm d} s\bigg)
            \bigg\|_{{L^{q}(\Omega;\mathbb{R})}}.
        \end{aligned}
    \end{equation}
    \end{enumerate}
    \end{theorem}
    \indent \textbf{Proof}. In the following exposition, we write $\eta_s$ to represent $\eta_{z,s}$ for short.
    For any $u\in[0,T]$, by the H{\"o}lder inequality one infers that
    \begin{equation}
    \begin{aligned}
       &\sup_{t\in[0,u]} \big\||X_{t\wedge\tau_N}-Y_{t \wedge \tau_N}|\big\|_{L^{v}(\Omega;\mathbb{R})} 
       \\
        &\leq 
       \sup_{t\in[0,u]} \Big\|\tfrac{|X_{t \wedge \tau_N} - Y_{t \wedge \tau_N}|}{\exp (\int_0^{t \wedge \tau_N} {\frac{1}{p}{\eta_r}} {\rm d}r)} \Big\|_{L^{p}(\Omega;\mathbb{R})}\cdot \Big\| {\exp \Big(\int_0^{\tau_N} {\tfrac{1}{p}{\eta_r}} {\rm{d}}r\Big)}\Big\|_{{L^{q}(\Omega;\mathbb{R})}}
    \end{aligned}
    \end{equation}
    and 
    \begin{equation}\label{eq:holder_inequal_with_sup}
        \begin{aligned}
            &\big\|\sup_{t\in[0,u]} |X_{t \wedge \tau_N}-Y_{t \wedge \tau_N}|\big\|_{{L^{v}(\Omega;\mathbb{R})}}
            \\
            &\leq    \Big\|\sup_{t\in[0,u]}\tfrac{|X_{t \wedge \tau_N} - Y_{t \wedge \tau_N}|}{\exp (\int_0^{t \wedge \tau_N} {\frac{1}{2}{\eta_r}} {\rm d}r)} \Big\|_{L^{p}(\Omega;\mathbb{R})}\cdot \Big\| {\exp \Big(\int_0^{\tau_N} {\tfrac{1}{2}{\eta_r}} {\rm{d}}r\Big)}\Big\|_{{L^{q}(\Omega;\mathbb{R})}}.
        \end{aligned}
    \end{equation}
    Therefore, it suffices to estimate terms {\color{red}
    $\mathbb{S}_1,\mathbb{S}_2,\mathbb{S}_3$} in \eqref{eq:pre_estimate_no_sup} and {\color{red} $\mathbb{T}_1,\mathbb{T}_2,\mathbb{T}_3$} in \eqref{eq:pre_estimate_with_sup}.
    Observe that the term $\mathbb{S}_1$ vanishes, due to the condition (c) in Theorem \ref{thm:main_thm} {\color{black} and  the growth of $g$}. For $\mathbb{S}_3$,
    by Young's inequality,
    \begin{equation}\label{eq:estimate_S_D_without_sup}
        \begin{aligned}
            \mathbb{S}_3 &\leq \mathbb{E}\bigg[\int_0^{t}C_{p,\varepsilon} \mathbbm{1}_{s\leq \tau_N}\tfrac{|{X_s} - {Y_s}|^{p}
            }{{\exp (\int_0^s {{\eta_r}} {\rm d}r)}}+C_p \mathbbm{1}_{s\leq \tau_N}\|g(Y_s)-b(s)\|^p {\rm d}s\bigg]
            \\
            &\leq C\int_0^{t}\mathbb{E}\bigg[\tfrac{|{X_{s\wedge \tau_N}} - {Y_{s\wedge \tau_N}}|^{p}
            }{{\exp (\int_0^{s\wedge \tau_N} {{\eta_r}} {\rm d}r)}}\bigg]{\rm d}s
            +
            C\mathbb{E}\bigg[\int_0^{t} \mathbbm{1}_{s\leq \tau_N}\|g(Y_s)-b(s)\|^p {\rm d}s\bigg].
        \end{aligned}
    \end{equation}
    Concerning $\mathbb{S}_2$, one can expand $f(Y_s)-f(Y_{\lfloor s \rfloor_N})$ by $\rm It\hat{o}$'s formula and then use the Young inequality and condition (c) in Theorem \ref{thm:main_thm} to infer
        \begin{align*}
            \mathbb{S}_2 &=
            \mathbb{E}\bigg[\int_0^{t \wedge \tau_N} {\tfrac{p|{X_s} - {Y_s}{|^{p - 2}}{\big\langle {X_s} - {Y_s},\int_{\lfloor s \rfloor_N}^s\langle f'(Y_r),a(r)\rangle +\frac{1}{2} \operatorname{trace}(b(r)^*\operatorname{Hess}_x(f(Y_r))b(r)){\rm{d}}r\big\rangle
            }}{{\exp (\int_0^s {{\eta_r}} {\rm d}r)}}} {\rm d}s\bigg]
            \\
            &\quad +\mathbb{E}\bigg[\int_0^{t \wedge \tau_N} {\tfrac{p|{X_s} - {Y_s}{|^{p - 2}}{\big\langle {X_s} - {Y_s},\int_{\lfloor s \rfloor_N}^s \langle f'(Y_r),b(r){\rm{d}}W_r\rangle\big\rangle
            }}{{\exp (\int_0^s {{\eta_r}} {\rm d}r)}}} {\rm d}s\bigg]
            \\
            &\quad +\mathbb{E}\bigg[\int_0^{t \wedge \tau_N} {\tfrac{p|{X_s} - {Y_s}{|^{p - 2}}{\langle {X_s} - {Y_s},f(Y_{\lfloor s \rfloor_N})-a(s)\rangle
            }}{{\exp (\int_0^s {{\eta_r}} {\rm d}r)}}} {\rm d}s\bigg]
            \\
            &\leq 
            \mathbb{E}\bigg[\int_0^{t \wedge \tau_N} {\tfrac{p|{X_s} - {Y_s}{|^{p - 1}}\big|\int_{\lfloor s \rfloor_N}^s\langle f'(Y_r),a(r)\rangle +\frac{1}{2} \operatorname{trace}(b(r)^*\operatorname{Hess}_x(f(Y_r))b(r)){\rm{d}}r\big|}{{\exp (\int_0^s {{\eta_r}} {\rm d}r)}}} {\rm d}s\bigg]
            \\
            &\quad+
            \mathbb{E}\bigg[\int_0^{t \wedge \tau_N} {\tfrac{p|{X_s} - {Y_s}{|^{p - 2}}{\big\langle {X_s} - {Y_s},\int_{\lfloor s \rfloor_N}^s \langle f'(Y_r),b(r){\rm d}W_r\rangle \big\rangle
            }}{{\exp (\int_0^s {{\eta_r}} {\rm d}r)}}} {\rm d}s\bigg]
            \\
            &\quad +\mathbb{E}\bigg[\int_0^{t \wedge \tau_N} {\tfrac{p|{X_s} - {Y_s}{|^{p - 1}}{|f(Y_{\lfloor s \rfloor_N})-a(s)|
            }}{{\exp (\int_0^s {{\eta_r}} {\rm d}r)}}} {\rm d}s\bigg]
            \\
            &\leq C\int_0^{t}\mathbb{E}\bigg[\tfrac{|{X_{s\wedge \tau_N}} - {Y_{s\wedge \tau_N}}|^{p}
            }{{\exp (\int_0^{s\wedge \tau_N} {{\eta_r}} {\rm d}r)}}\bigg]{\rm d}s
            +
            C\mathbb{E}\bigg[\int_0^{t} \mathbbm{1}_{s\leq \tau_N}|f(Y_{\lfloor s \rfloor_N})-a(s)|^p {\rm d}s\bigg]
            \\
            &\quad +
            \mathbb{E}\bigg[\int_0^{t \wedge \tau_N} {\tfrac{p|{X_s} - {Y_s}{|^{p - 2}}{\langle {X_s} - {Y_s},\int_{\lfloor s \rfloor_N}^s \langle f'(Y_r),b(r){\rm d}W_r\rangle \rangle
            }}{{\exp (\int_0^s {{\eta_r}} {\rm d}r)}}} {\rm d}s\bigg]+Ch^p.
            \stepcounter{equation}\tag{\theequation}
        \label{eq:estimate_SB_temp}
        \end{align*}
%
    To estimate the last but one term for $p\geq 4$, we expand the left item of the inner product by $\rm It\hat{o}$'s formula and $\rm It\hat{o}$'s product rule to obtain
\begin{align*}
            &\tfrac{p|{X_s} - {Y_s}|^{p - 2}({X_s} - {Y_s})}{\exp (\int_0^s {{\eta_r}} {\rm d}r)}
            \\
            &=
            \tfrac{p|{X_{\lfloor s \rfloor_N}} - {Y_{\lfloor s \rfloor_N}}|^{p - 2}({X_{\lfloor s \rfloor_N}} - {Y_{\lfloor s \rfloor_N}})}{\exp (\int_0^{\lfloor s \rfloor_N} {{\eta_r}} {\rm d}r)}
             +
            \int_{\lfloor s \rfloor_N}^s \tfrac{p|X_r-Y_r|^{p-2}(f(X_r)-a(r))}{\exp (\int_0^{r} {\color{red}\eta _{\iota} {\rm d}\iota})}{\rm{d}}r
            \\
            &\ \ \ +\int_{\lfloor s \rfloor_N}^s \tfrac{-p|X_r-Y_r|^{p-2}(X_r-Y_r)\eta_r}{\exp (\int_0^{r} {\color{red} \eta _{\iota} {\rm d}\iota})}{\rm{d}}r
            +
            \int_{\lfloor s \rfloor_N}^s \tfrac{p|X_r-Y_r|^{p-2}(g(X_r)-b(r))}{\exp (\int_0^{r} {\color{red}\eta _{\iota} {\rm d}\iota})}{\rm{d}}W_r
            \\
            &\ \ \ +\int_{\lfloor s \rfloor_N}^s \tfrac{p(p-2)(X_r-Y_r)}{\exp (\int_0^{r} {\color{red}\eta _{\iota} {\rm d}\iota})}|X_r-Y_r|^{p-4}\langle X_r-Y_r,f(X_r)-a(r)\rangle{\rm{d}}r
            \\
            &\ \ \ +\int_{\lfloor s \rfloor_N}^s \tfrac{p(p-2)(X_r-Y_r)}{\exp (\int_0^{r} {{\eta _{\iota}}} {\rm d}\iota)}|X_r-Y_r|^{p-4}\langle X_r-Y_r,\big(g(X_r)-b(r)\big){\rm{d}}W_r\rangle
            \\
            &\ \ \ +\int_{\lfloor s \rfloor_N}^s \tfrac{p(p-2)(p-4)(X_r-Y_r)}{2\exp (\int_0^{r} {{\eta _{\iota}}} {\rm d}\iota)}|X_r-Y_r|^{p-6}\big|(X_r-Y_r)\big(g(X_r)-b(r)\big)\big|^2{\rm{d}}r
            \\
            &\ \ \ +\int_{\lfloor s \rfloor_N}^s \tfrac{p(p-2)(X_r-Y_r)}{2\exp (\int_0^{r} {{\eta _{\iota}}} {\rm d}\iota)}|X_r-Y_r|^{p-4}\big\|g(X_r)-b(r)\big\|^2{\rm{d}}r
            \\
            &\ \ \ +\int_{\lfloor s \rfloor_N}^s \tfrac{p(p-2)|X_r-Y_r|^{p-4}}{\exp (\int_0^{r} {{\eta _{\iota}}} {\rm d}\iota)}\big(g(X_r)-b(r)\big)\big(g(X_r)-b(r)\big)^*(X_r-Y_r){\rm{d}}r.
         \stepcounter{equation}\tag{\theequation}
        \label{eq:expand_inner_product_without_sup}
    \end{align*}
Collecting some terms in \eqref{eq:expand_inner_product_without_sup} one can deduce
    \begin{align*}
       &\mathbb{E}\bigg[\int_0^{t \wedge \tau_N} {\Big\langle \tfrac{p|{X_s} - {Y_s}|^{p - 2}({X_s} - {Y_s})}{\exp (\int_0^s {{\eta_r}} {\rm d}r)}},\int_{\lfloor s \rfloor_N}^s \langle f'(Y_r), b(r){\rm d}W_r\rangle \Big\rangle {\rm d}s\bigg]
       \\
       &\leq
       \mathbb{E}\bigg[\int_0^{t \wedge \tau_N}\Big\langle \tfrac{p|{X_{\lfloor s \rfloor_N}} - {Y_{\lfloor s \rfloor_N}}|^{p - 2}({X_{\lfloor s \rfloor_N}} - {Y_{\lfloor s \rfloor_N}})}{\exp (\int_0^{\lfloor s \rfloor_N} {{\eta_r}} {\rm d}r)},\int_{\lfloor s \rfloor_N}^s \langle f'(Y_r), b(r){\rm d}W_r\rangle \Big\rangle {\rm d}s\bigg]
       \\
       &\ \ \ +
       \mathbb{E}\bigg[\int_0^{t \wedge \tau_N} \int_{\lfloor s \rfloor_N}^s \tfrac{p(p-1)|X_r-Y_r|^{p-2}|f(X_r)-a(r)|}{\exp (\int_0^{r} {{\eta _{\iota}}} {\rm d}\iota)}{\rm{d}}r\Big|\int_{\lfloor s \rfloor_N}^s \langle f'(Y_r), b(r){\rm d}W_r\rangle \Big| {\rm d}s\bigg]
       \\
       &\ \ \ +
       \mathbb{E}\bigg[\int_0^{t \wedge \tau_N}\Big\langle \int_{\lfloor s \rfloor_N}^s \tfrac{-p|X_r-Y_r|^{p-2}(X_r-Y_r)\eta_r}{\exp (\int_0^{r} {{\eta _{\iota}}} {\rm d}\iota)}{\rm{d}}r,\int_{\lfloor s \rfloor_N}^s \langle f'(Y_r), b(r){\rm d}W_r\rangle \Big\rangle {\rm d}s\bigg]
       \\
       &\ \ \ +
       \mathbb{E}\bigg[\int_0^{t \wedge \tau_N}\Big\langle \int_{\lfloor s \rfloor_N}^s \tfrac{p|X_r-Y_r|^{p-2}(g(X_r)-b(r))}{\exp (\int_0^{r} {{\eta _{\iota}}} {\rm d}\iota)}{\rm{d}}W_r,\int_{\lfloor s \rfloor_N}^s \langle f'(Y_r), b(r){\rm d}W_r\rangle \Big\rangle {\rm d}s\bigg]
       \\
       &\ \ \ +
       \mathbb{E}\bigg[\int_0^{t \wedge \tau_N}\Big\langle \int_{\lfloor s \rfloor_N}^s \tfrac{(X_r-Y_r)|X_r-Y_r|^{p-4}}{\frac{1}{p(p-2)}\exp (\int_0^{r} {{\eta _{\iota}}} {\rm d}\iota)}\big\langle X_r-Y_r,\big(g(X_r)-b(r)\big){\rm{d}}W_r\big\rangle,\int_{\lfloor s \rfloor_N}^s \langle f'(Y_r), b(r){\rm d}W_r\rangle \Big\rangle {\rm d}s\bigg]
       \\
       &\ \ \ +
       \mathbb{E}\bigg[\int_0^{t \wedge \tau_N} \int_{\lfloor s \rfloor_N}^s \tfrac{p(p-1)(p-2)|X_r-Y_r|^{p-3}}{2\exp (\int_0^{r} {{\eta _{\iota}}} {\rm d}\iota)}\big\|g(X_r)-b(r)\big\|^2{\rm{d}}r\Big|\int_{\lfloor s \rfloor_N}^s \langle f'(Y_r), b(r){\rm d}W_r\rangle \Big| {\rm d}s\bigg]
       \\
       &=:B_1+B_2+B_3+B_4+B_5+B_6.
     \stepcounter{equation}\tag{\theequation}
    \end{align*}
Next we estimate $B_i, i =1,2,...,6$ term by term.  Firstly, it is trivial to see $B_1=0$  by the condition (a) and (c) in Theorem \ref{thm:main_thm}.
For  $B_2$, using the fact $f\in \mathcal{C}^1_{\mathcal{P}}(\mathbb{R}^d,\mathbb{R}^{d})$, H{\"o}lder's inequality, Young's inequality and condition (c) in Theorem \ref{thm:main_thm}, we have 
    \begin{align*}
       B_2
       &\leq 
       C\int_0^{t } \int_{\lfloor s \rfloor_N}^s\mathbb{E}\bigg[ \mathbbm{1}_{s\leq \tau_N}\tfrac{|X_r-Y_r|^{p-1}(1+|X_r|+|Y_r|)^{c_f}}{\exp (\int_0^{r} {{\eta _{\iota}}} {\rm d}\iota)} 
       {\color{red} 
       \big|\int_{\lfloor s \rfloor_N}^s \big\langle f'(Y_{\iota}),b({\iota}){\rm{d}}W_{\iota}\big\rangle\big|
       }
       \bigg] {\rm{d}}r{\rm d}s
       \\
       &\ \ \ +
        C\int_0^{t} \int_{\lfloor s \rfloor_N}^s\mathbb{E}\bigg[\mathbbm{1}_{s\leq \tau_N}\tfrac{|X_r-Y_r|^{p-2}|f(Y_r)-f(Y_{\lfloor r \rfloor_N})|}{\exp (\int_0^{r} {{\eta _{\iota}}} {\rm d}\iota)}
        {\color{red} 
       \big|\int_{\lfloor s \rfloor_N}^s \big\langle f'(Y_{\iota}),b({\iota}){\rm{d}}W_{\iota}\big\rangle\big|
       }
        \bigg] {\rm{d}}r{\rm d}s
       \\
       &\ \ \ +
       C\int_0^{t} \int_{\lfloor s \rfloor_N}^s\mathbb{E}\bigg[\mathbbm{1}_{s\leq \tau_N}\tfrac{|X_r-Y_r|^{p-2}|f(Y_{\lfloor r \rfloor_N})-a(r)|}{\exp (\int_0^{r} {{\eta _{\iota}}} {\rm d}\iota)}
       {\color{red} 
       \big|\int_{\lfloor s \rfloor_N}^s \big\langle f'(Y_{\iota}),b({\iota}){\rm{d}}W_{\iota}\big\rangle\big|
       }
       \bigg] {\rm{d}}r{\rm d}s
       \\
       &\leq 
       C\int_0^{t} \sup_{r\in[{\lfloor s \rfloor_N},s]}\mathbb{E}\bigg[\mathbbm{1}_{s\leq \tau_N}\tfrac{|X_r-Y_r|^{p}}{\exp (\int_0^{r} {{\eta _{\iota}{\rm d}\iota)}}}\bigg]{\rm d}s
       \\
       &\ \ \ +
       C\int_0^{t} \bigg(\int_{\lfloor s \rfloor_N}^s \Big\| \mathbbm{1}_{s\leq \tau_N} \big(1+|X_r|+|Y_r|\big)^{c_f}
       {\color{red} 
       \big|\int_{\lfloor s \rfloor_N}^s \big\langle f'(Y_{\iota}),b({\iota}){\rm{d}}W_{\iota}\big\rangle\big|
       }
       \Big\|_{{L^{p}(\Omega;\mathbb{R})}}{\rm d}r\bigg)^p{\rm d}s
       \\
       &\ \ \ +
       C\int_0^{t} \bigg(\int_{\lfloor s \rfloor_N}^s \Big\| \mathbbm{1}_{s\leq \tau_N} \big|f(Y_r)-f(Y_{\lfloor r \rfloor_N})\big|\cdot
       {\color{red} 
       \big|\int_{\lfloor s \rfloor_N}^s \big\langle f'(Y_{\iota}),b({\iota}){\rm{d}}W_{\iota}\big\rangle\big|
       }
       \Big\|_{{L^{p/2}(\Omega;\mathbb{R})}}{\rm d}r\bigg)^{p/2}{\rm d}s
       \\
       &\ \ \ +
       C\int_0^{t} \bigg(\int_{\lfloor s \rfloor_N}^s \Big\| \mathbbm{1}_{s\leq \tau_N} \big|f(Y_{\lfloor r \rfloor_N})-a(r)\big|\cdot
       {\color{red} 
       \big|\int_{\lfloor s \rfloor_N}^s \big\langle f'(Y_{\iota}),b({\iota}){\rm{d}}W_{\iota}\big\rangle\big|
       }
       \Big\|_{{L^{p/2}(\Omega;\mathbb{R})}}{\rm d}r\bigg)^{p/2}{\rm d}s 
       \\
       &\leq 
       C\int_0^{t} \sup_{r\in[0,s]}\mathbb{E}\bigg[\tfrac{|X_{r \wedge \tau_N}-Y_{r \wedge \tau_N}|^{p}}{\exp (\int_0^{r \wedge \tau_N} {{\eta _{\iota}{\rm d}\iota)}}}\bigg]{\rm d}s +
       Ch^{\frac{p}{2}-1}\mathbb{E}\bigg[\int_0^{t} \int_{\lfloor s \rfloor_N}^s\mathbbm{1}_{r\leq \tau_N}|f(Y_{\lfloor r \rfloor_N})-a(r)|^p{\rm d}r {\rm d}s\bigg]+ Ch^p. 
             \stepcounter{equation}\tag{\theequation}
        \label{eq:B2_estimate_without_sup}
    \end{align*}
With the aid of H{\"o}lder's inequality, Young's inequality and condition (c) in Theorem \ref{thm:main_thm}, we estimate $B_3$ as follows:
\begin{equation}\label{eq:B3_estimate_without_sup}
    \begin{aligned}
       B_3&\leq
       \int_0^{t} \int_{\lfloor s \rfloor_N}^s \mathbb{E}\bigg[\mathbbm{1}_{s\leq \tau_N}\tfrac{p|X_r-Y_r|^{p-1}\eta_r}{\exp (\int_0^{r} {{\eta _{\iota}}} {\rm d}\iota)}
       {\color{red} 
       \big|\int_{\lfloor s \rfloor_N}^s \big\langle f'(Y_{\iota}),b({\iota}){\rm{d}}W_{\iota}\big\rangle\big|
       }
       \bigg] {\rm{d}}r{\rm d}s
       \\
       &\leq 
       C\int_0^{t} \sup_{r\in[{\lfloor s \rfloor_N},s]}\mathbb{E}\bigg[\mathbbm{1}_{s\leq \tau_N}\tfrac{|X_r-Y_r|^{p}}{\exp (\int_0^{r} {{\eta _{\iota}{\rm d}\iota)}}}\bigg]{\rm d}s
       \\
       &\ \ \ +
       C\int_0^{t} \bigg(\int_{\lfloor s \rfloor_N}^s \Big\| \mathbbm{1}_{s\leq \tau_N} \eta_r
       {\color{red} 
       \big|\int_{\lfloor s \rfloor_N}^s \big\langle f'(Y_{\iota}),b({\iota}){\rm{d}}W_{\iota}\big\rangle\big|
       }
       \Big\|_{{L^{p}(\Omega;\mathbb{R})}}{\rm d}r\bigg)^{p}{\rm d}s
       \\
       &\leq 
       C\int_0^{t} \sup_{r\in[0,s]}\mathbb{E}\bigg[\tfrac{|X_{r \wedge \tau_N}-Y_{r \wedge \tau_N}|^{p}}{\exp (\int_0^{r \wedge \tau_N} {{\eta _{\iota}{\rm d}\iota)}}}\bigg]{\rm d}s + Ch^{3p/2}. 
    \end{aligned}
\end{equation}
Using the property of stochastic integral, H{\"o}lder's inequality, Young's inequality and condition (c) in Theorem \ref{thm:main_thm} shows the estimate of
$B_4$:
    \begin{align*}
       B_4&=
       \int_0^{t}  \int_{\lfloor s \rfloor_N}^s \mathbb{E}\bigg[\mathbbm{1}_{s\leq \tau_N} \Big\langle \tfrac{p|X_r-Y_r|^{p-2}(g(X_r)-b(r))}{\exp (\int_0^{r} {{\eta _{\iota}}} {\rm d}\iota)},f'(Y_r)b(r)\Big\rangle_{HS}\bigg] {\rm d}r{\rm d}s
       \\
       &\leq \int_0^{t}  \int_{\lfloor s \rfloor_N}^s \mathbb{E} \bigg[\mathbbm{1}_{s\leq \tau_N}\tfrac{p|X_r-Y_r|^{p-2}\|g(X_r)-g(Y_r)\|}{\exp (\int_0^{r} {{\eta _{\iota}}} {\rm d}\iota)}\big\|f'(Y_r)b(r)\big\|\bigg] {\rm d}r{\rm d}s
       \\
       &\ \ \ +
       \int_0^{t}  \int_{\lfloor s \rfloor_N}^s \mathbb{E}\bigg[\mathbbm{1}_{s\leq \tau_N} \tfrac{p|X_r-Y_r|^{p-2}\|g(Y_r)-b(r)\|}{\exp (\int_0^{r} {{\eta _{\iota}}} {\rm d}\iota)}\big\|f'(Y_r)b(r)\big\|\bigg] {\rm d}r{\rm d}s
       \\
       &\leq 
       C\int_0^{t} \sup_{r\in[{\lfloor s \rfloor_N},s]}\mathbb{E}\bigg[\mathbbm{1}_{s\leq \tau_N}\tfrac{|X_r-Y_r|^{p}}{\exp (\int_0^{r} {{\eta _{\iota}{\rm d}\iota)}}}\bigg]{\rm d}s
       \\
       &\ \ \ +
       C\int_0^{t} \bigg(\int_{\lfloor s \rfloor_N}^s \Big\| \mathbbm{1}_{s\leq \tau_N} \big(1+|X_r|+|Y_r|\big)^{c_g}\big\|f'(Y_r)b(r)\big\|\Big\|_{{L^{p}(\Omega;\mathbb{R})}}{\rm d}r\bigg)^p{\rm d}s
       \\
       &\ \ \ +
       C\int_0^{t} \bigg(\int_{\lfloor s \rfloor_N}^s \Big\| \mathbbm{1}_{s\leq \tau_N} \big\|g(Y_r)-b(r)\big\|\big\|f'(Y_r)b(r)\big\|\Big\|_{{L^{p/2}(\Omega;\mathbb{R})}}{\rm d}r\bigg)^{p/2}{\rm d}s
       \\
       &\leq 
       C\int_0^{t} \sup_{r\in[0,s]}\mathbb{E}\bigg[\tfrac{|X_{r \wedge \tau_N}-Y_{r \wedge \tau_N}|^{p}}{\exp (\int_0^{r \wedge \tau_N} {{\eta _{\iota}{\rm d}\iota)}}}\bigg]{\rm d}s 
       \\
       & \ \ \
       +
       Ch^{\frac{p-2}{2}}\int_0^{t} \int_{\lfloor s \rfloor_N}^s \bigg(\mathbb{E}\Big[ \mathbbm{1}_{r\leq \tau_N}\big\|g(Y_r)-b(r)\big\|^p\Big]\bigg)^{1/2}{\rm d}r {\rm d}s+ Ch^p. 
        \stepcounter{equation}\tag{\theequation}
        \label{eq:B4_estimate_without_sup}
    \end{align*}
In a similar way, one can handle $B_5$ as follows:
\begin{equation}\label{eq:B5_estimate_without_sup}
    \begin{aligned}
       B_5
       &\leq \int_0^{t}\int_{\lfloor s \rfloor_N}^s \mathbb{E}\bigg[\mathbbm{1}_{s\leq \tau_N}\big\| \tfrac{|X_r-Y_r|^{p-2}}{\exp (\int_0^{r} {{\eta _{\iota}}} {\rm d}\iota)}\big(g(X_r)-b(r)\big)\big\|\big\|f'(Y_r)b(r)\big\|\bigg] {\rm d}r{\rm d}s
       \\
       &\leq 
       C\int_0^{t} \sup_{r\in[0,s]}\mathbb{E}\bigg[\tfrac{|X_{r \wedge \tau_N}-Y_{r \wedge \tau_N}|^{p}}{\exp (\int_0^{r \wedge \tau_N} {{\eta _{\iota}{\rm d}\iota)}}}\bigg]{\rm d}s 
       \\ 
       & \quad
       +
       Ch^{p/2-1}\int_0^{t} \int_{\lfloor s \rfloor_N}^s \bigg(\mathbb{E}\Big[ \mathbbm{1}_{r\leq \tau_N}\big\|g(Y_r)-b(r)\big\|^p\Big]\bigg)^{1/2}{\rm d}r {\rm d}s+ Ch^p.
    \end{aligned}
\end{equation}
Thanks to the  H{\"o}lder inequality, the Young inequality and condition (c) in Theorem \ref{thm:main_thm}, we treat $B_6$ in the following way:
\begin{equation}\label{eq:B6_estimate_without_sup}
    \begin{aligned}
       B_6&\leq
       C\int_0^{t} \int_{\lfloor s \rfloor_N}^s \mathbb{E}\bigg[\mathbbm{1}_{s\leq \tau_N}\tfrac{|X_r-Y_r|^{p-1}}{\exp (\int_0^{r} {{\eta _{\iota}}} {\rm d}\iota)}\big(1+|X_r|+|Y_r|\big)^{2c_g}
       {\color{red} 
       \big|\int_{\lfloor s \rfloor_N}^s \big\langle f'(Y_{\iota}),b({\iota}){\rm{d}}W_{\iota}\big\rangle\big|
       }
       \bigg] {\rm{d}}r{\rm d}s
       \\
       &\ \ \ +C\int_0^{t} \int_{\lfloor s \rfloor_N}^s \mathbb{E}\bigg[\mathbbm{1}_{s\leq \tau_N}\tfrac{|X_r-Y_r|^{p-3}}{\exp (\int_0^{r} {{\eta _{\iota}}} {\rm d}\iota)}\big\|g(Y_r)-b(r)\big\|^2
       {\color{red} 
       \big|\int_{\lfloor s \rfloor_N}^s \big\langle f'(Y_{\iota}),b({\iota}){\rm{d}}W_{\iota}\big\rangle\big|
       }
       \bigg] {\rm{d}}r{\rm d}s
       \\
        &\leq 
       C\int_0^{t} \sup_{r\in[0,s]}\mathbb{E}\bigg[\tfrac{|X_{r \wedge \tau_N}-Y_{r \wedge \tau_N}|^{p}}{\exp (\int_0^{r \wedge \tau_N} {{\eta _{\iota}{\rm d}\iota)}}}\bigg]{\rm d}s +
       Ch^{\frac{p-2}{2}}\int_0^{t} \int_{\lfloor s \rfloor_N}^s \bigg(\mathbb{E}\Big[ \mathbbm{1}_{r\leq \tau_N}\big\|g(Y_r)-b(r)\big\|^p\Big]\bigg)^{\frac{2}{3}}{\rm d}r {\rm d}s+ Ch^p.
    \end{aligned}
\end{equation}
Gathering \eqref{eq:estimate_SB_temp}, \eqref{eq:B2_estimate_without_sup}, \eqref{eq:B3_estimate_without_sup}, \eqref{eq:B4_estimate_without_sup}, \eqref{eq:B5_estimate_without_sup} and \eqref{eq:B6_estimate_without_sup} yields
\begin{equation}\label{eq:S_B_estimate_without_sup}
    \begin{aligned}
        \mathbb{S}_2 &\leq C\int_0^{t} \sup_{r\in[0,s]}\mathbb{E}\bigg[\tfrac{|X_{r \wedge \tau_N}-Y_{r \wedge \tau_N}|^{p}}{\exp (\int_0^{r \wedge \tau_N} {{\eta _{\iota}{\rm d}\iota)}}}\bigg]{\rm d}s +
       Ch^{\frac{p-2}{2}}\mathbb{E}\bigg[\int_0^{t} \int_{\lfloor s \rfloor_N}^s\mathbbm{1}_{s\leq \tau_N}|f(Y_{\lfloor s \rfloor_N})-a(r)|^p{\rm d}r {\rm d}s\bigg]+ Ch^p
       \\
       &\ \ \ +Ch^{\frac{p-2}{2}}\int_0^{t} \int_{\lfloor s \rfloor_N}^s \bigg(\mathbb{E}\Big[ \mathbbm{1}_{r\leq \tau_N}\big\|g(Y_r)-b(r)\big\|^p\Big]\bigg)^{\frac12}{\rm d}r {\rm d}s+C\mathbb{E}\bigg[\int_0^{t} \mathbbm{1}_{s\leq \tau_N}|f(Y_{\lfloor s \rfloor_N})-a(s)|^p {\rm d}s\bigg].
    \end{aligned}
\end{equation}
Then, combining  \eqref{eq:estimate_S_D_without_sup} with \eqref{eq:S_B_estimate_without_sup} and by Gronwall's inequality we arrive at the first assertion \eqref{eq:perturbation_theory_no_sup}. 
Now it remains to validate \eqref{eq:perturbation_theory_with_sup}.
With regard to $\mathbb{T}_1$, by Lemma \ref{lem:BDG_continuous_inequality}, the inner product inequality, the H{\"o}lder inequality and the elementary inequality, one deduces
        \begin{align*}
            \mathbb{T}_1
            &\leq C_p\left( \int_0^u \sum\limits_{i=1}^{m}\Big\| \mathbbm{1}_{s\leq \tau_N}{\tfrac{{\langle {X_s} - {Y_s},{g^{(i)}({X_s}) -b^{(i)}(s)}\rangle }}{{\exp (\int_0^s {{\eta_r}} {\rm{d}}r)}}}\Big\|^2_{{L^{p/2}(\Omega;\mathbb{R})}}{\rm d}s\right)^{1/2}
            \\
            &\leq C_p \sum\limits_{i=1}^{m}\left( \int_0^u \Big\| \mathbbm{1}_{s\leq \tau_N}{\tfrac{{|{X_s} - {Y_s}| }}{{\exp (\int_0^s {\frac{1}{2}{\eta_r}} {\rm{d}}r)}}}\Big\|^2_{{L^{p}(\Omega;\mathbb{R})}}\cdot \right.\left.\Big\| \mathbbm{1}_{s\leq \tau_N}{\tfrac{{|{g^{(i)}({X_s}) -b^{(i)}(s)}|}}{{\exp (\int_0^s {\frac{1}{2}{\eta_r}} {\rm{d}}r)}}}\Big\|^2_{{L^{p}(\Omega;\mathbb{R})}}{\rm d}s\right)^{1/2}
            \\
            &\leq \sup_{t\in[0,u]}\Big\| {\tfrac{|{{X_{t\wedge\tau_N}} - {Y_{t\wedge\tau_N}}| }}{{\exp (\int_0^{t\wedge\tau_N} {\frac{1}{2}{\eta_r}} {\rm{d}}r)}}}\Big\|_{{L^{p}(\Omega;\mathbb{R})}}\cdot C_p \sum\limits_{i=1}^{m} \left( \int_0^u  \Big\|\mathbbm{1}_{s\leq \tau_N}{\tfrac{{| {g^{(i)}({X_s}) -b^{(i)}(s)} |}}{{\exp (\int_0^s {\frac{1}{2}{\eta_r}} {\rm{d}}r)}}}\Big\|^2_{{L^{p}(\Omega;\mathbb{R})}}{\rm d}s\right)^{1/2}
            \\
            &\leq \tfrac{1}{4}\sup_{t\in[0,u]}\Big\| {\tfrac{|{{X_{t\wedge\tau_N}} - {Y_{t\wedge\tau_N}} |}}{{\exp (\int_0^{t\wedge\tau_N} {\frac{1}{2}{\eta_r}} {\rm{d}}r)}}}\Big\|^2_{{L^{p}(\Omega;\mathbb{R})}} 
            \\
            &\ \ \ \ +C_p \sum\limits_{i=1}^{m} \left( \int_0^u  \Big\|\mathbbm{1}_{s\leq \tau_N}{\tfrac{{| {g^{(i)}({X_s}) -g^{(i)}({Y_s})+ g^{(i)}({Y_s})-b^{(i)}(s)} |}}{{\exp (\int_0^s {\frac{1}{2}{\eta_r}} {\rm{d}}r)}}}\Big\|^2_{{L^{p}(\Omega;\mathbb{R})}}{\rm d}s\right), 
            \stepcounter{equation}\tag{\theequation}
            \label{eq:A_term_01}
        \end{align*}
    where
    \begin{equation}\label{eq:A_term_02}
        \begin{aligned}
            &\int_0^u  \Big\|\mathbbm{1}_{s\leq \tau_N}{\tfrac{{| {g^{(i)}({X_s}) -g^{(i)}({Y_s})+ g^{(i)}({Y_s})-b^{(i)}(s)} |}}{{\exp (\int_0^s { 
            {\color{red}
            \frac{1}{2}{\eta_r}}
            }
            {\rm{d}}r)}}}\Big\|^2_{{L^{p}(\Omega;\mathbb{R})}}{\rm d}s 
            \\
            &\leq 2C\int_0^u  \Big\|\tfrac{{|{{X_{s\wedge\tau_N}} -{Y_{s\wedge\tau_N}}}|}}{{\exp (\int_0^{s\wedge\tau_N} 
            {\color{red}
            \frac{1}{2}{\eta_r}}
             {\rm{d}}r)}}\Big\|^2_{{L^{p}(\Omega;\mathbb{R})}}{\rm d}s
            +
            2\int_0^u\Big\|\mathbbm{1}_{s\leq \tau_N}\tfrac{|g^{(i)}({Y_s})-b^{(i)}(s)|}{\exp (\int_0^s  
            {\color{red}
            \frac{1}{2}{\eta_r}}
            {\rm{d}}r)}\Big\|^2_{{L^{p}(\Omega;\mathbb{R})}}{\rm d}s.
        \end{aligned}
    \end{equation}
    Combining \eqref{eq:A_term_01} with \eqref{eq:A_term_02} yields
 \begin{equation}\label{eq:estimate_T_A_with_sup}
        \begin{aligned}
            \mathbb{T}_1 &\leq \tfrac{1}{4}\sup_{t\in[0,u]}\Big\| {\tfrac{{|{X_{t\wedge\tau_N}} - {Y_{t\wedge\tau_N}}| }}{{\exp (\int_0^{t\wedge\tau_N} {\frac{1}{2}{\eta_r}} {\rm{d}}r)}}}\Big\|^2_{{L^{p}(\Omega;\mathbb{R})}}
            \\
            &\ \ +C\int_0^u  \sup_{s\in[0,t]}\Big\|{\tfrac{{{|{X_{s\wedge\tau_N}} -{Y_{s\wedge\tau_N}}|}}}{{\exp (\int_0^{s\wedge\tau_N} {\frac{1}{2}{\eta_r}} {\rm{d}}r)}}}\Big\|^2_{{L^{p}(\Omega;\mathbb{R})}}{\rm d}t
            +C\int_0^u\Big\|\mathbbm{1}_{s\leq \tau_N}\big({g({Y_s})-b(s)}\big) \Big\|^2_{{L^{p}(\Omega;\mathbb{R}^{d\times m})}}{\rm d}s.
        \end{aligned}
    \end{equation}
 When it comes to $\mathbb{T}_3$, one can easily show
    \begin{equation} \label{eq:estimate_T_D_with_sup}
        \begin{aligned}
            \mathbb{T}_3 
            &\leq C\int_0^u \Big\|\mathbbm{1}_{s\leq \tau_N}\big(g({Y_s}) - b(s)\big)\Big\|^2_{{L^{p}(\Omega;\mathbb{R}^{d\times m})}} {\rm{d}}s.
        \end{aligned}
    \end{equation}
    Note that the term $\mathbb{T}_2$ needs to be treated carefully. Firstly, using the same arguments as  $\mathbb{S}_2$ shows
    \begin{equation}\label{eq:estimate_T_B_temp}
        \begin{aligned}
            \mathbb{T}_2 &\leq \Big\|\sup_{t\in[0,u]}\int_0^{t \wedge \tau_N} {{\Big\langle \tfrac{2({X_s} - {Y_s})}{\exp (\int_0^s {{\eta_r}} {\rm{d}}r)},\int_{\lfloor s \rfloor_N}^s\langle f'(Y_r),a(r)\rangle +\tfrac{1}{2} \operatorname{trace}\big(b(r)^*\operatorname{Hess}_x(f(Y_r))b(r)\big){\rm{d}}r \Big\rangle }} {\rm{d}}s\Big\|_{{L^{\frac{p}{2}}(\Omega;\mathbb{R})}}
            \\
            &\ \ \ +
            \Big\|\sup_{t\in[0,u]}\int_0^{t \wedge \tau_N} {{2\Big\langle \tfrac{{X_s} - {Y_s}}{\exp (\int_0^s {{\eta_r}} {\rm{d}}r)},\int_{\lfloor s \rfloor_N}^s \langle f'(Y_r),b(r){\rm d}W_r\rangle  \Big\rangle }} {\rm{d}}s\Big\|_{{L^{p/2}(\Omega;\mathbb{R})}}
            \\ 
            &\quad+\Big\|\sup_{t\in[0,u]}\int_0^{t \wedge \tau_N} {{2\Big\langle \tfrac{{X_s} - {Y_s}}{\exp (\int_0^s {{\eta_r}} {\rm{d}}r)},f(Y_{\lfloor s \rfloor_N})-a(s) \Big\rangle }} {\rm{d}}s\Big\|_{{L^{p/2}(\Omega;\mathbb{R})}}
            \\
            &\leq C\int_0^u \sup_{s\in[0,t]}\Big\|{\tfrac{{|{X_{s \wedge \tau_N}} - {Y_{s \wedge \tau_N}}|}}{{\exp (\int_0^{s \wedge \tau_N} {\frac{1}{2}{\eta_r}} {\rm{d}}r)}}}\Big\|^2_{{L^{p}(\Omega;\mathbb{R})}} {\rm{d}}t
            +
            \int_0^u \Big\|\mathbbm{1}_{s\leq \tau_N}\big|{f(Y_{\lfloor s \rfloor_N}) - a(s)}\big|\Big\|^2_{{L^{p}(\Omega;\mathbb{R})}} {\rm{d}}s
            \\
            &\quad+\Big\|\sup_{t\in[0,u]}\int_0^{t \wedge \tau_N} {{2\Big\langle \tfrac{{X_s} - {Y_s}}{\exp (\int_0^s {{\eta_r}} {\rm{d}}r)},\int_{\lfloor s \rfloor_N}^s \langle f'(Y_r),b(r){\rm d}W_r\rangle  \Big\rangle }} {\rm{d}}s\Big\|_{{L^{p/2}(\Omega;\mathbb{R})}}+Ch^2.
      \end{aligned}
    \end{equation}
%
%
    To estimate the last but one term for $p\geq 4$, we expand the left item in the inner product by $\rm It\hat{o}$'s formula and $\rm It\hat{o}$'s product rule to acquire
    \begin{equation}
        \begin{aligned} 
            \tfrac{{X_s} - {Y_s}}{\exp (\int_0^s {{\eta_r}} {\rm{d}}r)}
            &=
            \tfrac{{X_{\lfloor s \rfloor_N}} - {Y_{\lfloor s \rfloor_N}}}{\exp (\int_0^{\lfloor s \rfloor_N} {{\eta_r}} {\rm{d}}r)}
            +
             \int_{\lfloor s \rfloor_N}^s \tfrac{f({X_r}) - {a(r)}}{\exp (\int_0^r {{\eta_{\iota}}} {\rm{d}}\iota)}{\rm d}r
            \\
            &\ \ \ +\int_{\lfloor s \rfloor_N}^s \tfrac{g({X_r}) - {b(r)}}{\exp (\int_0^r {{\eta_{\iota}}} {\rm{d}}\iota)}{\rm d}W_r
            +
            \int_{\lfloor s \rfloor_N}^s \tfrac{(X_r-Y_r)(-\eta_r)}{\exp (\int_0^r {{\eta_{\iota}}} {\rm{d}}\iota)}{\rm d}r.
        \end{aligned}
    \end{equation}
 As a consequence,
    \begin{equation}\label{eq:B1_to_B4_with_sup}
        \begin{aligned}
            &\Big\|\sup_{t\in[0,u]}\int_0^{t \wedge \tau_N} {{2\Big\langle \tfrac{{X_s} - {Y_s}}{\exp (\int_0^s {{\eta_r}} {\rm{d}}r)},\int_{\lfloor s \rfloor_N}^s \big\langle f'(Y_r),b(r){\rm{d}}W_r\big\rangle \Big\rangle }}{\rm{d}}s\Big\|_{{L^{p/2}(\Omega;\mathbb{R})}}
            \\
            &\leq \Big\|\sup_{t\in[0,u]}\int_0^{t \wedge \tau_N} {{2\Big\langle \tfrac{{X_{\lfloor s \rfloor_N}} - {Y_{\lfloor s \rfloor_N}}}{\exp (\int_0^{\lfloor s \rfloor_N} {{\eta_r}} {\rm{d}}r)},\int_{\lfloor s \rfloor_N}^s \big\langle f'(Y_r),b(r){\rm{d}}W_r\big\rangle \Big\rangle }}{\rm{d}}s\Big\|_{{L^{p/2}(\Omega;\mathbb{R})}}
            \\
            &\ \ \ + \Big\|\sup_{t\in[0,u]}\int_0^{t \wedge \tau_N} {{2\Big\langle \int_{\lfloor s \rfloor_N}^s \tfrac{f({X_r}) - {a(r)}}{\exp (\int_0^r {{\eta_{\iota}}} {\rm{d}}\iota)}{\rm d}r,\int_{\lfloor s \rfloor_N}^s \big\langle f'(Y_r),b(r){\rm{d}}W_r\big\rangle \Big\rangle }}{\rm{d}}s\Big\|_{{L^{p/2}(\Omega;\mathbb{R})}}
            \\
            &\ \ \ +\Big\|\sup_{t\in[0,u]}\int_0^{t \wedge \tau_N} {{2\Big\langle \int_{\lfloor s \rfloor_N}^s \tfrac{g({X_r}) - {b(r)}}{\exp (\int_0^r {{\eta_{\iota}}} {\rm{d}}\iota)}{\rm d}W_r,\int_{\lfloor s \rfloor_N}^s \big\langle f'(Y_r),b(r){\rm{d}}W_r\big\rangle \Big\rangle }}{\rm{d}}s\Big\|_{{L^{p/2}(\Omega;\mathbb{R})}}
            \\
            &\ \ \ +\Big\|\sup_{t\in[0,u]}\int_0^{t \wedge \tau_N} {{2\Big\langle \int_{\lfloor s \rfloor_N}^s \tfrac{(X_r-Y_r)(-\eta_r)}{\exp (\int_0^r {{\eta_{\iota}}} {\rm{d}}\iota)}{\rm d}r,\int_{\lfloor s \rfloor_N}^s \big\langle f'(Y_r),b(r){\rm{d}}W_r\big\rangle \Big\rangle }}{\rm{d}}s\Big\|_{{L^{p/2}(\Omega;\mathbb{R})}}
            \\
            &=: \widetilde{B}_{1}+\widetilde{B}_{2}+\widetilde{B}_{3}+\widetilde{B}_{4}.
        \end{aligned}
    \end{equation}
    Let us estimate these four items in \eqref{eq:B1_to_B4_with_sup} separately. We first split $\widetilde{B}_{1}$ into two parts:
    \begin{equation}
        \begin{aligned}
            \widetilde{B}_{1} 
            &\leq\Big\| \sup_{t\in[0,u]} \big|\sum\limits_{k=0}^{n_t-1}\int_{t_k}^{t_{k+1}} 2\mathbbm{1}_{s\leq \tau_N}{{\Big\langle \tfrac{{X_{\lfloor s \rfloor_N}} - {Y_{\lfloor s \rfloor_N}}}{\exp (\int_0^{\lfloor s \rfloor_N} {{\eta_r}} {\rm{d}}r)},\int_{\lfloor s \rfloor_N}^s \big\langle f'(Y_r),b(r){\rm{d}}W_r\big\rangle \Big\rangle }}{\rm{d}}s\big|\Big\|_{{L^{p/2}(\Omega;\mathbb{R})}}
            \\
            &\ \ \ +\Big\|\sup_{t\in[0,u]}\int_{\lfloor t \rfloor_N}^t 2\mathbbm{1}_{s\leq \tau_N}{{\Big\langle \tfrac{{X_{\lfloor s \rfloor_N}} - {Y_{\lfloor s \rfloor_N}}}{\exp (\int_0^{\lfloor s \rfloor_N} {{\eta_r}} {\rm{d}}r)},\int_{\lfloor s \rfloor_N}^s \big\langle f'(Y_r),b(r){\rm{d}}W_r\big\rangle \Big\rangle }}{\rm{d}}s\Big\|_{{L^{p/2}(\Omega;\mathbb{R})}}
            \\
            &=:\widetilde{B}_{11}+\widetilde{B}_{12},
        \end{aligned}
    \end{equation}
    where we denote $n_{t}:=\lfloor t \rfloor_N/h$.
    By the condition (c) in Theorem \ref{thm:main_thm}, it follows that
    \[
        \zeta_n
        : =
        \sum\limits_{k=0}^{n-1}\int_{t_k}^{t_{k+1}} 2\mathbbm{1}_{s\leq \tau_N}{{\Big\langle \tfrac{{X_{\lfloor s \rfloor_N}} - {Y_{\lfloor s \rfloor_N}}}{\exp (\int_0^{\lfloor s \rfloor_N} {{\eta_r}} {\rm{d}}r)},\int_{\lfloor s \rfloor_N}^s \big\langle f'(Y_r),b(r){\rm{d}}W_r\big\rangle \Big\rangle }}{\rm{d}}s
    \]
    is a discrete martingale.
    The Doob discrete martingale inequality, Lemma \ref{lem:BDG_discrete_inequality}, H{\"o}lder's inequality and the condition (c) in Theorem \ref{thm:main_thm} imply that 
    \begin{equation}
        \begin{aligned}
            \widetilde{B}_{11} & \leq C_p\Big\| \sum\limits_{k=0}^{n_u-1}\int_{t_k}^{t_{k+1}} 2\mathbbm{1}_{s\leq \tau_N}{{\Big\langle \tfrac{{X_{\lfloor s \rfloor_N}} - {Y_{\lfloor s \rfloor_N}}}{\exp (\int_0^{\lfloor s \rfloor_N} {{\eta_r}} {\rm{d}}r)},\int_{\lfloor s \rfloor_N}^s \big\langle f'(Y_r),b(r){\rm{d}}W_r\big\rangle \Big\rangle }}{\rm{d}}s\Big\|_{{L^{p/2}(\Omega;\mathbb{R})}}
            \\
            &\leq C_p \left(\sum\limits_{k=0}^{n_u-1} \Big\| \int_{t_k}^{t_{k+1}} 2\mathbbm{1}_{s\leq \tau_N}{{\Big\langle \tfrac{{X_{\lfloor s \rfloor_N}} - {Y_{\lfloor s \rfloor_N}}}{\exp (\int_0^{\lfloor s \rfloor_N} {{\eta_r}} {\rm{d}}r)},\int_{\lfloor s \rfloor_N}^s \big\langle f'(Y_r),b(r){\rm{d}}W_r\big\rangle \Big\rangle }}{\rm{d}}s\Big\|^2_{{L^{p/2}(\Omega;\mathbb{R})}}\right)^{1/2}
            \\
            &\leq C_p \left(h \int_{0}^{\lfloor u \rfloor_N} \Big\|
            \tfrac{|X_{\lfloor s \rfloor_N \wedge \tau_N}-Y_{\lfloor s \rfloor_N \wedge \tau_N}|}{{\exp (\int_0^{\lfloor s \rfloor_N \wedge \tau_N} {\frac{1}{2}{\eta _{r}}} {\rm{d}}r)}}
            \Big\|^2_{{L^{p}(\Omega;\mathbb{R})}} \Big\|\big|\int_{\lfloor s \rfloor_N}^s \big\langle f'(Y_r),b(r){\rm{d}}W_r\big\rangle\big |\Big\|^2_{{L^{p}(\Omega;\mathbb{R})}}{\rm{d}}s\right)^{1/2}
            \\
            &\leq \sup_{s\in [0,u]} \Big\|{{
            \tfrac{|X_{{s \wedge \tau_N}}-Y_{{s \wedge \tau_N}}|}{{\exp (\int_0^{{s \wedge \tau_N}} {\frac{1}{2}{\eta _{r}}} {\rm{d}}r)}}
            \Big\|_{{L^{p}(\Omega;\mathbb{R})}}}} C_p \left(h \int_{0}^{\lfloor u \rfloor_N}
           \Big\|\big|\int_{\lfloor s \rfloor_N}^s \big\langle f'(Y_r),b(r){\rm{d}}W_r\big\rangle\big |\Big\|^2_{{L^{p}(\Omega;\mathbb{R})}}{\rm{d}}s\right)^{1/2}
            \\
            &\leq \tfrac{1}{8}\sup_{s\in [0,u]} \Big\|{{
            \tfrac{|X_{{s \wedge \tau_N}}-Y_{{s \wedge \tau_N}}|}{{\exp (\int_0^{{s \wedge \tau_N}} {\frac{1}{2}{\eta _{r}}} {\rm{d}}r)}}
            \Big\|^2_{{L^{p}(\Omega;\mathbb{R})}}}}+Ch^2.
        \end{aligned}
    \end{equation}
    With the help of H{\"o}lder's inequality, an elementary inequality and the condition (c) in Theorem \ref{thm:main_thm}, 
    for $p\geq4$ one can estimate $\widetilde{B}_{12}$ as
    \begin{equation}
        \begin{aligned}
            \widetilde{B}_{12}&\leq
            \bigg(\mathbb{E}\Big[\sup_{t\in[0,u]}\Big|\int_{\lfloor t \rfloor_N}^t {{2\mathbbm{1}_{s\leq \tau_N}\Big\langle \tfrac{{X_{\lfloor s \rfloor_N}} - {Y_{\lfloor s \rfloor_N}}}{\exp (\int_0^{\lfloor s \rfloor_N} {{\eta_r}} {\rm{d}}r)},\int_{\lfloor s \rfloor_N}^s \big\langle f'(Y_r),b(r){\rm{d}}W_r\big\rangle \Big\rangle }}{\rm{d}}s\Big|^{p/2}\Big]\bigg)^{2/p}
            \\
            &\leq Ch^{1-2/p}\bigg(\mathbb{E}\Big[\int_0^u \mathbbm{1}_{s\leq \tau_N}\Big|{\Big\langle \tfrac{{X_{\lfloor s \rfloor_N}} - {Y_{\lfloor s \rfloor_N}}}{\exp (\int_0^{\lfloor s \rfloor_N} {{\eta_r}} {\rm{d}}r)},\int_{\lfloor s \rfloor_N}^s \big\langle f'(Y_r),b(r){\rm{d}}W_r\big\rangle \Big\rangle }\Big|^{p/2}{\rm{d}}s\Big]\bigg)^{2/p}
            \\
            &\leq Ch^{1-2/p}\bigg(\int_0^u \Big\|\tfrac{|X_{\lfloor s \rfloor_N \wedge \tau_N}-Y_{\lfloor s \rfloor_N \wedge \tau_N}|}{\exp (\int_0^{\lfloor s \rfloor_N \wedge \tau_N} {\frac{1}{2}{\eta _{r}}} {\rm{d}}r)}\Big\|^{p/2}_{L^p(\Omega;\mathbb{R})}\Big\|\big|\int_{\lfloor s \rfloor_N}^s \big\langle f'(Y_r),b(r){\rm{d}}W_r\big\rangle\big|\Big\|^{p/2}_{L^p(\Omega;\mathbb{R})}{\rm{d}}s\bigg)^{2/p}
            \\
            &\leq \tfrac{1}{8}\sup_{s\in[0,u]}\Big\|\tfrac{|X_{s \wedge \tau_N}-Y_{s \wedge \tau_N}|}{\exp (\int_0^{s \wedge \tau_N} {\frac{1}{2}{\eta _{r}}} {\rm{d}}r)}\Big\|^{2}_{L^p(\Omega;\mathbb{R})}+Ch^{2-4/p}\left(\int_0^u
            \Big\|\big|\int_{\lfloor s \rfloor_N}^s \big\langle f'(Y_r),b(r){\rm{d}}W_r\big\rangle\big|\Big\|^{p/2}_{L^{p}(\Omega;\mathbb{R})}{\rm{d}}s\right)^{4/p}
            \\
            &\leq \tfrac{1}{8}\sup_{s\in[0,u]}\Big\|\tfrac{|X_{s \wedge \tau_N}-Y_{s \wedge \tau_N}|}{\exp (\int_0^{s \wedge \tau_N} {\frac{1}{2}{\eta _{r}}} {\rm{d}}r)}\Big\|^{2}_{L^p(\Omega;\mathbb{R})}+Ch^2.
        \end{aligned}
    \end{equation}
    Hence one concludes that 
    \begin{equation}\label{eq:B1_estimate_with_sup}
        \begin{aligned}
             \widetilde{B}_{1}
            &\leq \tfrac{1}{4}\sup_{s\in [0,u]} \Big\|{{
            \tfrac{|X_{s \wedge \tau_N}-Y_{s \wedge \tau_N}|}{{\exp (\int_0^{s \wedge \tau_N} {\frac{1}{2}{\eta _{r}}} {\rm{d}}r)}}
            \Big\|^2_{{L^{p}(\Omega;\mathbb{R})}}}}+Ch^2.
        \end{aligned}
    \end{equation}
Similar to the estimate of ${B}_{2}$, one treat 
$\widetilde{B}_{2}$ as follows:
    \begin{equation}\label{eq:B2_estimate_with_sup}
        \begin{aligned}
            \widetilde{B}_{2}
            &\leq C\int_0^{u} {{\int_{\lfloor s \rfloor_N}^s \Big\|\mathbbm{1}_{s\leq \tau_N}\tfrac{|{X_r} - {Y_r}|(1+|X_r|+|Y_r|)^{c_f}}{\exp (\int_0^r {{\eta_{\iota}}} {\rm{d}}\iota)}
           {\color{red} 
           \big|\int_{\lfloor s \rfloor_N}^s \big\langle f'(Y_{\iota}),b({\iota}){\rm{d}}W_{\iota}\big\rangle\big|
           }
            \Big\|_{{L^{p/2}(\Omega;\mathbb{R})}}{\rm d}r }}{\rm{d}}s
            \\
            &\ \ \ +C\int_0^{u} {{\int_{\lfloor s \rfloor_N}^s \Big\|\mathbbm{1}_{s\leq \tau_N}\tfrac{|f({Y_r})-f({Y_{\lfloor r \rfloor_N}})|}{\exp (\int_0^r {{\eta_{\iota}}} {\rm{d}}\iota)}
            {\color{red} 
            \big|\int_{\lfloor s \rfloor_N}^s \big\langle f'(Y_{\iota}),b({\iota}){\rm{d}}W_{\iota}\big\rangle\big|
            }
            \Big\|_{{L^{p/2}(\Omega;\mathbb{R})}}{\rm d}r }}{\rm{d}}s
            \\
            &\ \ \ +C\int_0^{u} {{\int_{\lfloor s \rfloor_N}^s \Big\|\mathbbm{1}_{s\leq \tau_N}\tfrac{|f({Y_{\lfloor r \rfloor_N}}) - a(r)|}{\exp (\int_0^r {{\eta_{\iota}}} {\rm{d}}\iota)}
            {\color{red} 
            \big|\int_{\lfloor s \rfloor_N}^s \big\langle f'(Y_{\iota}),b({\iota}){\rm{d}}W_{\iota}\big\rangle\big|
            }
            \Big\|_{{L^{p/2}(\Omega;\mathbb{R})}} {\rm d}r}}{\rm{d}}s
            \\
            &\leq C\int_0^{u}\sup_{r\in[0,s]}\Big\|\tfrac{|{X_{r \wedge \tau_N}} - {Y_{r \wedge \tau_N}}|}{\exp (\int_0^{r \wedge \tau_N} {\frac{1}{2}{\eta_{\iota}}} {\rm{d}}\iota)}\Big\|^2_{{L^{p}(\Omega;\mathbb{R})}}{\rm{d}}s +Ch^2
            \\
            &\ \ \ +Ch^{1/2}\int_0^{u} {{\int_{\lfloor s \rfloor_N}^s \Big\|\mathbbm{1}_{s\leq \tau_N}|f({Y_{\lfloor r \rfloor_N}}) - a(r)|\Big\|_{{L^{p}(\Omega;\mathbb{R})}} {\rm d}r}}{\rm{d}}s.
        \end{aligned}
    \end{equation}
With regard to $\widetilde{B}_{3}$, we employ H{\"o}lder's inequality, Young's inequality and Lemma \ref{lem:BDG_continuous_inequality} to derive
        \begin{align*}
            \widetilde{B}_{3}
            &\leq C\int_0^{u}  {{\Big\| \int_{\lfloor s \rfloor_N}^s \mathbbm{1}_{s\leq \tau_N}\tfrac{g({X_r}) -g({Y_r})}{\exp (\int_0^r {{\eta_{\iota}}} {\rm{d}}\iota)}{\rm d}W_r\Big\|_{{L^{p}(\Omega;\mathbb{R}^d)}} \Big\|\int_{\lfloor s \rfloor_N}^s \big\langle f'(Y_r),b(r){\rm{d}}W_r\big\rangle \Big\rangle }}\Big\|_{{L^{p}(\Omega;\mathbb{R}^d)}}{\rm{d}}s
            \\
            &\ \ \ + C\int_0^{u}  {{\Big\| \int_{\lfloor s \rfloor_N}^s \mathbbm{1}_{s\leq \tau_N}\tfrac{g({Y_r}) - {b(r)}}{\exp (\int_0^r {{\eta_{\iota}}} {\rm{d}}\iota)}{\rm d}W_r\Big\|_{{L^{p}(\Omega;\mathbb{R}^d)}} \Big\|\int_{\lfloor s \rfloor_N}^s \big\langle f'(Y_r),b(r){\rm{d}}W_r\big\rangle \Big\rangle }}\Big\|_{{L^{p}(\Omega;\mathbb{R}^d)}}{\rm{d}}s
            \\
            &\leq C\int_0^{u}  {{\Big\| \int_{\lfloor s \rfloor_N}^s \mathbbm{1}_{s\leq \tau_N}\tfrac{g({X_r}) -g({Y_r})}{\exp (\int_0^r {{\eta_{\iota}}} {\rm{d}}\iota)}{\rm d}W_r\Big\|^{4/3}_{{L^{p}(\Omega;\mathbb{R}^d)}} }}{\rm{d}}s+Ch^2
            \\
            &\ \ \ + Ch^{1/2}\int_0^{u}  {{\Big\| \int_{\lfloor s \rfloor_N}^s \mathbbm{1}_{s\leq \tau_N}\tfrac{g({Y_r}) - {b(r)}}{\exp (\int_0^r {{\eta_{\iota}}} {\rm{d}}\iota)}{\rm d}W_r\Big\|_{{L^{p}(\Omega;\mathbb{R}^d)}} }}{\rm{d}}s
            \\
            &\leq C\int_0^{u}  {{ \Big(\int_{\lfloor s \rfloor_N}^s \Big\|\mathbbm{1}_{s\leq \tau_N}\tfrac{|{X_r} -{Y_r}|}{\exp (\int_0^r {{\eta_{\iota}}} {\rm{d}}\iota)}\Big\|^{2}_{{L^{p}(\Omega;\mathbb{R})}} {\rm d}r\Big)^{2/3} }}{\rm{d}}s+Ch^2
            \\
            &\ \ \ + Ch^{1/2}\int_0^{u} \Big( {{\int_{\lfloor s \rfloor_N}^s  \Big\|\mathbbm{1}_{s\leq \tau_N}\big(g({Y_r}) - {b(r)}\big)\Big\|^2_{{L^{p}(\Omega;\mathbb{R}^{d\times m})}}{\rm d}r }}\Big)^{1/2}{\rm{d}}s
            \\
            &\leq C\int_0^{u}\sup_{r\in[0,s]}\Big\|\tfrac{|{X_{r \wedge \tau_N}} - {Y_{r \wedge \tau_N}}|}{\exp (\int_0^{r \wedge \tau_N} {\frac{1}{2}{\eta_{\iota}}} {\rm{d}}\iota)}\Big\|^2_{{L^{p}(\Omega;\mathbb{R})}}{\rm{d}}s +Ch^2
            \\
            &\ \ \ + Ch^{1/2}\int_0^{u} \Big( {{\int_{\lfloor s \rfloor_N}^s  \Big\|\mathbbm{1}_{s\leq \tau_N}\big(g({Y_r}) - {b(r)}\big)\Big\|^2_{{L^{p}(\Omega;\mathbb{R}^{d\times m})}}{\rm d}r }}\Big)^{1/2}{\rm{d}}s. 
            \stepcounter{equation}\tag{\theequation}
            \label{eq:B3_estimate_with_sup}
        \end{align*}
    Similar to the estimate of $B_3$, we bound $\widetilde{B}_{4}$ in the following way:
    \begin{equation}\label{eq:B4_estimate_with_sup}
        \begin{aligned}
            \widetilde{B}_{4}&\leq C\Big\|\int_0^u {{ \int_{\lfloor s \rfloor_N}^s \mathbbm{1}_{s\leq \tau_N}\tfrac{|X_r-Y_r|\eta_r}{\exp (\int_0^r {{\eta_{\iota}}} {\rm{d}}\iota)}
           {\color{red} 
           \big|\int_{\lfloor s \rfloor_N}^s \big\langle f'(Y_{\iota}),b({\iota}){\rm{d}}W_{\iota}\big\rangle\big|
           }
            {\rm d}r }}{\rm{d}}s\Big\|_{{L^{p/2}(\Omega;\mathbb{R})}}
            \\
            &\leq C\int_0^u {{ \int_{\lfloor s \rfloor_N}^s \Big\|\mathbbm{1}_{s\leq \tau_N}\tfrac{|X_r-Y_r|}{\exp (\int_0^r {{\eta_{\iota}}} {\rm{d}}\iota)}\Big\|_{{L^{p}(\Omega;\mathbb{R})}}\Big\|}{\eta_{r}
           {\color{red} 
           \big|\int_{\lfloor s \rfloor_N}^s \big\langle f'(Y_{\iota}),b({\iota}){\rm{d}}W_{\iota}\big\rangle\big|
           }
            \Big\|_{{L^{p}(\Omega;\mathbb{R})}}{\rm d}r }}{\rm{d}}s
            \\
            &\leq C\int_0^{u}\sup_{r\in[0,s]}\Big\|\tfrac{|{X_{s \wedge \tau_N}} - {Y_{s \wedge \tau_N}}|}{\exp (\int_0^{s \wedge \tau_N}{\frac{1}{2}{\eta_{\iota}}} {\rm{d}}\iota)}\Big\|^2_{{L^{p}(\Omega;\mathbb{R})}}{\rm{d}}s +Ch^2.
        \end{aligned}
    \end{equation}
    Putting 
    \eqref{eq:estimate_T_B_temp},
    \eqref{eq:B1_estimate_with_sup}, \eqref{eq:B2_estimate_with_sup}, \eqref{eq:B3_estimate_with_sup} and \eqref{eq:B4_estimate_with_sup} together yields
\begin{equation}\label{eq:estimate_T_B_with_sup}
        \begin{aligned}
        \mathbb{T}_2
        &\leq C\int_0^u \sup_{s\in[0,t]}\Big\|{\tfrac{{|{X_{s \wedge \tau_N}} - {Y_{s \wedge \tau_N}}|}}{{\exp (\int_0^{s \wedge \tau_N} {\frac{1}{2}{\eta_r}} {\rm{d}}r)}}}\Big\|^2_{{L^{p}(\Omega;\mathbb{R})}} {\rm{d}}t+\int_0^u \Big\|\mathbbm{1}_{s\leq \tau_N}\big|{f(Y_{\lfloor s \rfloor_N}) - a(s)}\big|\Big\|^2_{{L^{p}(\Omega;\mathbb{R})}} {\rm{d}}s
        \\
        &\quad+\tfrac{1}{4}\sup_{s\in [0,u]} \Big\|{{
            \tfrac{|X_{s \wedge \tau_N}-Y_{s \wedge \tau_N}|}{{\exp (\int_0^{s \wedge \tau_N} {\frac{1}{2}{\eta _{r}}} {\rm{d}}r)}}
            \Big\|^2_{{L^{p}(\Omega;\mathbb{R})}}}}+Ch^{1/2}\int_0^{u} {{\int_{\lfloor s \rfloor_N}^s \Big\|\mathbbm{1}_{r\leq \tau_N}|f({Y_{\lfloor r \rfloor_N}}) - a(r)|\Big\|_{{L^{p}(\Omega;\mathbb{R})}} {\rm d}r}}{\rm{d}}s
        \\
        &\quad +
        Ch^{1/2}\int_0^{u} \Big( {{\int_{\lfloor s \rfloor_N}^s  \Big\|\mathbbm{1}_{r\leq \tau_N}\big(g({Y_r}) - {b(r)}\big)\Big\|^2_{{L^{p}(\Omega;\mathbb{R}^{d\times m})}}{\rm d}r }}\Big)^{1/2}{\rm{d}}s+Ch^2.
        \end{aligned}
    \end{equation}
    Then the proof is thus completed by combining \eqref{eq:estimate_T_A_with_sup},  \eqref{eq:estimate_T_D_with_sup}, \eqref{eq:estimate_T_B_with_sup} and Gronwall's inequality.\qed

In the rest of this article, we concentrate on applications of
the previously obtained perturbation estimates to 
identify the order-one strong convergence of numerical methods for SDEs with non-globally monotone coefficients. 

\section{Order-one pathwise uniformly strong convergence of the SITEM scheme 
}
\label{sect:application1}

In order to numerically solve SDEs \eqref{eq:typical_sde} on a uniform grid $\{t_n = n h\}_{0 \leq n \leq N}$ with stepsize $h = \tfrac{T}{N}$,
a class of stopped increment-tamed EM (SITEM) method was proposed in \cite{hutzenthaler2018exponential}:
    \begin{equation}\label{eq:stop_tamed_method}
        Y_{t}=Y_{t_{n}}+\mathbbm{1}_{|Y_{t_{n}}|<\exp \left(\left|\ln (h)\right|^{1 / 2}\right)}\left[\tfrac{f\left(Y_{t_{n}}\right)\left(t-t_{n}\right)+g\left(Y_{t_{n}}\right)\left(W_{t}-W_{t_{n}}\right)}{1+\left|f\left(Y_{t_{n}}\right)\left(t-t_{n}\right)+g\left(Y_{t_{n}}\right)\left(W_{t}-W_{t_{n}}\right)\right|^{\delta}}\right], 
        \:
        Y_0=X_0,
        \
        \delta \geq 2,
        \
        t \in [t_n, t_{n+1}],
    \end{equation}
which was shown to inherit exponential integrability properties of original SDEs.
Combining the  perturbation theory obtained in \cite{Hutzenthaler2020}  with exponential integrability properties of both numerical solution and exact solution,
    the authors of \cite{Hutzenthaler2020} successfully 
    identified the order $\tfrac12$ strong convergence of 
    the SITEM method. An interesting question arises as to whether higher 
    convergence rate than order $\tfrac12$ can be obtained
    when high-order (e.g., Milstein-type) schemes 
    are used.  
    {\color{black}This is also expected  
    by  \cite{Hutzenthaler2020} 
    (see Remark 3.1 therein).}
    Unfortunately, following \cite[Theorem 1.2]{Hutzenthaler2020}, the convergence rates
    of any schemes would not exceed order $\tfrac12$,
which is nothing but the order of the H\"{o}lder regularity of the approximation process.
{\color{black} In the present section, we aim to 
 fill this gap and
reveal order-one strong convergence of the SITEM method
for
some particular SDEs with non-globally monotone coefficients,}
for which the Euler type method coincides with
the Milstein method and thus the order-one 
convergence is expected.
    To begin with, we define a stopping time $\tau^e_{N}:\Omega \rightarrow \{t_0,t_1,...,t_N\}$ as
    \begin{equation}
        \tau^e_{N}:=\inf \left\{\{T\} \cup \{ t\in {t_0,...,t_N}:|Y_{t}|\geq \exp (\left|\ln (h)\right|^{1 / 2})\} \right\}.
    \end{equation}
    Equipped with the stopping time, we can introduce the continuous version of 
    \eqref{eq:stop_tamed_method} as
    \begin{equation}\label{eq:stop_tamed_sde}
        Y_{t}=X_0+\int_0^t\mathbbm{1}_{s<\tau^e_{N}} a(s){\rm d}s+\int_0^t\mathbbm{1}_{s<\tau^e_{N}} b(s){\rm d}W_s.
    \end{equation}
    Here, for $s\in[t_k,t_{k+1})$, $a(s)$ and $b(s)$ are given by
    \begin{equation}
        a(s):=\psi^{[1]}(Z_s)f(Y_{t_k})+\tfrac{1}{2} \sum_{j=1}^{m}\psi^{[2]}(Z_s)\Big(g(Y_{t_k})e_j,g(Y_{t_k})e_j\Big), 
        \quad
        b(s):=\psi^{[1]}(Z_s)g(Y_{t_k}),
    \end{equation}
    where  $e_1=(1,...,0)^*,...,e_m=(0,...,1)^*$ are the Euclidean orthonormal basis of $\mathbb{R}^m$,
    \begin{equation}
        Z_s := f(Y_{t_k})(s-t_k)+g(Y_{t_k})(W_s-W_{t_k})
    \end{equation}
    {\color{black} and for  a fixed $\delta\geq2,$}
    \begin{equation} \label{eq:defn-psi}
        \psi(x):=x(1+|x|^\delta)^{-1}, x\in \mathbb{R}^d.
    \end{equation}
    By \cite[Theorem 2.9]{hutzenthaler2018exponential}
    and in the notation of \eqref{eq:F-derivatives},
    for any $z,u \in \mathbb{R}^d$ we have
    \begin{equation} \label{eq:psi_1_deri}
        \psi^{[1]}(z) u=
        \left\{\begin{array}{ll}u 
        & : z=0 ,
        \\ 
        \tfrac{u}{1+|z|^{\delta}}-\tfrac{\delta z|z|^{(\delta-2)}\langle z, u\rangle}{(1+|z|^{\delta})^{2}} & : z \neq 0,
        \end{array}\right.
    \end{equation} 
    and
    \begin{equation}\label{eq:psi_2_deri}
        \psi^{[2]}(z)(u, u)=\left\{\begin{array}{ll}0 & : z=0,
       \\ 
       \tfrac{2 \delta^{2}|z|^{2(\delta-2)} z|\langle z, u\rangle|^{2}}{(1+|z|^{\delta})^{3}}-\tfrac{\delta|z|^{(\delta-2)}\left[2 u\langle z, u\rangle+z|u|^{2}\right]}{(1+|z|^\delta)^{2}} & 
       \\ 
       -\tfrac{\delta(\delta-2)|z|^{(\delta-4)} z|\langle z, u\rangle|^{2}}{(1+|z|^{\delta})^{2}} & : z \neq 0.\end{array}\right.
    \end{equation}
Moreover, one can show the following properties of $\psi^{[1]},
\psi^{[2]}$, which are needed in the error analysis.
    \begin{lemma}\label{lem:tame_term_estimate}
        Let  $\delta\geq2$ and let $\psi$ be defined by \eqref{eq:defn-psi}. 
        Then for all $x\in \mathbb{R}^d,$
        \begin{equation}
            \begin{aligned}
            \|\psi^{[1]}(x)\|_{L(\mathbb{R}^d,\mathbb{R}^d)}
            &
            \leq
            1
            + \tfrac{\delta}{4}, \ \|\psi^{[1]}(x)-I\|_{L(\mathbb{R}^d,\mathbb{R}^d)}\leq \left[(1+ \tfrac{\delta}{4})\wedge (\delta+1)|x|^{\delta}\right],
            \\
            & \sup_{u\in\mathbb{R}^d,|u|\leq1} |\psi^{[2]}(x)(u, u)|\leq \left[(3{\delta}^2+{\delta})\wedge (3{\delta}^2+{\delta})|x|^{{\delta}-1}\right].
            \end{aligned}
        \end{equation} 
    \end{lemma}
    \indent \textbf{Proof}. By \eqref{eq:psi_1_deri} and \eqref{eq:psi_2_deri}, it is clear that
    \begin{equation}
        \begin{aligned}
            &\|\psi^{[1]}(x)\|_{L(\mathbb{R}^d,\mathbb{R}^d)}\leq 1 \vee \Big(\tfrac{1}{1+|x|^{\delta}} + \tfrac{\delta|x|^{\delta}}{(1+|x|^{\delta})^2}\Big)
            \leq 1+\tfrac{\delta}{4},
            \\
            & \|\psi^{[1]}(x)-I\|_{L(\mathbb{R}^d,\mathbb{R}^d)}
            \leq \Big(\tfrac{|x|^{\delta}}{1+|x|^{\delta}} +
            \tfrac{\delta|x|^{\delta}}{(1+|x|^{\delta})^2}\Big) \leq \left[(1+ \tfrac{\delta}{4})\wedge \big(\delta+1\big)|x|^{\delta}\right],
            \\
            & \sup_{u\in\mathbb{R}^d,|u|\leq1} |\psi^{[2]}(x)(u, u)|
            \leq 
            \tfrac{2\delta^2|x|^{2\delta-1}}{(1+|x|^\delta)^3}
            +
            \tfrac{(\delta^2+\delta)|x|^{\delta-1}}{(1+|x|^\delta)^2}
            \leq 
            \left[
            (3\delta^2+\delta)\wedge (3\delta^2+\delta)|x|^{\delta-1}
            \right].
        \end{aligned}
    \end{equation}
    
    Now we are ready to state the main convergence result of this section.
    
%
\begin{theorem}\label{thm:stopped_tamed_method}
        Let 
        $f:\mathbb R^d\rightarrow \mathbb R^d,g:\mathbb R^d\rightarrow \mathbb R^{d \times m} $ be measurable functions and let $f\in C^2{(\mathbb R^d,\mathbb R^d})$. 
        Let $f\in \mathcal{C}^1_{\mathcal{P}}(\mathbb{R}^d,\mathbb{R}^{d})$ 
        and let $g\in C^2(\mathbb{R}^d,\mathbb{R}^{d\times m})$ be Lipschitz satisfying, {\color{blue} for all $k_1,k_2 \in \{1,...,d\}, j_1, j_2
        \in \{1,...,m\}$},
        \begin{equation}
        \tfrac{\partial g^{(k_2,j_2)}}{\partial x_{k_1}} g^{(k_1,j_1)} = 0.
        \end{equation}
        Let 
        $U_0 \in \mathcal{C}^3_{\mathcal{D}}\big(\mathbb{R}^d,[0,\infty)\big)$ and $U_1 \in \mathcal{C}^1_{\mathcal{P}}\big(\mathbb{R}^d,[0,\infty)\big)$. 
        Let a class of  SITEM methods be defined by \eqref{eq:stop_tamed_sde} with $\delta\geq 3$ and 
        let $c,v,T\in(0,\infty),q,q_1,q_2\in (0,\infty],\alpha \in[0,\infty), p\geq 4 $.
        For all $x,y\in \mathbb R^d$, assume additionally that
        \begin{enumerate}[{\rm(1)}]
            \item 
                there exist constants $L,\kappa\geq0$ such that for any $i=1,...,d,j=1,...,m$,
                $$
                    \|\operatorname{Hess}_x(f^{(i)}(x))\|
                    \ {\textstyle \bigvee}\ 
                     \|\operatorname{Hess}_x(g^{(ij)}(x))\|
                    \leq 
                    L(1+|x|)^{\kappa};
                $$
            \item
                $|x|^{1/c} \leq c(1+U_0(x))$\ and\ $\mathbb{E}\left[e^{U_0(X_0)}\right]<\infty;$
            \item
                $(\mathcal{A}_{f,g} U_0 )(x)+\tfrac{1}{2}|g(x)^{*}(\nabla U_0(x))|^2+U_1(x) \leq c+\alpha U_0(x);$
            \item
                $
                \langle x-y,f(x)-f(y)\rangle \leq \left[c + \tfrac{U_0(x)+U_0(y)}{2q_1Te^{\alpha T}}+\tfrac{U_1(x)+U_1(y)}{2q_2e^{\alpha T}}\right]|x-y|^2. 
                $
        \end{enumerate}
        Then for $\tfrac{1}{q}=\tfrac{1}{q_1}+\tfrac{1}{q_2},\tfrac{1}{v}=\tfrac{1}{p}+\tfrac{1}{q}$ the approximation \eqref{eq:stop_tamed_sde} 
        used to solve \eqref{eq:typical_sde} admits 
        \begin{equation}
    \label{eq:thm_stopped_tamed_method_conclu001}
            \Big\| \sup_{t\in[0,T]}|X_t-Y_t| \Big\|_{L^{v}(\Omega;\mathbb{R})}
            \leq 
            C h, \ h \rightarrow 0.
        \end{equation}  
If the condition $(4)$ in Theorem {\rm \ref{thm:stopped_tamed_method}} is replaced by 
the following one: 
\begin{enumerate}[{\rm(4')}]
            \item 
for any $\eta>0$, there exists a constant $K_{\eta}$ such that  \begin{equation}\label{eq:strengthed_condition_item_stopped_tame}
            \langle x-y,f(x)-f(y)\rangle \leq \big[K_{\eta}+\eta\big(U_0(x)+U_0(y)+U_1(x)+U_1(y)\big)\big]|x-y|^2,
        \end{equation}
\end{enumerate}
        then  for any $v>0$ we have
\begin{equation}\label{eq:strengthed_conclusion_stopped_tame}
            \Big\| \sup_{t\in[0,T]}|X_t-Y_t| \Big \|_{L^{v}(\Omega;\mathbb{R})}\leq Ch
            , \ h \rightarrow 0.
        \end{equation}  
    \end{theorem}
 {\color{red}
 Before we come to the proof, let us give some comments on functions $U_0, U_1$ and some parameters used in the above theorem. We first mention that, similar conditions have been used in \cite{Hutzenthaler2020}. Here the non-negative function $U_0$ plays a role of Lyapunov function for (stochastic) differential equations (see conditions (2), (3)).
 From conditions (4) or (4'), one observes that $U_0$ is also used to control the growth of the derivative of the drift $f$. For some models such as the Brownian dynamics, stochastic van der Pol oscillator and stochastic Duffing-van der Pol oscillator, $U_0$ alone is, however, not able to control the growth and an additional non-negative function $U_1$ is introduced in conditions (3), (4). As shown later, different models require different choices of functions $U_0, U_1$ such that conditions (2), (3), (4) or (4') in Theorem \ref{thm:stopped_tamed_method} are all satisfied.
 In view of Lemma \ref{lem:tame_term_estimate}, we require the method parameter $\delta \geq 3$ to guarantee the convergence order $ 1 \wedge \tfrac{\delta-1}{2} \geq 1$. In addition, the parameters $q,q_1,q_2$ and $v,p,q$ are two sets of conjugate numbers for the use of H\"{o}lder's inequality. Now we start the proof.
    }
    
    \indent \textbf{Proof of Theorem \ref{thm:stopped_tamed_method}}. The proof relies on the use of 
    Theorem \ref{thm:main_thm} and in what follows
    we check all the conditions there. Firstly, let $\tau_N=\tau^e_N $ and it is obvious that for $s\in[0,T]$
    \begin{equation} 
        \{s\leq \tau^e_{N}\}=\left\{ 
        \begin{array}{l}
        \{\lfloor s \rfloor_N \leq \tau^e_{N}\},\lfloor s \rfloor_N=s;\\
        \{\lfloor s \rfloor_N <\tau^e_{N}\},\lfloor s \rfloor_N<s,
        \end{array} \right.
    \end{equation}
    which implies $ \{s\leq \tau^e_{N}\}\in \mathcal{F}_{\lfloor s \rfloor_N}$. 
    {\color{red}
    By the virtue of the condition $(4)$, H{\"o}lder's inequality, Jensen's  inequality (see \cite[Lemma 2.22]{cox2024local} and the fact that $U_0(x),U_1(x)\geq 0$ one derives that 
    }
    \begin{equation}\label{eq:term(4)_to_moment_estimate}
        \begin{aligned}
        &\bigg\|{\exp \bigg( \int_0^{\tau^e_{N}} {\left[\tfrac{\left\langle X_{s}-Y_{s}, f\left(X_{s}\right)-f\left(Y_{s}\right)\right\rangle+\frac{1+\varepsilon}{2}\left\|g\left(X_{s}\right)-g\left(Y_{s}\right)\right\|^{2}}{\left|X_{s}-Y_{s}\right|^{2}}\right]^{+}} {\rm{d}}s\bigg)}\bigg\|_{L^{q}(\Omega;\mathbb{R})}
        \\
        &\leq C_{\varepsilon,\alpha,c,q_1,q_2,T}\bigg\|\exp \bigg( \int_0^{\tau^e_{N}} \left[\tfrac{U_0(X_s)+U_0(Y_s)}{2q_1Te^{\alpha T }}+\tfrac{U_1(X_s)+U_1(Y_s)}{2q_2e^{\alpha T }}\right] {\rm{d}}s\bigg)\bigg\|_{L^{q}(\Omega;\mathbb{R})}
        \\
        &\leq C\Big(\mathbb{E}\Big[\exp\Big(\int_0^{\tau^e_{N}}\tfrac{U_0(X_s)+U_0(Y_s)}{2Te^{\alpha T }}{\rm d}s\Big)\Big]\Big)^{1/q_1}\cdot\Big(\mathbb{E}\Big[\exp\Big(\int_0^{\tau^e_{N}}\tfrac{U_1(X_s)+U_1(Y_s)}{2e^{\alpha T }}{\rm d}s\Big)\Big]\Big)^{1/q_2}
        \\
        &\leq C\sup_{s\in[0,T]} \Big(\mathbb{E}\Big[\exp\Big(\tfrac{U_0(X_s)}{e^{\alpha s }}\Big)\Big]\Big)^{\tfrac{1}{2q_1}}
        \cdot 
        \sup_{s\in[0,T]}\Big(\mathbb{E}\Big[\exp\Big(\tfrac{U_0(Y_s)}{e^{\alpha s }}\Big)\Big]\Big)^{\tfrac{1}{2q_1}}
        \\
        &\ \ \ \cdot \sup_{s\in[0,T]}\Big(\mathbb{E}\Big[\exp\Big(\int_0^{s \wedge\tau^e_{N}}\tfrac{U_1(X_u)}{e^{\alpha u }}{\rm d}u\Big)\Big]\Big)^{\tfrac{1}{2q_2}}\cdot \sup_{s\in[0,T]}
        \Big(\mathbb{E}\Big[\exp\Big(\int_0^{s\wedge\tau^e_{N}}\tfrac{U_1(Y_u)}{e^{\alpha u }}{\rm d}u\Big)\Big]\Big)^{\tfrac{1}{2q_2}}
        \\
        &
        {\color{red}
        \leq 
        C\sup_{s\in[0,T]} \Big(\mathbb{E}\Big[\exp\Big(\tfrac{U_0(X_s)}{e^{\alpha s }}+\int_0^{s}\tfrac{U_1(X_u)}{e^{\alpha u }}{\rm d}u\Big)\Big]\Big)^{\tfrac{1}{2q_1}}
        \cdot 
        \sup_{s\in[0,T]}\Big(\mathbb{E}\Big[\exp\Big(\tfrac{U_0(Y_s)}{e^{\alpha s }}+\int_0^{s\wedge\tau^e_{N}}\tfrac{U_1(Y_u)}{e^{\alpha u }}{\rm d}u\Big)\Big]\Big)^{\tfrac{1}{2q_1}}  
        }
        \\
        &\quad 
        {\color{red}
        \cdot
        \sup_{s\in[0,T]} \Big(\mathbb{E}\Big[\exp\Big(\tfrac{U_0(X_s)}{e^{\alpha s }}+\int_0^{s}\tfrac{U_1(X_u)}{e^{\alpha u }}{\rm d}u\Big)\Big]\Big)^{\tfrac{1}{2q_2}}
        \cdot 
        \sup_{s\in[0,T]}\Big(\mathbb{E}\Big[\exp\Big(\tfrac{U_0(Y_s)}{e^{\alpha s }}+\int_0^{s\wedge\tau^e_{N}}\tfrac{U_1(Y_u)}{e^{\alpha u }}{\rm d}u\Big)\Big]\Big)^{\tfrac{1}{2q_2}}
        }
        \\
        &
        \leq C\sup_{s\in[0,T]} \Big(\mathbb{E}\Big[\exp\Big(\tfrac{U_0(X_s)}{e^{\alpha s }}+\int_0^{s}\tfrac{U_1(X_u)}{e^{\alpha u }}{\rm d}u\Big)\Big]\Big)^{\tfrac{1}{2q}}
        \cdot 
        \sup_{s\in[0,T]}\Big(\mathbb{E}\Big[\exp\Big(\tfrac{U_0(Y_s)}{e^{\alpha s }}+\int_0^{s\wedge\tau^e_{N}}\tfrac{U_1(Y_u)}{e^{\alpha u }}{\rm d}u\Big)\Big]\Big)^{\tfrac{1}{2q}}  
        \\
        &< \infty.
        \end{aligned}
    \end{equation}
   
    Here the last inequality stands due to the exponential integrability property for both exact solution $\{X_s\}_{s\in[0,T]}$ and numerical solution $\{Y_s\}_{s\in[0,T]}$ (see \cite[Corollary 2.4]{cox2024local} and \cite[Corollary 2.10]{hutzenthaler2018exponential}).
    For any $p\geq 4$, by  \eqref{eq:term(4)_to_moment_estimate} and noting $U_1 \in \mathcal{C}^1_{\mathcal{P}}\big(\mathbb{R}^d,[0,\infty)\big)$ and $|x|^{1/c} \leq c(1+U_0(x))$, we get
    \begin{equation}
        \begin{aligned}
        &\sup_{s\in[0,T]}\left\|\mathbbm{1}_{s\leq \tau^e_{N}}\left[\tfrac{\left\langle X_{s}-Y_{s}, f\left(X_{s}\right)-f\left(Y_{s}\right)\right\rangle+\frac{1+\varepsilon}{2}\left\|g\left(X_{s}\right)-g\left(Y_{s}\right)\right\|^{2}}{\left|X_{s}-Y_{s}\right|^{2}}\right]^{+}\right\|_{L^{3p}(\Omega;\mathbb{R})}
        \\
        &\leq \sup_{s\in[0,T]}\Big\| \left[C + \tfrac{U_0(X_s)+U_0(Y_s)}{2q_1Te^{\alpha T}}+\tfrac{U_1(X_s)+U_1(Y_s)}{2q_2e^{\alpha T}}\right]\Big\|_{L^{3p}(\Omega;\mathbb{R})}
        \\
        &\leq C+C\Big[\sup_{s\in[0,T]}\|U_0(X_s)\|_{L^{3p}(\Omega;\mathbb{R})}+\sup_{s\in[0,T]}\|U_1(X_s)\|_{L^{3p}(\Omega;\mathbb{R})}
        \\
        &\ \ \ \ +\sup_{s\in[0,T]}\|U_0(Y_s)\|_{L^{3p}(\Omega;\mathbb{R})}+\sup_{s\in[0,T]}\|U_1(Y_s)\|_{L^{3p}(\Omega;\mathbb{R})}\Big]
        <\infty,
        \end{aligned}
    \end{equation}
    which confirms condition (b) in Theorem \ref{thm:main_thm}.
    By \eqref{eq:term(4)_to_moment_estimate} and condition $(1)$ in Theorem \ref{thm:stopped_tamed_method},
    \begin{equation}
        \begin{aligned}
                &\sup_{s\in[0,T]}\|X_s\|_{L^{6pc_g \vee3pc_f\vee3p}(\Omega;\mathbb{R}^{d})} {\textstyle \bigvee}
                \sup_{s\in[0,T]}\|Y_s\|_{L^{{6pc_g \vee3pc_f}}(\Omega;\mathbb{R}^{d})}< \infty
        \end{aligned}
    \end{equation}
    and for any $i=1,...d,\sup_{s\in[0,T]}\|\operatorname{Hess}_x(f^{(i)}(Y_s))\|_{L^{3p}(\Omega;\mathbb{R}^{d\times d})}<\infty$.
    By Lemma \ref{lem:tame_term_estimate},
    \begin{equation}
        \begin{aligned}
            &\sup_{s\in[0,T]}\|\mathbbm{1}_{s<\tau^e_{N}}a(s)\|_{L^{3p}(\Omega;\mathbb{R}^{d})}{\textstyle \bigvee} \sup_{s\in[0,T]}\|\mathbbm{1}_{s<\tau^e_{N}}b(s)\|_{L^{3p}(\Omega;\mathbb{R}^{d\times m})} < \infty.
        \end{aligned}
    \end{equation}
    This verifies condition (c) in Theorem \ref{thm:main_thm}, which in turn implies
    \begin{equation}\label{eq:estimate_convergence_rate}
        \begin{aligned}
                &\big\|\sup_{t\in[0,T]} |X_{t \wedge \tau^e_{N}}-Y_{t \wedge \tau^e_{N}}|\big\|_{{\color{red}L^{v}(\Omega;\mathbb{R})}}\leq C\Bigg[
                h^2 
                +
                h^{\frac{1}{2}}\int_0^{T} {{\int_{\lfloor s \rfloor_N}^s \Big\|\mathbbm{1}_{r\leq \tau^e_{N}}|f({Y_{\lfloor r \rfloor_N}}) - \mathbbm{1}_{r< \tau^e_{N}}a(r)|\Big\|_{{L^{p}(\Omega;\mathbb{R})}} {\rm d}r}}{\rm{d}}s
                \\
                & 
                \quad +
                \int_0^T\Big\|\mathbbm{1}_{s\leq \tau^e_{N}}\|{g({Y_s})-\mathbbm{1}_{s< \tau^e_{N}}b(s)}\| \Big\|^2_{{L^{p}(\Omega;\mathbb{R})}}{\rm d}s
                +
                \int_0^T \Big\|\mathbbm{1}_{s\leq \tau^e_{N}}\big|{f(Y_{\lfloor s \rfloor_N}) - \mathbbm{1}_{s< \tau^e_{N}}a(s)}\big|\Big\|^2_{{L^{p}(\Omega;\mathbb{R})}} {\rm{d}}s 
                \\
                & 
                \quad +
                h^{\frac{1}{2}}\int_0^{T} \Big( {{\int_{\lfloor s \rfloor_N}^s\Big\|\mathbbm{1}_{r\leq \tau^e_{N}}\|{g({Y_r})-\mathbbm{1}_{r< \tau^e_{N}}b(r)}\| \Big\|^2_{{L^{p}(\Omega;\mathbb{R})}}{\rm d}r }}\Big)^{\frac{1}{2}}{\rm{d}}s\Bigg]^{\frac{1}{2}}.
        \end{aligned}
    \end{equation}
    By the property of Lebesgue integral and Lemma \ref{lem:tame_term_estimate}, one can show
    \begin{equation}\label{eq:f(Y_s_)-a(s)}
        \begin{aligned}
            &\int_0^T\Big\|\mathbbm{1}_{s\leq \tau^e_{N}}\big|{f(Y_{\lfloor s \rfloor_N}) - \mathbbm{1}_{s<\tau^e_{N}}a(s)}\big|\Big\|^2_{{L^{p}(\Omega;\mathbb{R})}}{\rm d}s
            \\
            &=\int_0^T\Big\|\mathbbm{1}_{s< \tau^e_{N}}\big|{f(Y_{\lfloor s \rfloor_N}) - a(s)}\big|\Big\|^2_{{L^{p}(\Omega;\mathbb{R})}}{\rm d}s
            \\
            &\leq \int_0^T\Big\| \big|{f(Y_{\lfloor s \rfloor_N}) - \psi^{[1]}(Z_s)f(Y_{\lfloor s \rfloor_N})-\tfrac{1}{2} \sum_{j=1}^{m}\psi^{[2]}(Z_s)\Big(g(Y_{{\lfloor s \rfloor_N}})e_j,g(Y_{{\lfloor s \rfloor_N}})e_j\Big)}\big|\Big\|^2_{{L^{p}(\Omega;\mathbb{R})}}{\rm d}s
            \\
            &\leq Ch^{\delta-1}.
        \end{aligned}
    \end{equation}
    Moreover,
    \begin{equation}
        \begin{aligned}
            &\int_0^T\Big\|\mathbbm{1}_{s\leq \tau^e_{N}}\big\|{g(Y_{s}) - \mathbbm{1}_{s<\tau^e_{N}}b(s)}\big\|\Big\|^2_{{L^{p}(\Omega;\mathbb{R})}}{\rm d}s
            \\
            &\leq 
            C\int_0^T
                \Big\|\mathbbm{1}_{s< \tau^e_{N}}\big\|g(Y_{s}) - {g(Y_{\lfloor s \rfloor_N})}\big\|\Big\|^2_{{L^{p}(\Omega;\mathbb{R})}}
              +
              \Big\|\mathbbm{1}_{s< \tau^e_{N}}\big\|{g(Y_{\lfloor s \rfloor_N})}
                - 
              \psi^{[1]}(Z_s) {g(Y_{\lfloor s \rfloor_N})}             
              \big\|\Big\|^2_{{L^{p}(\Omega;\mathbb{R})}}
            {\rm d}s
            \\
            &\leq  C\int_0^T
                \Big\|\big\|g(Y_{s}) - {g(Y_{\lfloor s \rfloor_N})}\big\|\Big\|^2_{{L^{p}(\Omega;\mathbb{R})}}
                 {\rm d}s
                 +Ch^3.
        \end{aligned}
    \end{equation}
    For $i=1,...,d,j=1,...,m$, by the It\^{o} formula and 
    recalling $\tfrac{\partial g^{(k_2,j_2)}}{\partial x_{k_1}} g^{(k_1,j_1)} = 0$, we arrive at
   \begin{equation}
        \begin{aligned}
            &\int_0^T
                \Big\|\big\|g^{(ij)}(Y_{s}) - {g^{(ij)}(Y_{\lfloor s \rfloor_N})}\big\|\Big\|^2_{{L^{p}(\Omega;\mathbb{R})}}
                 {\rm d}s
        \\
        &\leq 
        C\int_0^T
            \Big\|\big\|
                \int_{\lfloor s \rfloor_N}^s
                    \big \langle g^{(ij)\prime}(Y_{r}),\mathbbm{1}_{r< \tau^e_N}a(r)\big \rangle
                    +\tfrac{1}{2}\operatorname{trace}\big(\mathbbm{1}_{r< \tau^e_N}b(r)^*\operatorname{Hess}_x(g^{(ij)}(Y_{r}))b(r)\big)
                     {\rm d}r
            \big\|\Big\|^2_{{L^{p}(\Omega;\mathbb{R})}}
        {\rm d}s
        \\
        &\quad +
        C\int_0^T
            \Big\|\big\|
                \int_{\lfloor s \rfloor_N}^s
                    \big \langle g^{(ij)\prime}(Y_{r}),\mathbbm{1}_{r< \tau^e_N}
                    \psi^{[1]}(Z_r) {g(Y_{\lfloor r \rfloor_N})}  
                     {\rm d}W_r
                     \big \rangle
            \big\|\Big\|^2_{{L^{p}(\Omega;\mathbb{R})}}
        {\rm d}s
        \\
        &\leq C\int_0^T
            \Big\|\big\|
                \int_{\lfloor s \rfloor_N}^s
                    \big \langle g^{(ij)\prime}(Y_{r})-g^{(ij)\prime}(Y_{\lfloor r \rfloor_N}),\mathbbm{1}_{r< \tau^e_N}
                    \psi^{[1]}(Z_r) {g(Y_{\lfloor r \rfloor_N})}  
                     {\rm d}W_r
                     \big\rangle
            \big\|\Big\|^2_{{L^{p}(\Omega;\mathbb{R})}}
        {\rm d}s
        \\
        &\quad +
        C\int_0^T
            \Big\|\big\|
                \int_{\lfloor s \rfloor_N}^s
                    \big \langle g^{(ij)\prime}(Y_{\lfloor r \rfloor_N}),\mathbbm{1}_{r< \tau^e_N}
                    \psi^{[1]}(Z_r) {g(Y_{\lfloor r \rfloor_N})}  
                     {\rm d}W_r
                     \big\rangle
            \big\|\Big\|^2_{{L^{p}(\Omega;\mathbb{R})}}
        {\rm d}s+Ch^2
        \\
        &\leq 
        C\int_0^T
            \Big\|\big\|
                \int_{\lfloor s \rfloor_N}^s
                \Big(
                     g^{(ij)\prime}(Y_{\lfloor r \rfloor_N})^*\mathbbm{1}_{r< \tau^e_N}
                    \psi^{[1]}(Z_r) {g(Y_{\lfloor r \rfloor_N})} 
                    - g^{(ij)\prime}(Y_{\lfloor r \rfloor_N})^*\mathbbm{1}_{r< \tau^e_N}
                    {g(Y_{\lfloor r \rfloor_N})}
                \Big)
                     {\rm d}W_r
            \big\|\Big\|^2_{{L^{p}(\Omega;\mathbb{R})}}
        {\rm d}s
        \\
        &\quad +Ch^2
        \\
        &\leq Ch^2.
        \end{aligned}
    \end{equation} 
    Therefore one obtains
    \begin{equation}\label{eq:g(Ys)-b(s)}
        \begin{aligned}
            &\int_0^T\Big\|\mathbbm{1}_{s\leq \tau^e_{N}}\big\|{g(Y_{s}) - \mathbbm{1}_{s<\tau^e_{N}}b(s)}\big\|\Big\|^2_{{L^{p}(\Omega;\mathbb{R})}}{\rm d}s \leq Ch^2.
        \end{aligned}
    \end{equation}
    The same arguments used in \eqref{eq:f(Y_s_)-a(s)} and \eqref{eq:g(Ys)-b(s)} can be applied to estimate 
    the second and {\color{red} fifth} terms on the right-hand side of \eqref{eq:estimate_convergence_rate}. Hence we deduce
    that for any $\tfrac{1}{v}=\tfrac{1}{p}+\tfrac{1}{q}$ and $\delta\geq3$
    \begin{equation}
        \begin{aligned}
        &\Big\| \sup_{t\in[0,T]} |X_{t \wedge \tau^e_{N}}-Y_{t \wedge \tau^e_{N}}| \Big\|_{{L^{v}(\Omega;\mathbb{R})}}
        &\leq Ch.
        \end{aligned}
    \end{equation}
    Observe  that 
    \begin{equation}
        \begin{aligned}
        &\big\|
        \sup_{t\in[0,T]}|X_t-Y_t|
        \big\|_{{L^{v}(\Omega;\mathbb{R})}}
        \\
        &\leq \big\|
        \sup_{t\in[0,T]}\mathbbm{1}_{\color{red}\tau^e_{N}<t}|X_t-Y_t|
        \big\|_{{L^{v}(\Omega;\mathbb{R})}}+\big\|\sup_{t\in[0,T]}\mathbbm{1}_{\color{red}\tau^e_{N}\geq t}|X_t-Y_t|\big\|_{{L^{v}(\Omega;\mathbb{R})}}
        \\
        &\leq \big\|\mathbbm{1}_{\tau^e_{N}<T}\big\|_{{L^{2v}(\Omega;\mathbb{R})}}
        \big\|\sup_{t\in[0,T]}|X_t-Y_t|\big\|_{{L^{2v}(\Omega;\mathbb{R})}}+\big\|\sup_{t\in[0,T]}|X_{t\wedge\tau^e_{N}}-Y_{t\wedge\tau^e_{N}}|\big\|_{{L^{v}(\Omega;\mathbb{R})}}.
        \end{aligned}
    \end{equation}
    Using Lemma \ref{lem:BDG_continuous_inequality} ensures
    \begin{equation}
        \begin{aligned}
        &\big\|\sup_{t\in[0,T]}|X_t-Y_t|\big\|_{{L^{2v}(\Omega;\mathbb{R})}}
        \\
        &\leq\Big\|\sup_{t\in[0,T]} \big| \int_0^t \big(f(X_s)-\mathbbm{1}_{s< \tau^e_{N}}a(s)\big){\rm d}s\big|
        +
        \sup_{t\in[0,T]} \big| \int_0^t \big(g(X_s)-\mathbbm{1}_{s< \tau^e_{N}}b(s)\big){\rm d}W_s\big|\Big\|_{{L^{2v}(\Omega;\mathbb{R})}}
        \\
        &\leq C,
        \end{aligned}
    \end{equation}
    and by the condition $(2)$ in Theorem \ref{thm:stopped_tamed_method}, the Markov inequality and $-\tfrac{1}{4!}x^4\geq -e^x,\ x\geq 0$ one infers
    \begin{equation}
        \begin{aligned}
            \mathbb{P}[\tau^e_{N}<T]
            &\leq \mathbb{P}
                \big[|Y_T|\geq          \exp(|\ln(h)|^{1/2})
                \big]
            \\
            &\leq \mathbb{P}\Big[\tfrac{1+U_0(Y_T)}{e^{\alpha T}}\geq \tfrac{1}{ce^{\alpha T}}\exp(\tfrac{1}{c}|\ln(h)|^{1/2})\Big]
            \\
            & \leq \mathbb{E}\Big[\exp\big(\tfrac{1+U_0(Y_T)}{e^{\alpha T}}\big)\Big] \exp\big(-\tfrac{1}{ce^{\alpha T}}\exp(\tfrac{1}{c}|\ln(h)|^{1/2})\big)
            \\
            &\leq C_1\exp\big(-\tfrac{|\ln(h)|^2}{24c^5e^{\alpha T}}\big).
        \end{aligned}
    \end{equation}
    For any $C_2>0$ and $h<1$ being small enough, 
    one knows  $|\ln(h)|^2\geq -\tfrac{1}{C_2}(2v\ln(h))$ and hence one gets for  small  $h < 1$,
    \[ 
      \exp\big( - \tfrac{|\ln(h)|^2}{24c^5e^{\alpha T}}\big)
      \leq h^{2v},
    \]
    which validates \eqref{eq:thm_stopped_tamed_method_conclu001}.
    Finally, note that if \eqref{eq:strengthed_condition_item_stopped_tame} holds, then for any $\gamma>0$,
    \begin{small}
    \begin{equation}
        \begin{aligned}
            &\bigg\|{\exp \bigg( \int_0^{\tau^e_{N}} {\left[\tfrac{\left\langle X_{s}-Y_{s}, f\left(X_{s}\right)-f\left(Y_{s}\right)\right\rangle+\frac{1+\varepsilon}{2}\left\|g\left(X_{s}\right)-g\left(Y_{s}\right)\right\|^{2}}{\left|X_{s}-Y_{s}\right|^{2}}\right]^{+}} {\rm{d}}s\bigg)}\bigg\|_{L^{\gamma}(\Omega;\mathbb{R})}< \infty.
        \end{aligned}
        \end{equation}
    \end{small}
    The assertion \eqref{eq:strengthed_conclusion_stopped_tame} can be acquired by repeating the above arguments, 
    which finishes  the proof.\qed

    In what follows we employ Theorem \ref{thm:stopped_tamed_method} to obtain the first-order strong convergence of the time-stepping scheme
    \eqref{eq:stop_tamed_method} for SDE models
    without globally monotone coefficients,
    taken from \cite{hutzenthaler2018exponential, Hutzenthaler2020}. 
    In the recent publication \cite{Hutzenthaler2020},
    the authors derived only a convergence rate of order $\tfrac12$ of the same scheme    \eqref{eq:stop_tamed_method}, even for the following additive noise driven SDE models 
    and multiplicative noise driven second order SDE models.
    Since the conditions of Theorem \ref{thm:stopped_tamed_method} are the same as {\color{red}
    those in \cite[Proposition 3.3]{Hutzenthaler2020}}, we just give the convergence results here and do not repeat the verification of the conditions.  Indeed, one can {\color{red}
    refer to  \cite[3.1.2,\ 3.1.6,\ 3.1.7,\ 3.1.3,\ 3.1.4]{Hutzenthaler2020} and \cite[Chapter 4]{cox2024local} }
    for details on the verification of the conditions for the following  different models. 
    The initial value $X_0$ of the following models is assumed to be  deterministic for simplicity.
    \\

    \textbf{Stochastic Lorenz equation with additive noise.}\ Let $d=m=3$ and $\alpha_1,\alpha_2,\alpha_3 \in [0,\infty)$. For $x=(x_1,x_2,x_3)^*\in \mathbb{R}^3$, we let 
    \begin{equation}
        f(x)=(\alpha_1(x_2-x_1),\alpha_2x_1-x_2-x_1x_3,x_1x_2-\alpha_3x_3)
    \end{equation} 
    and let $g( \cdot )$ be a constant matrix.
    Moreover, we take $U_0(x)=|x|^2$ and $U_1(x)=0$. 
    Then all conditions in Theorem \ref{thm:stopped_tamed_method} are fulfilled.
    Therefore, using the  SITEM method \eqref{eq:stop_tamed_sde} with $\delta\geq3$
    to solve the above stochastic Lorenz equation with additive noise yields that for any $r>0$, there exists a constant $C_r>0$ such that
    \begin{equation}
        \Big\| \sup_{t\in[0,T]}|X_t - Y_t| \Big \|_{L^{r}(\Omega;\mathbb{R})}\leq C_rh.
    \end{equation} 

    \textbf{Brownian dynamics.}\ Let $d=m\geq1,c,\beta\in(0,\infty)$ and $\theta \in [0,\tfrac{2}{\beta})$. 
    Assume that $V\in \mathcal{C}^3_{\mathcal{D}}\big(\mathbb{R}^d,[0,\infty)\big)\cap C^3\big(\mathbb{R}^d,[0,\infty)\big), V',\operatorname{Hess}_x \big(V^{(i)}\big)\in \mathcal{C}^1_{\mathcal{P}}\big(\mathbb{R}^d,\mathbb{R}^{d\times d}\big),i=1,...,d$ and 
    $\limsup_{r\searrow0} \sup_{z\in \mathbb{R}^d}
    $ $\tfrac{|z|^r}{1+V(z)}<\infty.$ For $x\in \mathbb{R}^d$, we set
    \begin{equation}
        f(x)=-(\nabla V)(x),g(x)=\sqrt{\beta}I_{\mathbb{R}^d\times\mathbb{R}^d}.
    \end{equation} 
    This equation is also termed as the overdamped Langevin dynamics in literature.
    In addition, we suppose that $(\Delta V)(x)\leq c+cV(x)+\theta\|(\nabla V)(x)\|^2$ and for any $\eta>0$
    \[
        \sup_{x,y\in \mathbb{R}^d,x\neq y} \Big[\tfrac{\langle x-y,(\nabla V)(y)-(\nabla V)(x)\rangle}{|x-y|^2}-\eta\big(V(x)+V(y)+|(\nabla V)(x)|^2+|(\nabla V)(y)|^2\big)\Big]<\infty.
    \]
Let $v \in (0,\tfrac{2}{\beta}-\theta)$, $U_0(x)=vV(x)$ and $U_1(x)=v(1-\tfrac{\beta}{2}(\theta+v))|(\nabla V)(x)|^2$. Then all conditions in Theorem \ref{thm:stopped_tamed_method} are fulfilled.
Therefore, applying the SITEM method \eqref{eq:stop_tamed_sde} ($\delta\geq3$) to the above Brownian dynamics yields 
    that for any $r>0$, there exists a constant $C_r>0$ such that
    \begin{equation}
        \Big \|
        \sup_{t\in[0,T]}|X_t-Y_t|
        \Big \|_{L^{r}(\Omega;\mathbb{R})}
        \leq 
        C_rh.
    \end{equation} 

\textbf{Langevin dynamics.}\ Let $d=2m\geq1,\gamma\in(0,\infty)$ and $\beta\in(0,\infty)$. 
    Assume that $V\in \mathcal{C}^3_{\mathcal{D}}\big(\mathbb{R}^m,[0,\infty)\big)\cap C^3\big(\mathbb{R}^m,[0,\infty)\big), V',\operatorname{Hess}_x \big(V^{(i)}\big)\in \mathcal{C}^1_{\mathcal{P}}\big(\mathbb{R}^m,\mathbb{R}^{m\times m}\big),i=1,...,m$ and 
    $\limsup$ $_{r\searrow0}
    $ $\sup_{z\in \mathbb{R}^m}\tfrac{|z|^r}{1+V(z)}<\infty.$ For $x=(x_1,x_2)^*\in \mathbb{R}^{2m},u\in\mathbb{R}^{m}$, we let
    \begin{equation} \label{eq:Langevin_dynamics}
        f(x)=(x_2,-(\nabla V)(x_1)-\gamma x_2),g(x)u=(0,\sqrt \beta u).
    \end{equation} 
    This equation is also termed as the underdamped Langevin dynamics in literature.
    In addition, we suppose that for any $\eta>0$
    \[
        \sup_{x,y\in \mathbb{R}^m,x\neq y} \Big[\tfrac{|(\nabla V)(x)-(\nabla V)(y)|}{|x-y|}-\eta\big(V(x)+V(y)+|x|^2+|y|^2\big)\Big]<\infty.
    \] 
Let $v \in (0,\infty)$, and for $x=(x_1,x_2)^*\in \mathbb{R}^{2m}, U_0(x)=\tfrac{v}{2}(|x_1|^2+|x_2|^2)+vV(x_1),U_1(x)=0$. 
Then all conditions in Theorem \ref{thm:stopped_tamed_method} are fulfilled.
Therefore using the  SITEM method \eqref{eq:stop_tamed_sde}($\delta\geq3$) for the Langevin dynamics
yields that for any $r>0$, there exists a constant $C_r>0$ such that
    \begin{equation}\label{eq:langev_dynam_conver_rate}
        \Big \|
        \sup_{t\in[0,T]}|X_t-Y_t|
        \Big\|_{L^{r}(\Omega;\mathbb{R})}
        \leq 
        C_rh.
    \end{equation} 

In \cite{cui2022density},
the authors proposed a splitting averaged vector field (AVF) scheme for the Langevin dynamics \eqref{eq:Langevin_dynamics}. Equipped with the exponential integrability properties of the implicit approximations $\{Y_n\}_{0 \leq n \leq N}$,
the authors of \cite{cui2022density} spent
a lot of efforts to analyze the pointwise strong error 
$
    \big( \sup_{0\leq n \leq N} 
    \mathbb{E} \big[ \| X_{t_n} - Y_{t_n} \|^p \big] \big)^{1/p}, p \geq 2.
$
As the first step, the pointwise strong convergence rate of order $\tfrac12$ was obtained,
which was later lifted to be order one by 
using technical arguments in the Malliavin calculus.
Instead, we analyze the pathwise uniformly 
strong error of an explicit time-stepping
scheme directly. 
By simply relying on the newly developed perturbation 
estimates in Section \ref{sect:perturbation-estimates}, we show a pathwise uniformly strong 
convergence rate of exact order one given by 
\eqref{eq:langev_dynam_conver_rate}.
{\color{blue}
It is worthwhile to mention that the authors of \cite{cui2022density} also proved the existence of the density function of the numerical solution produced by the splitting AVF scheme and provided the convergence rate of density functions for the scheme. 
Despite the same convergence rate, the splitting AVF scheme proposed by \cite{cui2022density}, as an implicit one, is expected to be more numerically stable than the explicit SITEM scheme, particularly for large step-sizes.
}

\textbf{Stochastic van der Pol oscillator.}\ Let $d=2,m\geq1, c,\alpha\in(0,\infty)$ and $\gamma,\beta\in[0,\infty)$. For $x=(x_1,x_2)^*\in \mathbb{R}^2,u \in \mathbb{R}^m$, we let
    \begin{equation}\label{eq:vander_Pol_oscilla}
        f(x)=(x_2,(\gamma-\alpha x^2_1)x_2-\beta x_1)^*,\ g(x)u=(0,\phi(x_1)u)^*,
    \end{equation}  
    where $\phi\in C^{2}(\mathbb{R}, \mathbb{R}^{1\times m})$ is a globally Lipschitz function.
    Let $v \in (0,\tfrac{\alpha}{2c})$, $U_0(x)=\tfrac{v}{2}|x|^2$ and $U_1(x)=v(\alpha-2cv)(x_1x_2)^2$. 
Then all conditions in Theorem \ref{thm:stopped_tamed_method} are fulfilled.
Therefore, applying the SITEM method \eqref{eq:stop_tamed_sde} ($\delta\geq3$) to the stochastic van der Pol oscillator  yields that, for any $r>0$, there exists a constant $C_r>0$ such that
    \begin{equation}
            \Big \|\sup_{t\in[0,T]}|X_t-Y_t| \Big \|_{L^{r}(\Omega;\mathbb{R})}\leq C_rh.
    \end{equation} 

    \textbf{Stochastic Duffing-van der Pol oscillator.}\ Let $d=2,m\geq1, \alpha_1,\alpha_2\in \mathbb{R}$ and $\alpha_3,c\in(0,\infty)$. For $x=(x_1,x_2)^*\in \mathbb{R}^2,u \in \mathbb{R}^m$, we let
    \begin{equation}
        f(x)=(x_2,\alpha_2x_2-\alpha_1x_1-\alpha_3x^2_1x_2-x^3_1)^*,\ g(x)u=(0,\phi(x_1)u)^*,
    \end{equation} 
    where $\phi\in C^{2}(\mathbb{R}; \mathbb{R}^{1\times m})$ is a globally Lipschitz function. 
    Let $v \in (0,\tfrac{\alpha_3}{c})$, $U_0(x)=\tfrac{v}{2}(\tfrac{x^4_1}{2}+x^2_2)$ and $U_1(x)=v(\alpha_3-cv)(x_1x_2)^2$. 
Then all conditions in Theorem \ref{thm:stopped_tamed_method} are fulfilled.
Therefore, applying the  SITEM method \eqref{eq:stop_tamed_sde} ($\delta\geq3$) to the stochastic Duffing-van der Pol oscillator yields that, for any $r>0$, there exists a constant $C_r>0$ such that
    \begin{equation}
        \Big\|
        \sup_{t\in[0,T]}|X_t-Y_t|
        \Big\|_{L^{r}(\Omega;\mathbb{R})}
        \leq C_rh.
    \end{equation}

    \textbf{Numerical experiments.}\ 
    {\color{blue}
    Now let us present some numerical experiments to test not only the strong convergence rate, but also the dynamic properties of the proposed method. We take the Langevin dynamics and the stochastic van der Pol oscillator as test examples. Let $T = 1$, $N = 2^{k}, k = 6, 7,..., 11$ and regard the fine approximations
    with $h_{\text{exact}} = 2^{-14}$
    as the "true" solution. 
    Also, we take $M = 5000$ Monte Carlo sample paths to approximate the expectation. 
    
    For the Langevin dynamics \eqref{eq:Langevin_dynamics}, we assign 
    $$
        m=1, \nabla V(x)=x^3-x, \gamma=1,\beta=2,X_0=(1,1)^*
    $$
    and  $\delta=3$ for the SITEM method.
    Such a type of potential $V(x_1) =\tfrac{1}{4}x_1^4-\tfrac{1}{2}x^2_1 $ is called double-well potential.
    Figure \ref{fig:Langevin_convergence_rate} displays the mean-square approximation errors of the SITEM method, the implicit splitting AVF method in \cite{cui2022density} and the implicit Euler method in \cite{talay2002stochastic}. Numerical results show that, the three methods all have a strong convergence rate of order one and the splitting AVF method is slightly better in terms of computational error. 
    \begin{figure}[H]
        \centering
        \includegraphics[scale=0.21]{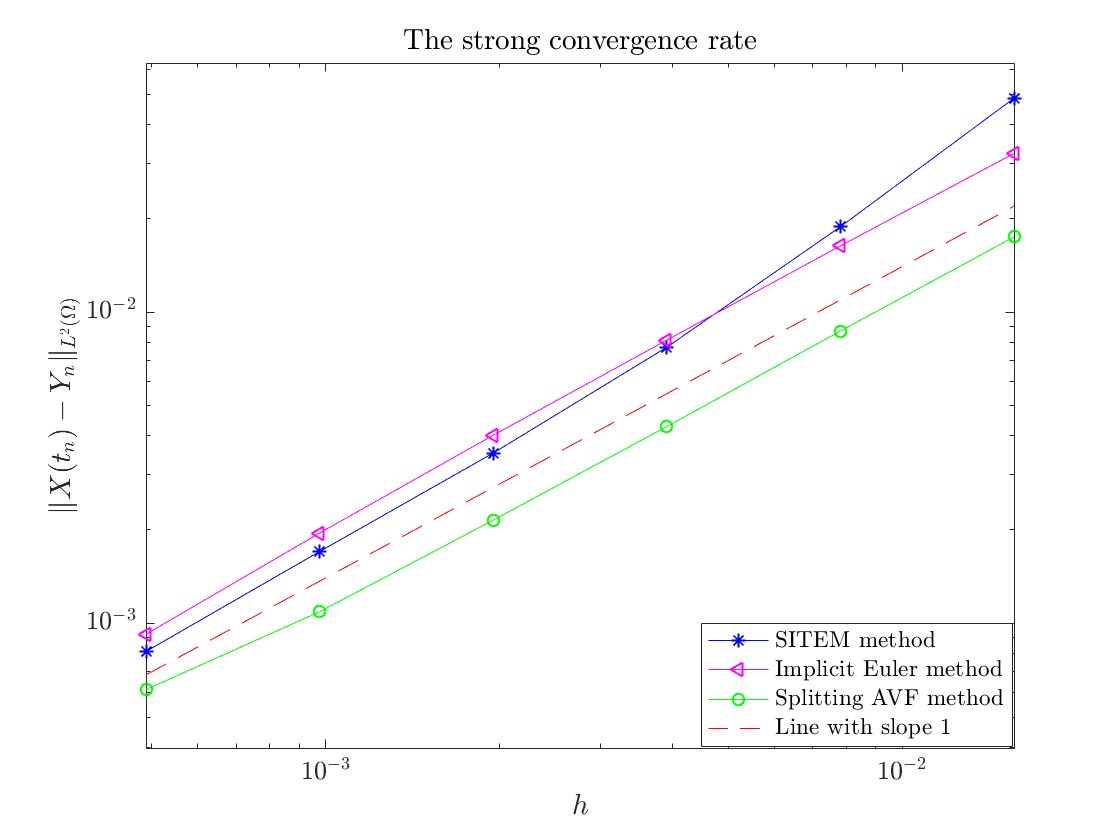}
        \caption{\small A comparison of strong convergence rates for  Langevin dynamics\label{fig:Langevin_convergence_rate}}
    \end{figure}
    Furthermore, it is known that 
    (see, e.g., \cite{mattingly2002ergodicity}), the Langevin dynamics admits a unique invariant distribution 
    \begin{equation}
        \mathbf{p}(x_1,x_2)=\Gamma_{\mathbf{p}}\exp(-V(x_1))\exp(\tfrac{1}{2}|x_2|^2),
    \end{equation} 
    where $\Gamma_{\mathbf{p}}$ is a normalization constant.
    Therefore, this equation is always used to sample from a target probability distribution $\pi (x_1) \propto e^{-V (x_1)}$. Next we test the ability of the SITEM method to sample from the distribution. We take a large time endpoint $T=500$ and test the SITEM in this paper, the implicit splitting AVF method in \cite{cui2022density} and the implicit Euler method in \cite{talay2002stochastic}. Numerical results are depicted in Figure \ref{fig:Langevin_density}, using two different stepsizes with $h=2^{-7}, 2^{-4}$. There one can observe that, these three methods all perform very well in the case of small stepsize $h=2^{-7}$. As the stepsize increases to $h=2^{-4}$, the implicit Euler method produces better approximations than the other two methods. The SITEM method and the implicit splitting AVF method perform similarly and give acceptable approximations. It should be noted that both the splitting AVF method and the implicit Euler method are implicit time-stepping schemes and their computational costs are more expensive than the SITEM method in the high-dimensional setting $m >1$.
    
    \begin{figure}[H]
    \centering  
        \subfigure[]{
        \includegraphics[width=7.2cm,height = 6cm]{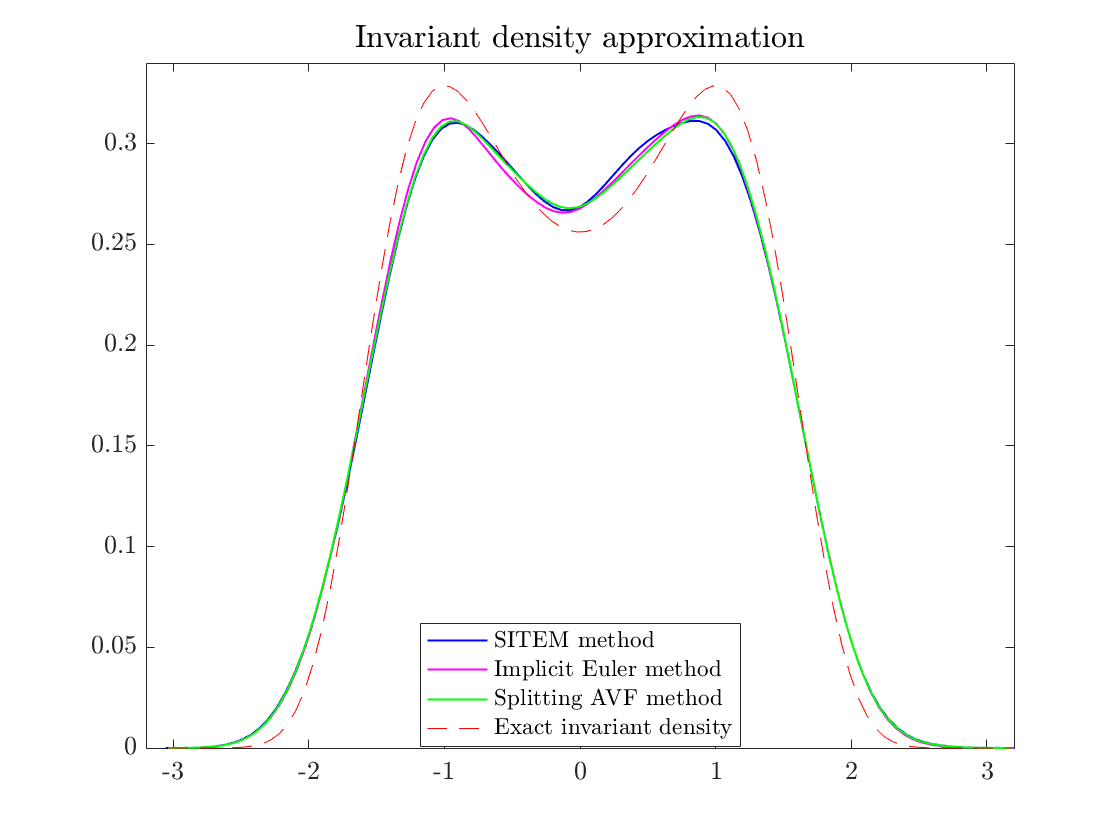}}
        \subfigure[]{
        \includegraphics[width=7.2cm,height = 6cm]{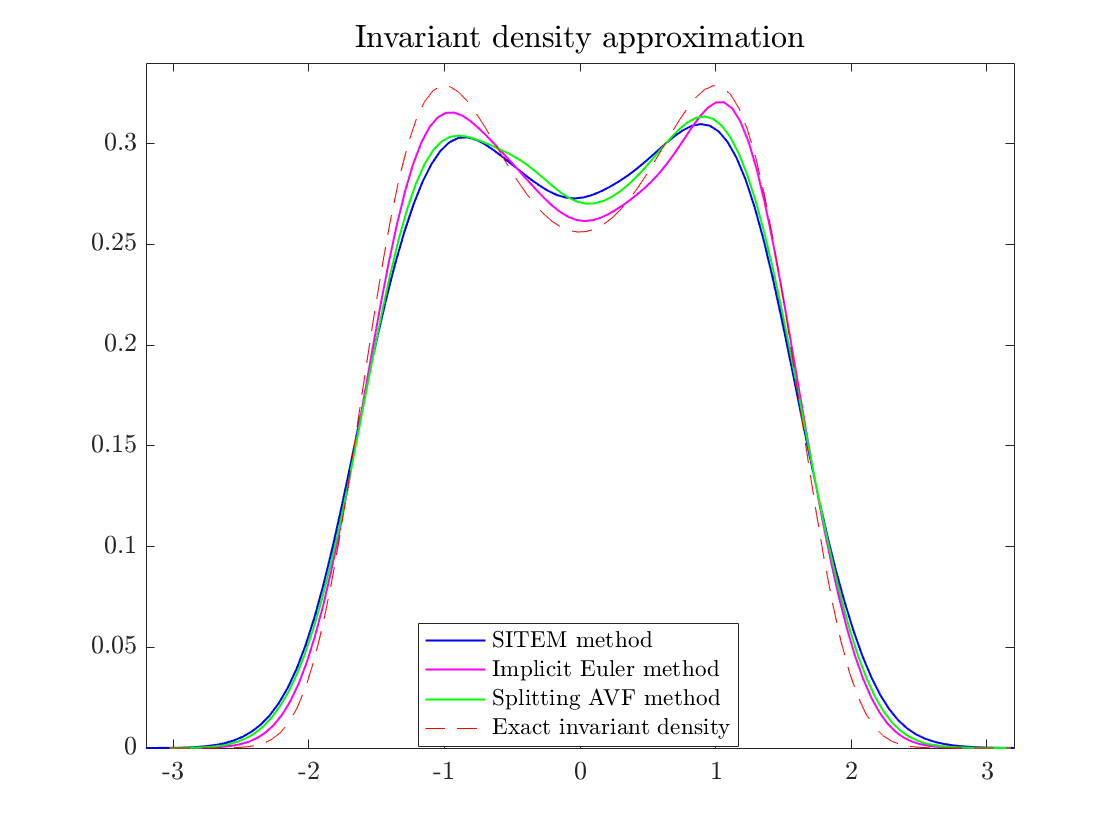}}
        \caption{\small Invariant density approximation: (a) $h=2^{-7}$; (b) $h=2^{-4}$}
        \label{fig:Langevin_density}
\end{figure}

We next turn to the stochastic van der Pol oscillator \eqref{eq:vander_Pol_oscilla} with coefficients 
    \begin{equation}\label{eq:vander_pol_parameter}
        m=1, \gamma=\alpha=0.2, \beta=1, X_0=(0.5,1.5)^*, \delta=3
    \end{equation}
in the case of the additive noise $\phi_1(x) \equiv \vartheta =\sqrt{0.1}$ and the multiplicative noise $\phi_2(x)=0.8x$. The mean-square approximation errors are presented in Figure \ref{fig:vander_pol_convergence_rate}, where one can observe order one convergence rate for both additive and multiplicative cases.
    \begin{figure}[H]
        \centering
        \includegraphics[scale=0.21]{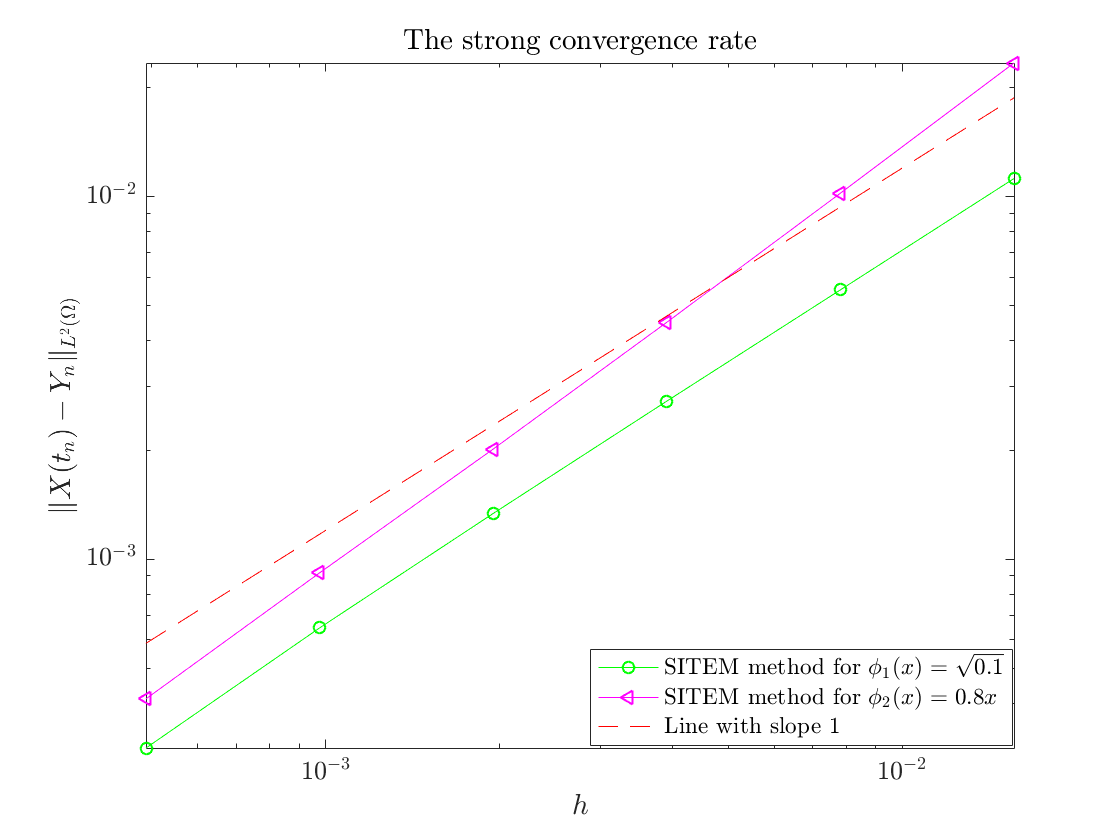}
        \caption{\small Strong convergence rate for stochastic van der Pol oscillator\label{fig:vander_pol_convergence_rate}}
    \end{figure}
\noindent Now let us focus on the case of additive noise $\phi(x) \equiv \vartheta=\sqrt{0.1}$. According to \cite[page 137]{to2000nonlinear}, one knows that, in our setting, the stochastic van der Pol oscillator model \eqref{eq:vander_Pol_oscilla} admits a stationary joint probability density
    \begin{equation}\label{eq:stationary_density_function}
    \mathbf{p}(x_{1},x_{2})=\Gamma_{\mathbf{p}} \exp\Big(-\tfrac{\alpha}{8\vartheta^2}
    \big((x_1^2+x^2_2)^2-8(x_1^2+x_2^2)\big)\Big),
   \end{equation} 
where $\Gamma_{\mathbf{p}}$ is a normalization constant.
Let $\mathbf{p}_1$ and $\mathbf{p}_2$ be marginal distribution of $\mathbf{p}$ in \eqref{eq:stationary_density_function}, defined by
\[
\mathbf{p}_1(x_{1}) : = \int_{\mathbb{R}} \mathbf{p} (x_1, x_2) \, {\rm d} x_2,
\quad
\mathbf{p}_2(x_{2}) : = \int_{\mathbb{R}} \mathbf{p} (x_1, x_2) \, {\rm d} x_1.
\]
It is obvious to observe that $\mathbf{p}(x_{1},x_{2})$ is symmetric with respect to variables $x_1,x_2$. As a consequence, $\mathbf{p}_1 \equiv \mathbf{p}_2$ and 
for polynomial functions $\Phi:\mathbb{R}\rightarrow \mathbb{R}$, 
\[
\mathbb{E} [ \Phi ( \xi_1 ) ] = \mathbb{E} [ \Phi ( \xi_2 ) ],
\quad
\xi_1 \sim \mathbf{p}_1(x_{1}), \
\xi_2 \sim \mathbf{p}_2(x_{2}).
\]
Moreover, for a fixed variable, $\mathbf{p}$ is an even function with respect to the other. Therefore, by particularly taking $\Phi=I$ (identity map) we obtain
\[
\mathbb{E} [  \xi_1 ] = \mathbb{E} [  \xi_2 ] = 0,
\quad
\xi_1 \sim \mathbf{p}_1(x_{1}), \
\xi_2 \sim \mathbf{p}_2(x_{2}).
\]
    This implies that, in the phase plane, ($\mathbb{E}[X_{1t}],\mathbb{E}[X_{2t}])$ will gradually tend to the trivial steady state $(0,0)$, as $t \rightarrow \infty$. 
    Unlike the deterministic case, the average oscillation period and limit cycle do not exist for the stochastic van der Pol oscillator, due to the presence of the noise.

In what follows we test the dynamics of numerical approximations produced by the SITEM method. Over the time interval $ [0,600]$,
Figures \ref{fig:evolution_expect_vander_Pol_small_h}, \ref{fig:evolution_expect_vander_Pol_big_h}  show the sample average trajectory and phase plane of the SITEM method for the stochastic van der Pol oscillator using two stepsizes $h = 2^{-7}, 2^{-4}$. From these figures, one can clearly see that, numerical approximations of ($\mathbb{E}[X_{1t}],\mathbb{E}[X_{2t}])$ produced by the SITEM method, even for a relatively large stepsize $h=2^{-4}$, tend to the trivial steady state $(0,0)$, as $t \rightarrow \infty$, reproducing the dynamics of the original model. Moreover, Figures \ref{fig:evolution_meansquare_expect_vander_Pol_small_h}, \ref{fig:evolution_meansquare_expect_vander_Pol_big_h} demonstrate numerical approximations of $\mathbb{E}[|X_{1t}|^2]$ and $\mathbb{E}[|X_{2t}|^2]$ with two stepsizes $h = 2^{-7}, 2^{-4}$, where one can observe that they all tend to some non-zero steady states, as $t \rightarrow \infty$. The above numerical experiments indicate that, using stepsizes of moderate size, the SITEM method is able to reproduce the dynamic property of the stochastic van der Pol oscillator.


    It is very interesting to mention that, for the deterministic van der Pol oscillator, i.e., $\phi(x)=0$, the exact solution is periodic and has a limit cycle. We use the SITEM method to numerically discretize it with various stepsizes. Our numerical results indicate that, both the period and limit cycle of the deterministic model are well reproduced by the SITEM method, even using a large stepsize $h=2^{-2}$.
\begin{figure}[H]
    \centering  
        \subfigure[]{
        \includegraphics[width=10.2cm,height = 5cm]{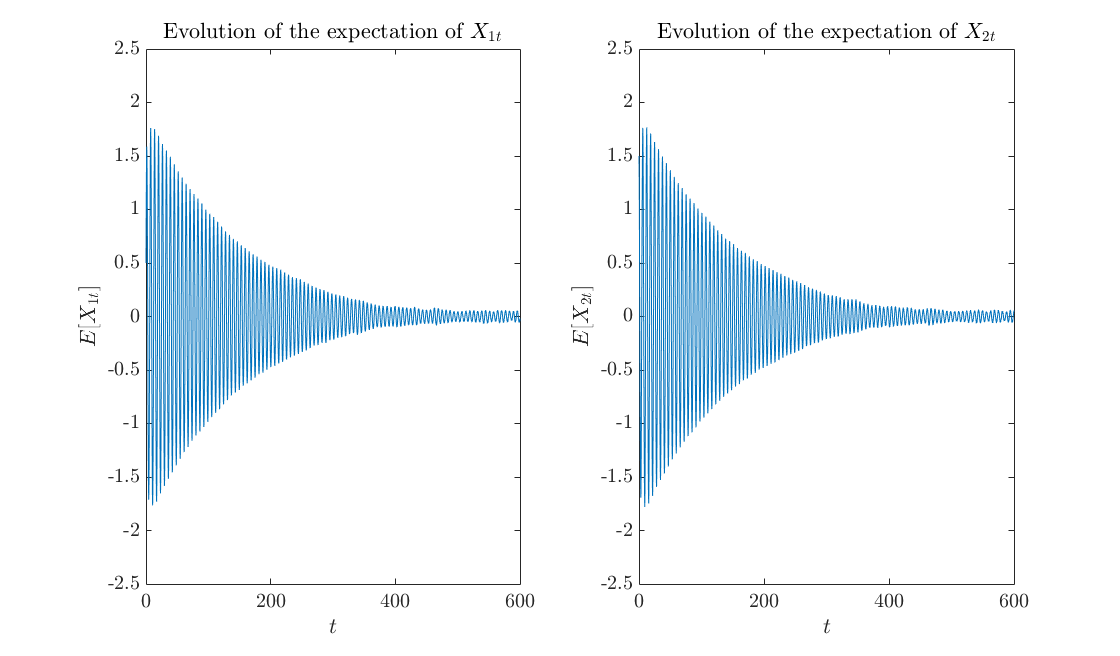}}
        \subfigure[]{
        \includegraphics[width=5.2cm,height = 5cm]{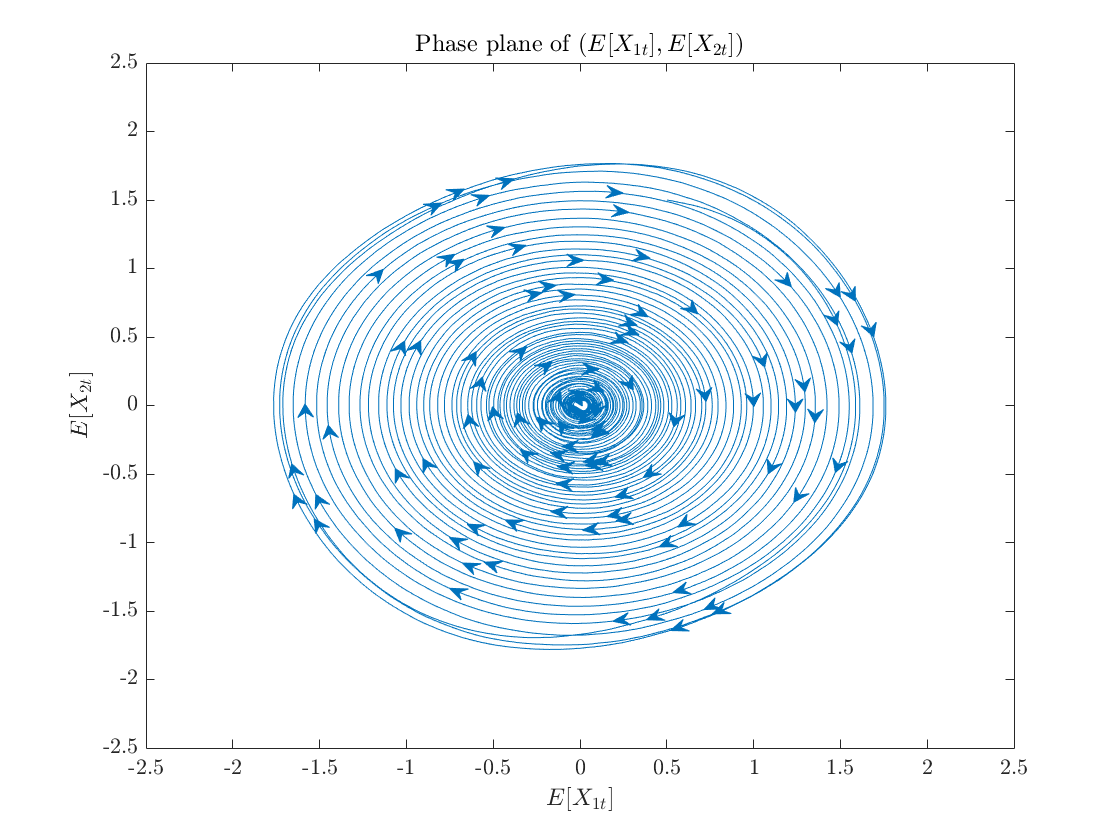}}
        \caption{\small  Sample average trajectory and phase plane  with $h=2^{-7}$}
        \label{fig:evolution_expect_vander_Pol_small_h} 
\end{figure}

   \begin{figure}[H]
    \centering  
        \subfigure[]{
        \includegraphics[width=10.2cm,height = 5cm]{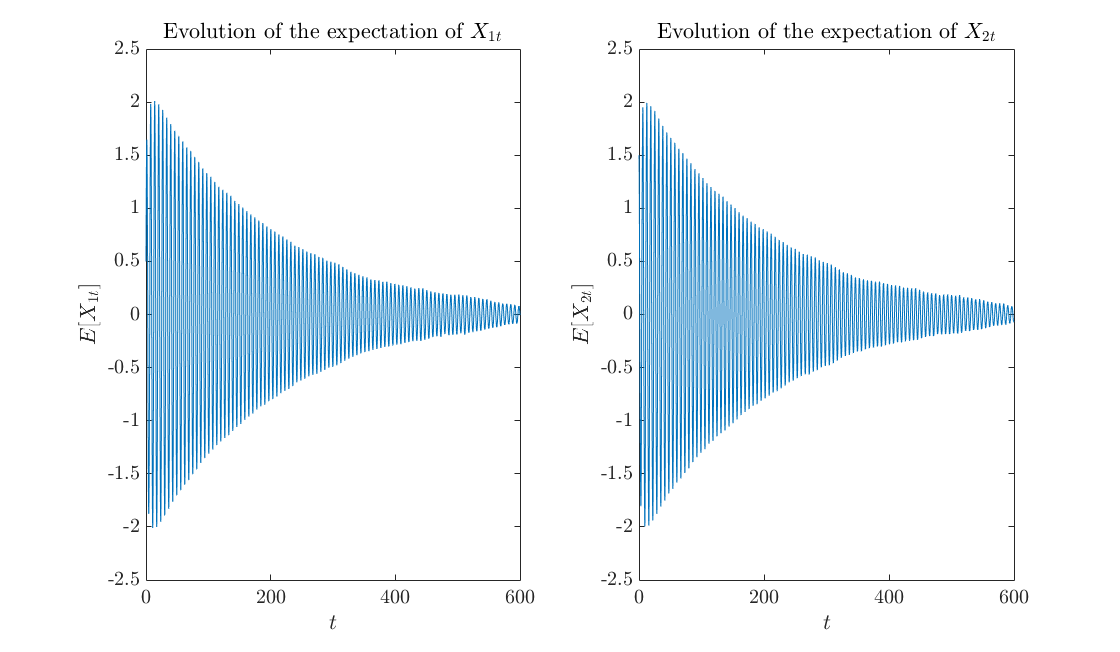}}
        \subfigure[]{
        \includegraphics[width=5.2cm,height = 5cm]{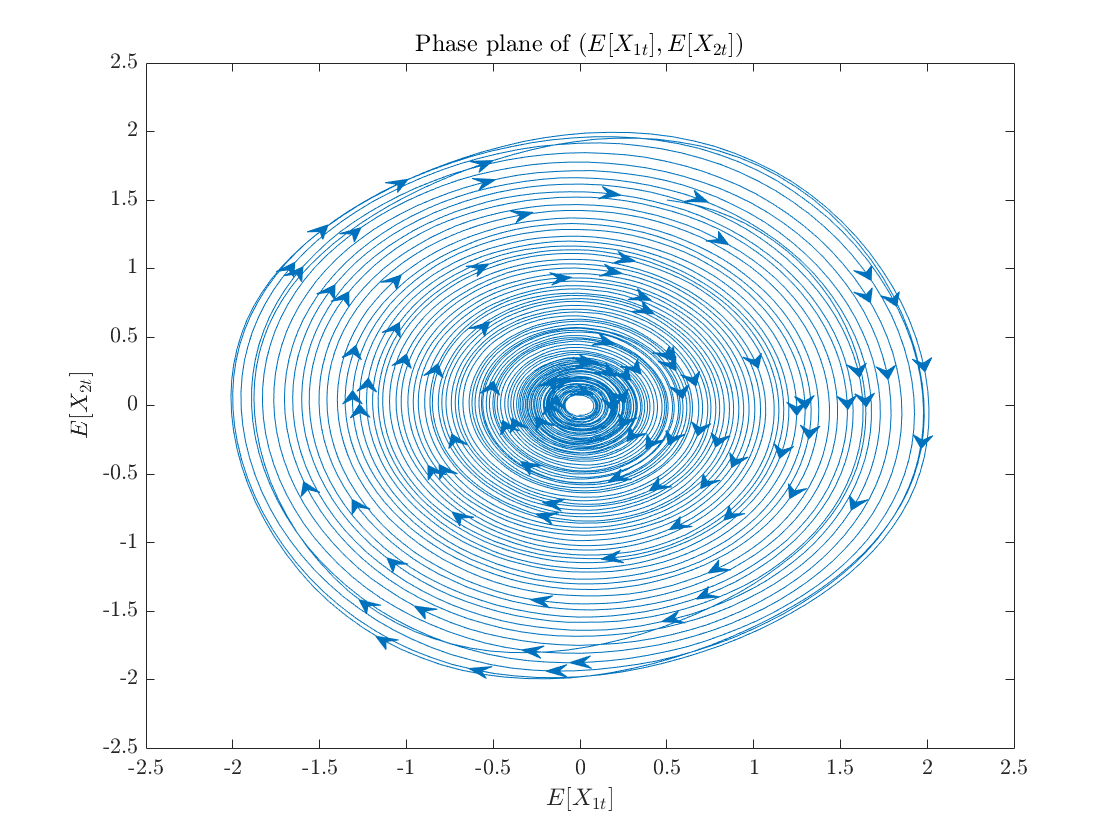}}
        \caption{\small Sample average trajectory and phase plane  with $h=2^{-4}$}
        \label{fig:evolution_expect_vander_Pol_big_h} 
\end{figure}

\begin{figure}[H]
        \centering
        \includegraphics[scale=0.28]{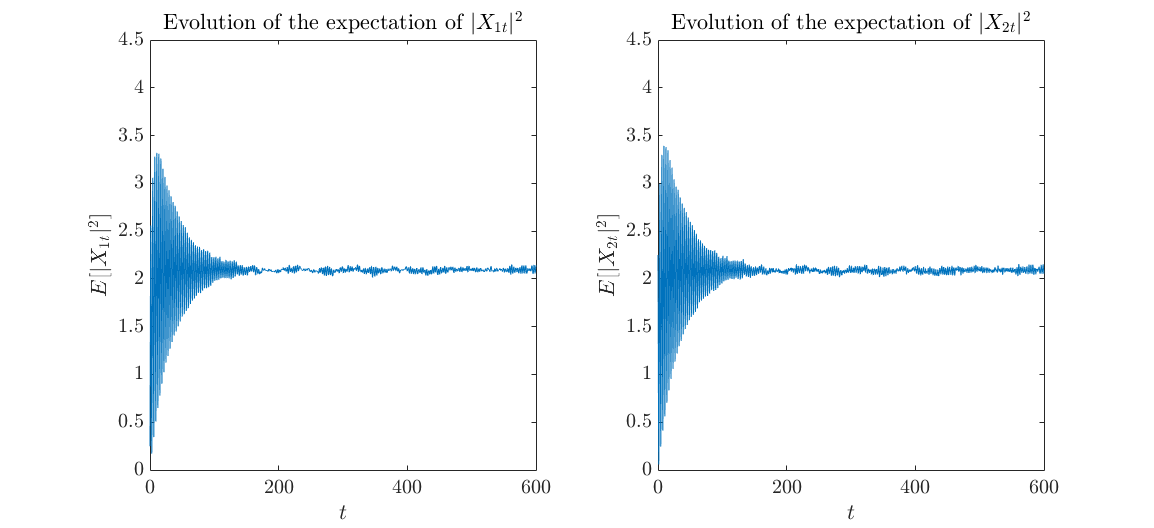}
        \caption{\small Numerical approximations of $\mathbb{E}[|X_{1t}|^2]$ and $\mathbb{E}[|X_{2t}|^2]$ with $h=2^{-7}$}
        \label{fig:evolution_meansquare_expect_vander_Pol_small_h} 
\end{figure}

\begin{figure}[H]
        \centering
        \includegraphics[scale=0.28]{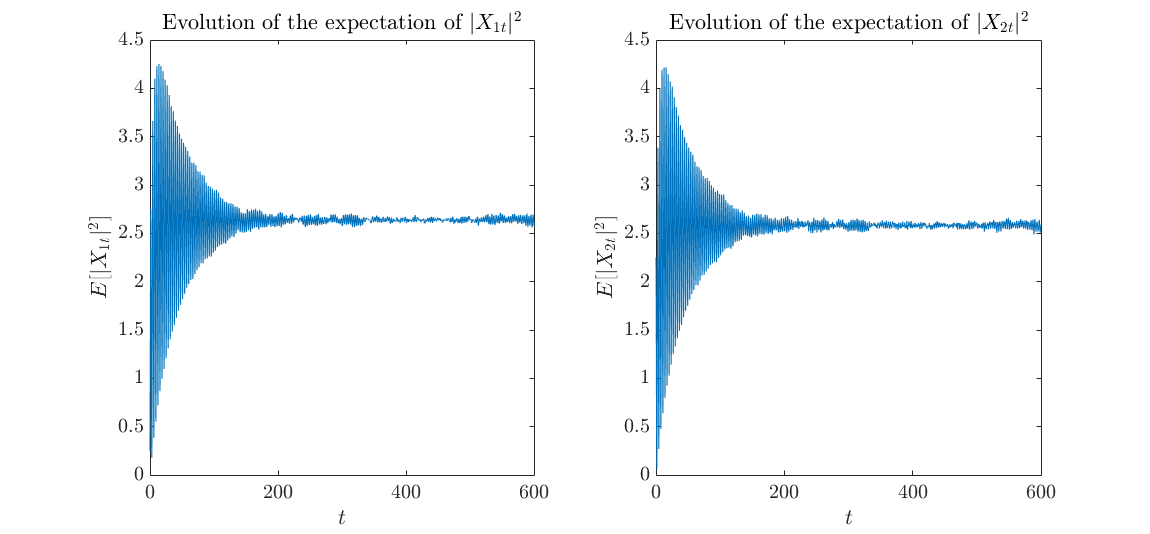}
        \caption{\small Numerical approximations of $\mathbb{E}[|X_{1t}|^2]$ and $\mathbb{E}[|X_{2t}|^2]$ with $h=2^{-4}$}
        \label{fig:evolution_meansquare_expect_vander_Pol_big_h} 
\end{figure}
    }

\section{A positive preserving Milstein type scheme for the stochastic LV competition model
with order-one pathwise uniformly strong convergence}
\label{sect:application3}
In this section we look at strong approximations
of the following  $d$-dimensional  stochastic Lotka--Volterra (LV) competition model for interacting multi-species in ecology \cite{mao2007stochastic,bahar2004stochastic,li2009population}:
{\color{blue}
    \begin{equation}\label{eq:L_V_model}
        \left\{ 
        \begin{array}{l}
            {\rm{d}}X_t 
            = 
            \text{diag}(X_t)(b - AX_t)
            {\rm{d}}t + \text{diag}(X_t)\sigma {\rm{d}}W_t,\ t\in [0,T],\\
            X_0 = x_0 \in {(\mathbb{R}^d)^ +},
        \end{array} \right.
    \end{equation}}
    where $T\in(0,\infty)$, 
    $b:=(b^{(i)})_{i=1,...,d}\in \mathbb{R}^d,\ A:=(a^{(ij)})_{i,j=1,...,d}\in \mathbb{R}^{d\times d},\ \sigma :=(\sigma^{(ij)})_{i=1,...,d,j=1,...,m}\in \mathbb{R}^{d\times m}$ and $(\mathbb{R}^d)^+ := \{x\in \mathbb{R}^d: x^{(1)}>0,..., x^{(d)}>0\}$. 
    Here the model is 
    driven by multi-dimensional noise and
    $\{ W_t\}_{t \in [0, T]}$ stands for
    a $m$-dimensional standard Brownian motion defined on $(\Omega,\mathcal{F},\{\mathcal{F}_t\}_{t\geq 0},\mathbb{P})$.
    For any vector $x \in \mathbb{R}^d$, we use $\text{diag}(x)$  to denote a $d \times d$ diagonal matrix whose principal diagonal is $x$.
    In order to show the well-posedness of the underlying model in $(\mathbb{R}^d)^+$, some assumptions are put on the elements of the matrix $A$.

\begin{assumption}\label{assump:LV_model}
        Every element of $A$ is non-negative and $\min_{1\leq i\leq d} \{a^{(ii)}\} >0.$ 
\end{assumption} 
Under the above assumption, the model \eqref{eq:L_V_model} has a unique global strong solution in $(\mathbb{R}^d)^+$.
    \begin{lemma}\label{lem:orig_LV_property}
        Under Assumption \ref{assump:LV_model}, there exists a unique  strong solution $\{X_t\}_{t\geq 0}$ for the equation \eqref{eq:L_V_model} staying in $(\mathbb{R}^d)^+$.
    \end{lemma}
    \indent \textbf{Proof}.  The well-posedness of the model in $(\mathbb R^d)^+$ has been
    established in \cite{bahar2004stochastic,mao2007stochastic}  for the scalar noise case. 
    For the present model driven by multi-dimensional noise, one can similarly prove it without any difficulty.\qed

Despite the existence and uniqueness of the positive solution to the stochastic LV model, the closed-form solution is not explicitly known
and efficient numerical approximations become an important tool
in applications.
Recently, several researchers proposed and 
analyzed positivity-preserving numerical 
schemes for such a typical multi-dimensional SDE 
model with  highly nonlinear and positive solution
(see, e.g., \cite{hong2021positivity,cai2023positivity,LI2023107260,mao2021positivity}).
%
%
It is worthwhile to point out that, under certain assumptions specified later, the highly nonlinear
drift coefficient $f (x) := \text{diag}( x ) 
( b - A x)$  and the linear diffusion coefficient
$g (x) := \text{diag}( x )  \sigma $
obey the Razumikhin-type growth condition
\[
\langle x,f(x) \rangle 
+
c \| g (x) \|^2
\le K ( 1 + |x |^2),
\]
but violate the
global monotonicity condition 
\begin{equation}
\label{eq:LVmodel-mono-condition}
\langle x - y,f(x) - f(y)\rangle  
+ 
c \|g(x) - g(y)\|^2 \le {K}|x - y{|^2},
\end{equation}
where {\color{blue} $x,y\in (\mathbb{R}^d)^+$}.
As already mentioned in the introduction part, 
the lack of the global monotonicity condition
causes an essential difficulty in obtaining
convergence rates of numerical approximations.
Usually, one needs to resort to exponential integrability of both the analytical and numerical solutions.
%
Very recently, the authors of \cite{LI2023107260} constructed a Lamperti transformed EM method for \eqref{eq:L_V_model} and used exponential integrability of both the analytical and numerical solutions to obtain
pointwise strong convergence rate of order $\tfrac12$ under some restrictive conditions. 
%

In this work we aim to propose a novel positivity preserving explicit Milstein-type method for the stochastic LV competition model and recover exactly order-one pathwise uniformly strong convergence of the new method under much relaxed conditions on the coefficients and stepsize, 
by relying on the use of previous perturbation estimates in Section \ref{sect:perturbation-estimates} (see Theorem \ref{thm:LV_convergence_rate}). To introduce 
the novel scheme, we regard the system \eqref{eq:L_V_model} as an interacting particle system of $d$ particles  evolving on the line.
For any $i\in\{1,...,d\}$, we consider a single particle of \eqref{eq:L_V_model} as follows:
    \begin{equation}\label{eq:one_dim_LV_model}
         {\rm{d}} X^{(i)}_t=\big(X^{(i)}_tb^{(i)}-X^{(i)}_t\sum_{j=1}^{d}a^{(ij)}X^{(j)}_t\big){\rm d}t+X^{(i)}_t\sum_{j=1}^m\sigma^{(ij)}  {\rm{d}}W^{(j)}_t.
    \end{equation}
In order to numerically approximate \eqref{eq:one_dim_LV_model}
on a uniform grid $\{t_n = n h\}_{0 \leq n \leq N}$ with stepsize $h = \tfrac{T}{N}$,
we propose the following linear-implicit (explicit) Milstein method starting from $Y_0=X_0$:
    \begin{equation}\label{eq:linear_mil_method}
        \begin{aligned} Y^{(i)}_{t_{n+1}}=&Y^{(i)}_{t_{n}}+Y^{(i)}_{t_{n+1}}\big(b^{(i)}-\sum_{j=1}^{d}a^{(ij)}Y^{(j)}_{t_{n}}\big)h+Y^{(i)}_{t_{n}}\sum_{j=1}^{m}\sigma^{(ij)}\Delta W^{(j)}_{t_n}
            \\
            &+\tfrac{1}{2}Y^{(i)}_{t_{n}}\sum_{j_1=1}^{m}\sum_{j_2=1}^{m}\sigma^{(ij_1)}\sigma^{(ij_2)}\Delta W^{(j_1)}_{t_n}\Delta W^{(j_2)}_{t_n} -\tfrac{1}{2}Y^{(i)}_{t_{n+1}}\sum_{j=1}^{m}(\sigma^{(ij)})^2h,
        \end{aligned} 
    \end{equation}
    where for $n=0,...,N-1,j=1,...,m, \Delta W^{(j)}_{t_n}:=\Delta W^{(j)}_{t_{n+1}}-\Delta W^{(j)}_{t_n}$. 
    Further, some elementary rearrangements turn \eqref{eq:linear_mil_method} into
    \begin{equation}\label{eq:rearrange_milstein_no_Qi}
        \begin{aligned}
        &\Big(1 - {b^{(i)}}h + h\sum\limits_{j = 1}^d {{a^{(ij)}}Y^{(j)}_{t_n}}  + \tfrac{ 1}{2}\sum_{j=1}^{m}{(\sigma^{(ij)})^2}h\Big)Y^{(i)}_{t_{n+1}} 
        \\
        &= Y^{(i)}_{t_n}\Big(1 + \sum_{j=1}^{m}{\sigma^{(ij)}}\Delta W^{(j)}_{t_n}
        +\tfrac{1}{2}\sum_{j_1=1}^{m}\sum_{j_2=1}^{m}\sigma^{(ij_1)}\sigma^{(ij_2)}\Delta W^{(j_1)}_{t_n}\Delta W^{(j_2)}_{t_n}\Big).
        \end{aligned} 
    \end{equation}
%
%
Under mild assumptions on the stepsize $h > 0$,
one can readily check that  the proposed scheme
is well-posed and positivity preserving.
%
\begin{proposition}[Positivity preserving]\label{propo:Positivity_milstein_method}
        Let Assumption \ref{assump:LV_model} be satisfied and let $X_0 = x_0 \in {(\mathbb{R}^d)^ +}$.
        For  some $\gamma>1$ and
        $0<h \leq \min_{1\leq i\leq d}
        \{
        \tfrac{1}{\gamma(b_i - \frac12\sum_{j=1}^{m}
                {(\sigma^{(ij)})^2})
                    \vee 0
                }
               \wedge
               T\} \ (\tfrac{1}{0}:=\infty)$,
        the proposed scheme \eqref{eq:linear_mil_method} is well-defined and has a unique positive solution
        in $(\mathbb R^d)^+$.
\end{proposition}
\indent \textbf{Proof}. 
Given $Y_{t_n} \in  (\mathbb{R}^d)^+$, by Assumption \ref{assump:LV_model} and the range of $h$  
    , it is clear to see 
    \[
        0<\Big(1 - {b^{(i)}}h + h\sum\limits_{j = 1}^d {{a^{(ij)}}Y^{(j)}_{t_n}}  + \tfrac{ 1}{2}\sum_{j=1}^{m}{(\sigma^{(ij)})^2}h\Big)^{-1}\leq \tfrac{\gamma}{\gamma-1},
    \]
    which implies  equation \eqref{eq:linear_mil_method} is well-defined and for any $ 1 \leq i \leq d$,
    \begin{equation}\label{eq:itera_milstein_no_Q_i}
      \begin{aligned}
        Y^{(i)}_{t_{n+1}}
        &=Y^{(i)}_{t_n}\Big(\tfrac{1}{2}\big(1+\sum_{j=1}^{m}{\sigma^{(ij)}}\Delta W^{(j)}_{t_n}\big)^2+\tfrac{1}{2}\Big)\Big(1 - {b^{(i)}}h + h\sum\limits_{j = 1}^d {{a^{(ij)}}Y^{(j)}_{t_n}}  + \tfrac{ 1}{2}\sum_{j=1}^{m}{(\sigma^{(ij)})^2}h\Big)^{-1}
        > 0.
      \end{aligned}
    \end{equation}
Since $Y_0 \in  (\mathbb{R}^d)^+$, one infers that $Y_{t_n} \in  (\mathbb{R}^d)^+$ for any $n=0,...,N-1$, which ensures that the linear-implicit Milstein scheme is positivity-preserving and the  proof is thus completed.\qed 
%
%

Before coming to the convergence analysis of 
the scheme, we would like to mention an interesting observation for the particular stochastic LV model.
Although the global monotonicity condition \eqref{eq:LVmodel-mono-condition} does not
hold, the drift coefficients of the model obey
a special form of locally monotonicity condition
\eqref{eq:Lv-model-f-locally-monot}, where
the control terms on the right-hand side (cf. $U_0 (x), U_1(x)$) depend only on $x$ and does not depend on $y$. This interesting finding
helps us to derive the convergence rate more easily,
without requiring the exponential integrability properties of the numerical approximations. 
As a direct consequence of Theorem \ref{thm:main_thm}, we formulate the following proposition for approximations of SDEs fulfilling the special case of locally monotonicity condition
\eqref{eq:Lv-model-f-locally-monot}.
%
\begin{proposition}\label{prop:LV_from_main_corol} 
    Let $f:\mathbb R^d\rightarrow \mathbb R^d,g:\mathbb R^d\rightarrow \mathbb R^{d \times m} $ be measurable functions and $f\in C^2{(\mathbb R^d,\mathbb R^d})$. 
    Further, let $f\in \mathcal{C}^1_{\mathcal{P}}(\mathbb{R}^d,\mathbb{R}^{d})$ with constants $K_f,c_{f}$ and $g$ be Lipschitz. 
    For a set $D_X \subset \mathbb{R}^d$, we assume that $X \colon \Omega \times [0,T]\rightarrow D_X$
    and $Y \colon \Omega \times [0,T]\rightarrow D_X$ 
      be defined by \eqref{eq:typical_sde} and \eqref{eq:approxi_sde} with continuous sample paths, respectively, satisfying  $\xi_Y=\xi_X=X_0$.  
     Let $U_0 \in C^2\big(\mathbb{R}^d,[0,\infty)\big), U_1 \in C\big(\mathbb{R}^d,[0,\infty)\big)$, and let $c,v,q,T\in (0,\infty),\alpha \in[0,\infty), p\geq 4 $. Besides, suppose that for all $x,y\in \mathbb R^d$,
    \begin{enumerate}[{\rm(1)}]
        \item 
            there exist constants $L,\kappa\geq0$ such that for any $i=1,...,d$,
                $$
                |U_0(x)|
                {\textstyle \bigvee}
                    \|\operatorname{Hess}_x(f^{(i)}(x))\|
                    \leq 
                    L(1+|x|)^{\kappa};
                $$
        \item
            $|x|^{1/c} \leq c(1+U_0(x))$ and $\mathbb{E}\left[e^{U_0(X_0)}\right]<\infty;$
        \item
            $(\mathcal{A}_{f,g} U_0 )(x)+\tfrac{1}{2}|g(x)^{*}(\nabla U_0(x))|^2+U_1(x) \leq c+\alpha U_0(x);$
        \item
            for any $\eta>0$, there exists a constant $K_{\eta}$ such that
                \begin{equation}
                \label{eq:Lv-model-f-locally-monot}
                    \langle x-y,f(x)-f(y)\rangle \leq \big[K_{\eta}+\eta(U_0(x)+U_1(x))\big]|x-y|^2,\ x,y\in D_X;
                \end{equation}
        \item
        {\color{black}
        for any $\theta \geq 1$,
        \begin{equation}
            \begin{aligned}
                    \sup_{s\in[0,T]}\|a(s)\|_{L^{\theta}(\Omega;\mathbb{R}^{d})}{\textstyle \bigvee} 
                    \sup_{s\in[0,T]}\|b(s)\|_{L^{\theta}(\Omega;\mathbb{R}^{d\times m})}
                    <K_{\theta},
            \end{aligned}
        \end{equation}
        where $K_{\theta}>0$ is independent of $h$.}
    \end{enumerate}
    Then for $\tfrac{1}{v}=\tfrac{1}{p}+\tfrac{1}{q}$, the approximation \eqref{eq:approxi_sde}  of \eqref{eq:typical_sde} admits 
    \begin{equation}\label{eq:other_methods}
        \begin{aligned}
            &\big\|\sup_{t\in[0,T]} |X_{t}-Y_{t}|\big\|_{{L^{v}(\Omega;\mathbb{R})}} \leq 
            C\bigg[h^2 
            + 
            \int_0^T\Big\|\mathbbm{1}_{s\leq \tau_N}\|{g({Y_s})-b(s)}\| \Big\|^2_{{L^{p}(\Omega;\mathbb{R})}}{\rm d}s
            \\
            &\quad\quad +h^{\frac{1}{2}}\int_0^{T} {{\int_{\lfloor s \rfloor_N}^s \Big\|\mathbbm{1}_{r\leq \tau_N}|f({Y_{\lfloor r \rfloor_N}}) - a(r)|\Big\|_{{L^{p}(\Omega;\mathbb{R})}} {\rm d}r}}{\rm{d}}s
            +\int_0^T \Big\|\mathbbm{1}_{s\leq \tau_N}\big|{f(Y_{\lfloor s \rfloor_N}) - a(s)}\big|\Big\|^2_{{L^{p}(\Omega;\mathbb{R})}} {\rm{d}}s
            \\
            &\quad\quad +
            h^{\frac{1}{2}}\int_0^{T} \Big( {{\int_{\lfloor s \rfloor_N}^s\Big\|\mathbbm{1}_{r\leq \tau_N}\|{g({Y_r})-b(r)}\| \Big\|^2_{{L^{p}(\Omega;\mathbb{R})}}{\rm d}r }}\Big)^{\frac{1}{2}}{\rm{d}}s\bigg]^{\frac{1}{2}},
        \end{aligned}
    \end{equation}
    where $C$ is independent of $h$.
    \end{proposition}
    \indent \textbf{Proof}. 
    Proposition \ref{prop:LV_from_main_corol} can be directly obtained by using Theorem \ref{thm:main_thm}. To see this, 
    setting $\tau_N=T$  in Theorem \ref{thm:main_thm} and utilizing 
    the same arguments as used in \eqref{eq:term(4)_to_moment_estimate} one deduces that
    \begin{equation}\label{eq:prove_in_proposi_LV_from_main_corol}
        \begin{aligned}
            &\bigg\|{\exp \bigg( \int_0^{T} {\left[\tfrac{\left\langle X_{s}-Y_{s}, f\left(X_{s}\right)-f\left(Y_{s}\right)\right\rangle+\frac{1+\varepsilon}{2}\left\|g\left(X_{s}\right)-g\left(Y_{s}\right)\right\|^{2}}{\left|X_{s}-Y_{s}\right|^{2}}\right]^{+}} {\rm{d}}s\bigg)}\bigg\|_{L^{q}(\Omega;\mathbb{R})}< \infty.
        \end{aligned}
    \end{equation}
   In view of \eqref{eq:prove_in_proposi_LV_from_main_corol} and conditions (1),(2),(4) in Proposition \ref{prop:LV_from_main_corol}, we arrive at
    \begin{equation}
        \begin{aligned}
            \sup_{s\in[0,T]}\left\|\mathbbm{1}_{s\leq \tau_N}\left[\tfrac{\left\langle X_{s}-Y_{s}, f\left(X_{s}\right)-f\left(Y_{s}\right)\right\rangle+\frac{1+\varepsilon}{2}\left\|g\left(X_{s}\right)-g\left(Y_{s}\right)\right\|^{2}}{\left|X_{s}-Y_{s}\right|^{2}}\right]^{+}\right\|_{L^{3p}(\Omega;\mathbb{R})} <\infty
        \end{aligned}
    \end{equation}
    and
    \begin{equation}
        \sup_{s\in[0,T]}\|X_s\|_{L^{6pc_g \vee3pc_f\vee3p}(\Omega;\mathbb{R}^{d})}<\infty.
    \end{equation}
    Combining Theorem \ref{thm:main_thm} 
    with conditions (1),(5) in Proposition \ref{prop:LV_from_main_corol} then completes the proof.\qed
    
%
    It is worth noting that the restrictions on $U_0$ and $U_1$ in Proposition \ref{prop:LV_from_main_corol} are more relaxed than those in Theorem \ref{thm:stopped_tamed_method} 
    due to the particular condition (4) (see \cite[Corollary 2.4]{cox2024local}).
%
%
In what follows, we utilize Proposition \ref{prop:LV_from_main_corol} to prove the strong convergence rate of the newly  developed 
Milstein type method. 
For simplicity of presentation, we denote
    \[
      \underline{a}
      := 
      \min\nolimits_{1\leq i\leq d} \{a^{(ii)}\} ,\ \bar{b}:= \max\nolimits_{1\leq i\leq d}\{|b^{(i)}|\};\ \bar \sigma:= \max\nolimits_{1\leq i\leq d,1\leq j\leq m}\{|\sigma^{(ij)}|\}
    \]
    and
    \begin{equation}\label{eq:Qi_def}
        Q^{(i)}_{t_n}:=1 - {b^{(i)}}h + h\sum\limits_{j = 1}^d {{a^{(ij)}}Y^{(j)}_{t_n}+ \tfrac{1}{2}\sum_{j=1}^{m}{(\sigma^{(ij)})^2}h},\ i\in \{ 1,...,d \}.
    \end{equation}
    Under conditions in Proposition
    \ref{propo:Positivity_milstein_method},
    one knows 
 \begin{equation}
 \label{eq:Q-inverse-boundedness}
 0<(Q^{(i)}_{t_n})^{-1}
            \leq 
            \min \Big \{
            \tfrac{\gamma}{\gamma-1},
            1+\tfrac{\gamma}{\gamma-1} \big |b_i - \tfrac{1}{2}\sum_{j=1}^{m}(\sigma^{(ij)})^2 \big|h
            \Big\}.
\end{equation}
In the notation of $Q^{(i)}_{t_n}$,  one can  
rewrite the scheme \eqref{eq:rearrange_milstein_no_Qi} as
{\color{red}
    \begin{equation}\label{eq:with_Qi_equation}
        \begin{aligned}
            Y^{(i)}_{t_{n+1}}
            & =
            Y^{(i)}_{t_n}
            \Big(
            1  + 
                \sum_{j=1}^{m}{\sigma^{(ij)}}\Delta W^{(j)}_{t_n}
               +
                \tfrac{1}{2}\sum_{j_1=1}^{m}\sum_{j_2=1}^{m}\sigma^{(ij_1)}\sigma^{(ij_2)}\Delta W^{(j_1)}_{t_n}\Delta W^{(j_2)}_{t_n}
            \Big)
            (Q^{(i)}_{t_n})^{-1}
            \\
            &=
             Y^{(i)}_{t_n}
             +
            Y^{(i)}_{t_n}\big( (Q^{(i)}_{t_n})^{-1}-1\big)
            \\
            &\quad +
             Y^{(i)}_{t_n}
            \Big(
                \sum_{j=1}^{m}{\sigma^{(ij)}}\Delta W^{(j)}_{t_n}
               +
                \tfrac{1}{2}\sum_{j_1=1}^{m}\sum_{j_2=1}^{m}\sigma^{(ij_1)}\sigma^{(ij_2)}\Delta W^{(j_1)}_{t_n}\Delta W^{(j_2)}_{t_n}
            \Big)
            (Q^{(i)}_{t_n})^{-1}
            \\
            &=
            Y^{(i)}_{t_n}
             +
            Y^{(i)}_{t_n}(Q^{(i)}_{t_n})^{-1}
            \big( 
            b^{(i)}h-h\sum\limits_{j = 1}^d {a^{(ij)}}Y^{(j)}_{t_n}
            \big)
            \\
            &\quad +
             Y^{(i)}_{t_n}
            \Big(
                \sum_{j=1}^{m}{\sigma^{(ij)}}\Delta W^{(j)}_{t_n}
               +
                \tfrac{1}{2}\sum_{j_1=1}^{m}\sum_{j_2=1}^{m}\sigma^{(ij_1)}\sigma^{(ij_2)}\Delta W^{(j_1)}_{t_n}\Delta W^{(j_2)}_{t_n}
                -
                \tfrac{1}{2}
                \sum_{j=1}^{m}{(\sigma^{(ij)})^2}h
            \Big)
            (Q^{(i)}_{t_n})^{-1}
            \\
            & =
             Y^{(i)}_{t_n}+
             \int_{t_n}^{t_{n+1}} 
             {\tfrac{{{b^{(i)}} - \sum\limits_{j = 1}^d {{a^{(ij)}}Y^{(j)}_{t_n}} }}{{Q^{(i)}_{t_n}}}Y^{(i)}_{t_n}{\rm{d}}s} 
            + \sum\limits_{{j_1} = 1}^m {\int_{t_n}^{t_{n+1}}  {
            \bigg(
            \tfrac{{\sigma^{(ij_1)} + \sigma^{(ij_1)}\sum\limits_{{j_2} = 1}^m {\sigma^{(ij_2)}\big(W^{(j_2)}_s - W^{(j_2)}_{\lfloor s \rfloor_N}\big)} }}{{Q^{(i)}_{t_n}}}Y^{(i)}_{t_n}
            \bigg)
            {\rm{d}}W^{(j_1)}_s} }.
        \end{aligned} 
    \end{equation}
}
    For any $t\in[0,T]$, one can thus define a continuous version of \eqref{eq:with_Qi_equation} as follows:
    \begin{equation}
    \begin{aligned}\label{eq:conti_glob_version}
            Y^{(i)}_t  
            & =
            Y^{(i)}_0
            +
            \int_{{0}}^t {\tfrac{{{b^{(i)}} - \sum\limits_{j = 1}^d {{a^{(ij)}}Y^{(j)}_{\lfloor s \rfloor_N}} }}{{Q^{(i)}_{\lfloor s \rfloor_N}}}Y^{(i)}_{\lfloor s \rfloor_N}{\rm{d}}s} 
            + \sum\limits_{{j_1} = 1}^m {\int_{{0}}^t {
            \bigg(
            \tfrac{{\sigma^{(ij_1)} + \sigma^{(ij_1)}\sum\limits_{{j_2} = 1}^m {\sigma^{(ij_2)}\big(W^{(j_2)}_s - W^{(j_2)}_{\lfloor s \rfloor_N}\big)} }}{{Q^{(i)}_{\lfloor s \rfloor_N}}}Y^{(i)}_{\lfloor s \rfloor_N}
            \bigg)
            {\rm{d}}W^{(j_1)}_s} }.
        \end{aligned}
    \end{equation}
    
To prove the strong convergence rate of the scheme, the  following lemma is also indispensable.
    \begin{lemma}[Bounded moments]\label{lem:LV_bounded_moment}
        Let all conditions in Proposition \ref{propo:Positivity_milstein_method} hold.
        Then for any $p \geq 1$, there exists a positive constant $C_p$ independent of  $h$, such that the numerical approximations produced by \eqref{eq:with_Qi_equation} obey
        \begin{equation}\label{eq:LV_bounded_moment}
            \sup_{1\leq n \leq N} \mathbb{E}[|Y_{t_n}|^p]\leq C_p.
        \end{equation}
    \end{lemma}
    \indent \textbf{Proof}. For any $p\in \mathbb{Z}^+$, by \eqref{eq:with_Qi_equation} and binomial expansion we have that for $i=1,...,d,\ j,j_1,j_2=1,...,m$,
    \begin{align*}
            &\mathbb{E}[|Y^{(i)}_{t_{n+1}}|^p]
            \\
            &=\mathbb{E}\big[|Y^{(i)}_{t_n}|^p|(Q^{(i)}_{t_n})^{-1}|^p\big]\times \mathbb{E}\Big[\Big(1 + \sum_{j=1}^{m}{\sigma^{(ij)}}\Delta W^{(j)}_{t_n}+\tfrac{1}{2}\sum_{j_1=1}^{m}\sum_{j_2=1}^{m}\sigma^{(ij_1)}\sigma^{(ij_2)}\Delta W^{(j_1)}_{t_n}\Delta W^{(j_2)}_{t_n}\Big)^p\Big]\\
            &\leq (1+Ch)^p\mathbb{E}[|Y^{(i)}_{t_n}|^p]\sum_{k=0}^p \binom{p}{k}\mathbb{E}\Big[\Big(\sum_{j=1}^{m}{\sigma^{(ij)}}\Delta W^{(j)}_{t_n}+\tfrac{1}{2}\sum_{j_1=1}^{m}\sum_{j_2=1}^{m}\sigma^{(ij_1)}\sigma^{(ij_2)}\Delta W^{(j_1)}_{t_n}\Delta W^{(j_2)}_{t_n}\Big)^k\Big]
            \\
            &\leq (1+Ch)^p\mathbb{E}[|Y^{(i)}_{t_n}|^p] 
            \\
            &\ \ \ \times \bigg(1+C_p\sum_{k=1}^p \mathbb{E}\Big[\Big(\sum_{j=1}^{m}{\sigma^{(ij)}}\Delta W^{(j)}_{t_n}+\tfrac{1}{2}\sum_{j_1=1}^{m}\sum_{j_2=1}^{m}\sigma^{(ij_1)}\sigma^{(ij_2)}\Delta W^{(j_1)}_{t_n}\Delta W^{(j_2)}_{t_n}\Big)^k\Big]\bigg),
      \stepcounter{equation}\tag{\theequation}
    \end{align*} 
    where $\binom{p}{k}:=\tfrac{p!}{k!(p-k)!}$ are the coefficients of binomial expansion.
    By iteration one concludes that for any $n=0, 1,...,N-1$,
    \begin{equation}
        \begin{aligned}
            \mathbb{E}[|Y^{(i)}_{t_{n+1}}|^p]&\leq (1+Ch)^{p+1}\mathbb{E}[|Y^{(i)}_{t_n}|^p]
            \leq (1+Ch)^{(p+1)(n+1)}
            \mathbb{E} [|Y^{(i)}_0|^p]
            \leq e^{2(p+1)CT}|x^{(i)}_0|^p,
        \end{aligned} 
    \end{equation}
    which finishes the proof for positive integer $p$. Thanks to the H\"{o}lder inequality, the inequality \eqref{eq:LV_bounded_moment} also holds true
    for any non-integer $p\geq 1$.\qed

Now we are well-prepared to show the order-one pathwise uniformly strong convergence of the proposed scheme.
    \begin{theorem}[Order-one pathwise uniformly strong convergence]
    \label{thm:LV_convergence_rate}
        Let Assumption \ref{assump:LV_model} be satisfied  and assume 
        $0<h \leq \min_{1\leq i\leq d}
        \{
        \tfrac{1}{\gamma(b_i - \frac12\sum_{j=1}^{m}
                {(\sigma^{(ij)})^2})
                    \vee 0
                }
               \wedge
               T\}$
        with $\gamma>1$. Let $X_s$ and $Y_s$ be the exact solution and numerical solution defined by \eqref{eq:one_dim_LV_model} and \eqref{eq:conti_glob_version},
               respectively. 
               Then for any $r>0$,
        \begin{equation}\label{eq:LV-model-main-result}
            \mathbb{E}
            \Big[
            \sup _{t\in [0,T]} |X_t-Y_t|^r 
            \Big] 
            \leq 
            Ch^{r}.
        \end{equation}
    \end{theorem}
    \indent \textbf{Proof}. Define
      \[
        U_{0}(x):=v(1+|x|^2)^{1/2},
        \quad
        U_{1}(x):=v(1+|x|^2)^{ - \frac{1}{2}}\big(\underline{a}d^{-\frac{1}{2}}-\tfrac{vm\bar\sigma^{2}}{2}\big)|x|^3+ v(1+|x|^2)^{1/2}.
      \] 
    Here $v>0$ is chosen to be some small constant such that $ (\underline{a}d^{-\frac{1}{2}}-\tfrac{vm\bar\sigma^{2}}{2})>0$, 
    which in turn ensures $U_{1}(x)$ is positive  and there exists an $\epsilon >0$ satisfying  $\epsilon |x|^2\leq U_{1}(x)$.
    It is easy to check that 
    conditions $(1)$ and $(2)$ in Proposition \ref{prop:LV_from_main_corol} are fulfilled. To validate condition $(3)$, we derive
        \begin{align*}
            &U_{0}^{\prime}(x) f(x)+\tfrac{1}{2} \operatorname{trace}\Big(g(x) g(x)^{*}\operatorname{Hess}_x(U_0(x))\Big)+\tfrac{1}{2}\left|g(x)^{*}(\nabla U_{0})(x)\right|^{2}+U_{1}(x)-\alpha U_{0}(x)
            \\
            &=\sum\limits_{i = 1}^d {v(1+|x|^2)^{-\frac{1}{2}}{(x^{(i)})^2}({b^{(i)}} - \sum\limits_{j = 1}^d {{a^{(ij)}}{x^{(j)}}} )}  + \tfrac{1}{2}\sum\limits_{i = 1}^d {\sum\limits_{k = 1}^m {v(|1+|x|^2)^{-\frac{1}{2}}{(x^{(i)})^2}(\sigma^{(ik)})^2} }  
            \\
            &\ \ \ - \tfrac{1}{2}\sum\limits_{i = 1}^d {\sum\limits_{j = 1}^d {\sum\limits_{k = 1}^m {v(|1+|x|^2)^{-\frac{3}{2}}{(x^{(i)})^2}{(x^{(j)})^2}{\sigma^{(ik)}}{\sigma^{(jk)}}} } }
            + {U_1}(x)   
            \\
            &\ \ \  + \tfrac{1}{2}\sum\limits_{k = 1}^m {\Big(\sum\limits_{i = 1}^d {v(1+|x|^2)^{ - \frac{1}{2}}{(x^{(i)})^2}{\sigma^{(ik)}}} } {\Big)^2}- \alpha v(1+|x|^2)^{  \frac{1}{2}}
            \\
            &\leq \sum\limits_{i = 1}^d {v(1+|x|^2)^{-\tfrac{1}{2}}{(x^{(i)})^2}} \bar b + \tfrac{1}{2}\sum\limits_{i = 1}^d {\sum\limits_{k = 1}^m {v(1+|x|^2)^{ - \frac{1}{2}}{(x^{(i)})^2}{{\bar \sigma }^2}  } } 
            \\
            &\ \ \ +\tfrac{1}{2}\sum\limits_{i = 1}^d {\sum\limits_{j = 1}^d {\sum\limits_{k = 1}^m {v(1+|x|^2)^{ - \frac{3}{2}}(x^{(i)})^2(x^{(j)})^2{{\bar \sigma }^2}} } }-\alpha v(1+|x|^2)^{  \frac{1}{2}}
            \\
            &\ \ \ + \tfrac{m}{2} {{v^2}{\bar \sigma }^2|x{|^4}(1+|x|^2)^{ - 1} - }  {v(1+|x|^2)^{ - \frac{1}{2}}\underline a \sum\limits_{i = 1}^d{(x^{(i)})^3} }+U_1(x)
            \\
            &\leq  v\bar b (1+|x|^2)^{\frac{1}{2}}+mv \bar\sigma^2(1+|x|^2)^{\frac{1}{2}} -\alpha v(1+|x|^2)^{\frac{1}{2}}+
            v(1+|x|^2)^{ - \frac{1}{2}}\big(\tfrac{mv}{2}{\bar \sigma }^2-\underline a d^{-\frac{1}{2}}\big)|x|^3
            +U_1(x)
            \\
            &= (\bar b+m\bar\sigma^2-\alpha +1)v(1+|x|^2)^{\frac{1}{2}}
            \leq 0,
        \stepcounter{equation}\tag{\theequation}
            \label{eq:LV_ex_inte_proof}
      \end{align*} 
      where $\alpha > 0$ is chosen to be large enough
      so that $\bar b+m\bar\sigma^2-\alpha +1 \leq 0$ 
      and the  condition $(3)$ in Proposition \ref{prop:LV_from_main_corol} is hence validated. Furthermore, for any $x,y\in (\mathbb{R}^d)^+$, it holds that
    \begin{equation}
        \begin{aligned}
            &
            \langle x-y, f(x)-f(y)\rangle
            =\langle x-y,\text{diag}(x)b-\text{diag}(x)Ax-\text{diag}(y)b+\text{diag}(y)Ay\rangle
            \\
            & \quad
            =\langle x-y,\text{diag}(x-y)b\rangle-\langle x-y,\text{diag}(x)A(x-y)\rangle-\langle x-y,\text{diag}(x-y)Ay\rangle
            \\
            &
            \quad
            \leq |b||x-y|^{2}+|x-y|^{2}\| \text{diag}(x)A \|-\sum_{i=1}^{d}(x^{(i)}-y^{(i)})^2\sum_{j=1}^{d}a^{(ij)}y^{(j)}
            \\
            & \quad
            \leq 
            ( |b|+ \Vert A \Vert |x| ) |x-y|^{2}.\notag
        \end{aligned}
    \end{equation}
This implies condition $(4)$ in Proposition \ref{prop:LV_from_main_corol} with
    $K_{\eta} \geq \tfrac{1}{4\eta \epsilon} $. 
    By observing that for 
    $ i=1,...,d,\ j=1,...,m,$
    \begin{equation}\label{eq:notation-a-b}
        a^{(i)}(s)=\tfrac{{{b^{(i)}} - \sum\limits_{j = 1}^d {{a^{(ij)}}Y^{(j)}_{\lfloor s \rfloor_N}} }}{{Q^{(i)}_{\lfloor s \rfloor_N}}}Y^{(i)}_{\lfloor s \rfloor_N},
        \quad
        b^{(ij)}(s)=\tfrac{{\sigma^{(ij)} + \sigma^{(ij)}\sum\limits_{{j_1} = 1}^m {\sigma^{(ij_1)}\big(W^{(j_1)}_s - W^{(j_1)}_{\lfloor s \rfloor_N}\big)} }}{{Q^{(i)}_{\lfloor s \rfloor_N}}}Y^{(i)}_{\lfloor s \rfloor_N},
    \end{equation}
    the condition $(5)$ in Proposition \ref{prop:LV_from_main_corol} is therefore satisfied  due to Lemma \ref{lem:LV_bounded_moment} and \eqref{eq:Q-inverse-boundedness}.
Now all conditions of Proposition \ref{prop:LV_from_main_corol} have been confirmed. As a consequence, we arrive at the assertion \eqref{eq:other_methods}.
Following the notation of \eqref{eq:notation-a-b}
and recalling $f (x) := \text{diag}( x ) 
( b - A x)$, $g (x) := \text{diag}( x )  \sigma $,
we use Lemma \ref{lem:LV_bounded_moment} and the H{\"older} inequality to show that, for any $p\geq4$,
    \begin{equation}\label{eq:LV_last_proof_01}
        \begin{aligned}
            \Big\|{f^{(i)}(Y_{\lfloor s \rfloor_N}) - a^{(i)}(s)}\Big\|_{{L^{p}(\Omega;\mathbb{R})}}
            &
            =\Big\| \big(1-\tfrac{1}{Q^{(i)}_{\lfloor s \rfloor_N}}\big)\big(b^{(i)}- \sum\limits_{j = 1}^d{{a^{(ij)}}Y^{(j)}_{\lfloor s \rfloor_N}\big)}Y^{(i)}_{\lfloor s \rfloor_N}\Big\|_{{L^{p}(\Omega;\mathbb{R})}}
            \\
            &\leq \Big\| \big(1-\tfrac{1}{Q^{(i)}_{\lfloor s \rfloor_N}}\big)\Big\|_{{L^{2p}(\Omega;\mathbb{R})}} \Big\|\big(b^{(i)}- \sum\limits_{j = 1}^d{{a^{(ij)}}Y^{(j)}_{\lfloor s \rfloor_N}\big)}Y^{(i)}_{\lfloor s \rfloor_N}\Big\|_{{L^{2p}(\Omega;\mathbb{R})}}
            \\ &
            \leq Ch.
      \end{aligned}
    \end{equation}
Meanwhile, expanding $Y^{(i)}(s)$ at $s={\lfloor s \rfloor_N}$ by It\^{o}'s formula gives
        \begin{align*}
            &\Big\|{g^{(ij)}({Y_s})-b^{(ij)}(s)} \Big\|_{{L^{p}(\Omega;\mathbb{R})}} 
            \\
            &
            = 
            \bigg\|Y^{(i)}_{\lfloor s \rfloor_N}{\sigma^{(ij)}} + \tfrac{{{b^{(i)}} - \sum\limits_{k = 1}^d {{a^{(ik)}}Y^{(k)}_{\lfloor s \rfloor_N}} }}{{Q^{(i)}_{\lfloor s \rfloor_N}}}Y^{(i)}_{\lfloor s \rfloor_N}(s - {\lfloor s \rfloor_N}){\sigma^{(ij)}}+\tfrac{{Y^{(i)}_{\lfloor s \rfloor_N}}{\sigma^{(ij)}}}{{Q^{(i)}_{\lfloor s \rfloor_N}}}\sum\limits_{{j_1} = 1}^m \Big({{\sigma^{(ij_1)}}} (W^{(j_1)}_s - W^{(j_1)}_{\lfloor s \rfloor_N})
            \\ 
            &\ \ \ +\tfrac{1}{2}\sum\limits_{{j_2} = 1}^m {{\sigma^{(ij_1)}}{\sigma^{(ij_2)}}(W^{(j_1)}_s - W^{(j_1)}_{\lfloor s \rfloor_N})(W^{(j_2)}_s - W^{(j_2)}_{\lfloor s \rfloor_N})}  - \tfrac{1}{2}(\sigma^{(ij_1)})^2(s - {\lfloor s \rfloor_N})\Big)
            \\
            &\ \ \ -\tfrac{{Y^{(i)}_{\lfloor s \rfloor_N}}}{{Q^{(i)}_{\lfloor s \rfloor_N}}}{\sigma^{(ij)}}\Big(1 + \sum\limits_{{j_1} = 1}^m {{\sigma^{(ij_1)}}(W^{(j_1)}_s - W^{(j_1)}_{\lfloor s \rfloor_N})} \Big)\bigg\|_{{L^{p}(\Omega;\mathbb{R})}}
            \\
            &
            = \bigg\|Y^{(i)}_{\lfloor s \rfloor_N}{\sigma^{(ij)}}\tfrac{Q^{(i)}_{\lfloor s \rfloor_N}-1}{{Q^{(i)}_{\lfloor s \rfloor_N}}}
            +
            \tfrac{{{b^{(i)}} - \sum\limits_{k = 1}^d {{a^{(ik)}}Y^{(k)}_{\lfloor s \rfloor_N}} }}{{Q^{(i)}_{\lfloor s \rfloor_N}}}Y^{(i)}_{\lfloor s \rfloor_N}(s - {\lfloor s \rfloor_N}){\sigma^{(ij)}}
            \\
            &\ \ \ +\tfrac{{Y^{(i)}_{\lfloor s \rfloor_N}{\sigma^{(ij)}}}}{{2Q^{(i)}_{\lfloor s \rfloor_N}}}\sum\limits_{{j_1} = 1}^m\Big(\sum\limits_{{j_2} = 1}^m {{\sigma^{(ij_1)}}{\sigma^{(ij_2)}}(W^{(j_1)}_s - W^{(j_1)}_{\lfloor s \rfloor_N})(W^{(j_2)}_s - W^{(j_2)}_{\lfloor s \rfloor_N})}  - (\sigma^{(ij_1)})^2(s - {\lfloor s \rfloor_N})\Big) 
            \bigg\|_{{L^{p}(\Omega;\mathbb{R})}}
            \\
            &\leq 
            \Big\|Y^{(i)}_{\lfloor s \rfloor_N}{\sigma^{(ij)}}\tfrac{Q^{(i)}_{\lfloor s \rfloor_N}-1}{{Q^{(i)}_{\lfloor s \rfloor_N}}}\Big\|_{{L^{p}(\Omega;\mathbb{R})}}
            +
            \bigg\| 
            \tfrac{{{b^{(i)}} - \sum\limits_{k = 1}^d {{a^{(ik)}}Y^{(k)}_{\lfloor s \rfloor_N}} }}{{Q^{(i)}_{\lfloor s \rfloor_N}}}Y^{(i)}_{\lfloor s \rfloor_N}(s - {\lfloor s \rfloor_N}){\sigma^{(ij)}}
            \bigg\|_{{L^{p}(\Omega;\mathbb{R})}}
            \\
            &\ \ \ 
            + \bigg\|\tfrac{Y^{(i)}_{{\lfloor s \rfloor_N}} \sigma^{(ij)}}{2 Q^{(i)}_{{\lfloor s \rfloor_N}}} \sum_{j_{1}=1}^{m}\bigg(\sum_{j_{2}=1}^{m} \sigma^{(ij_{1})} \sigma^{(ij_{2})}\big(W^{(j_{1})}_s-W^{(j_{1})}_{{\lfloor s \rfloor_N}}\big)\big(W^{(j_{2})}_s-W^{(j_{2})}_{{\lfloor s \rfloor_N}}\big)-(\sigma^{(ij_1)})^2\left(s-{\lfloor s \rfloor_N}\right)\bigg)\bigg\|_{{L^{p}(\Omega;\mathbb{R})}}
            \\
            &\leq Ch.
            \stepcounter{equation}\tag{\theequation}
            \label{eq:LV_last_proof_02}
        \end{align*}
    Combining \eqref{eq:LV_last_proof_01}-\eqref{eq:LV_last_proof_02} with \eqref{eq:other_methods} finally yields the desired assertion \eqref{eq:LV-model-main-result}.\qed

    \textbf{Numerical experiments.}
    {\color{blue}
      Now we present some numerical experiments to confirm the theoretical results. Let $d=m=2$ and consider the stochastic LV competition model with  coefficients 
        \begin{equation}\label{eq:numer_exper_coef}
            b=\left( {\begin{array}{*{20}{c}}
            1\\
                0.5
            \end{array}} \right),
            A=\left( {\begin{array}{*{20}{c}}
            1&{0.5}\\
            0&{0.5}
            \end{array}} \right),
           \sigma=\left( {\begin{array}{*{20}{c}}
            1&0\\
            0&{0.75}
            \end{array}} \right),
        Y(0)=\left( {\begin{array}{*{20}{c}}
            1\\
            3
        \end{array}} \right).
        \end{equation}    
Let $T = 1$, $N = 2^{k}, k = 6, 7,..., 11$ and regard the fine approximations with $h_{\text{exact}} = 2^{-14}$ as the "true" solution.
We consider the mean-square approximation errors and take $M = 5000$ Monte Carlo sample paths to approximate the expectation. A comparison of strong convergence rates for
our method, the Lamperti transformed EM method in \cite{LI2023107260} and the truncated EM method in \cite{mao2021positivity} is presented in Figure \ref{fig:LV_strong_convergence_rate}. One can clearly observe order-one convergence of  our method and the Lamperti transformed EM method, and order $\tfrac12$ convergence of the truncated EM method.
    \begin{figure}[H]
      \centering
      \includegraphics[scale=0.21]{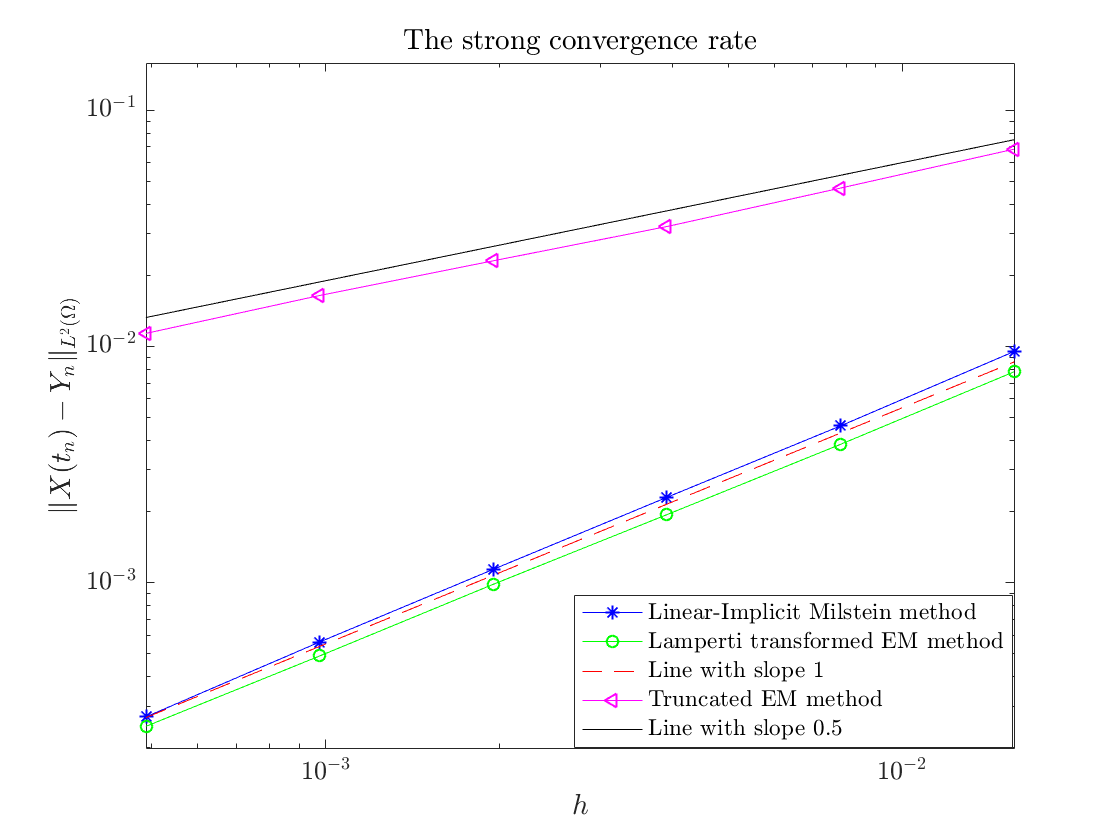}
      \caption{\small A comparison of strong convergence rates for LV competition  model
      \label{fig:LV_strong_convergence_rate}}
    \end{figure}
In addition to the strong convergence rate, we would also like to investigate the dynamic preservation of the proposed method. As shown by \cite{bahar2004stochastic,mao2007stochastic}, under Assumption \ref{assump:LV_model}, the exact solution $\{X_t\}_{t\geq 0}$ of \eqref{eq:L_V_model} admits an ultimate boundedness property,
    i.e., there exist two positive constants $C_1,C_2$ independent of $X_0$ such that 
    \begin{equation}\label{eq:LV_ultimate_bound}
             \limsup_{t\rightarrow  \infty} \mathbb{E}[|X_t|] \leq C_1 \ \
             \text{and}\ \
         \limsup_{t\rightarrow  \infty} \tfrac{1}{t} 
             \int_0^t \mathbb{E}[|X_s|^2] {\mathrm d}s \leq C_2.
         \end{equation}
Figure \ref{fig:LV_longtime_moment} displays moments over long-time interval $[0, T], T=500$ of the linear-implicit Milstein method using small stepsize $h=2^{-7}$ and large stepsize $h=1$. It is observed that, the numerical approximations produced by the linear-implicit Milstein method remain bounded after a long time, even for a large stepsize $h=1$. This recovers the property \eqref{eq:LV_ultimate_bound} of the exact solution, which can be theoretically explained as follows. 
For any $p \geq 1$ and $0<h < \min_{1\leq i\leq d}
        \{
        \tfrac{1}{(b_i - \frac12\sum_{j=1}^{m}
                {(\sigma^{(ij)})^2})
                    \vee 0
                }
            \}$,
recall \eqref{eq:Qi_def} and \eqref{eq:with_Qi_equation},
it holds that
 \[
      \begin{aligned}
        \mathbb E \big[|Y^{(i)}_{t_{n+1}}|^p\big]&=\mathbb E \big[\big|\tfrac{Y^{(i)}_{t_n}}{Q^{(i)}_{t_n}}\big|^p\big] \mathbb E\Big[ \Big(\tfrac{1}{2}\big(1+\sum_{j=1}^{m}{\sigma^{(ij)}}\Delta W^{(j)}_{t_n}\big)^2+\tfrac{1}{2}\Big)^p\Big]
        \leq C_{p,h}(\tfrac{1}{ha^{(ii)}})^p 
        <+\infty,
      \end{aligned}
    \]
    for any $n\in \mathbb{N}$. 
    
Another significant dynamic of the model \eqref{eq:L_V_model} is permanence and extinction. Permanence means that the species will persist and  extinction means that the species will eventually become extinct (see \cite[Definition 3.6, Section 4]{li2009population} for precise definitions). According to \cite[Theorem 3.7]{li2009population}, in our setting \eqref{eq:numer_exper_coef} the species will be permanent. In Figure \ref{fig:LV_permanent_property}, we draw $M=500$ sample paths of the proposed linear implicit Milstein scheme over time interval $[0,500]$ with $h=2^{-3}$. Evidently, the numerical approximations reproduce the dynamic of the permanence of the original model. 
When the noise intensity increases from $\sigma$ to $2 \sigma$, in view of \cite[Corollary 2]{li2009population}, the species will eventually become extinct. 
Figure \ref{fig:LV_extinction_property} presents $M=500$ sample paths of the proposed approximation method over time interval $[0,500]$ with $h=2^{-3}$. 
There one can see that the numerical approximations tend to zero, reproducing the dynamic of the extinct of the original model.
For both cases, all paths stay in $(\mathbb{R}^2)^+$, confirming the positivity preserving of the proposed method.

    \begin{figure}[H]
    \centering  
        \subfigure[]{
        \includegraphics[width=6.6cm,height = 5.4cm]{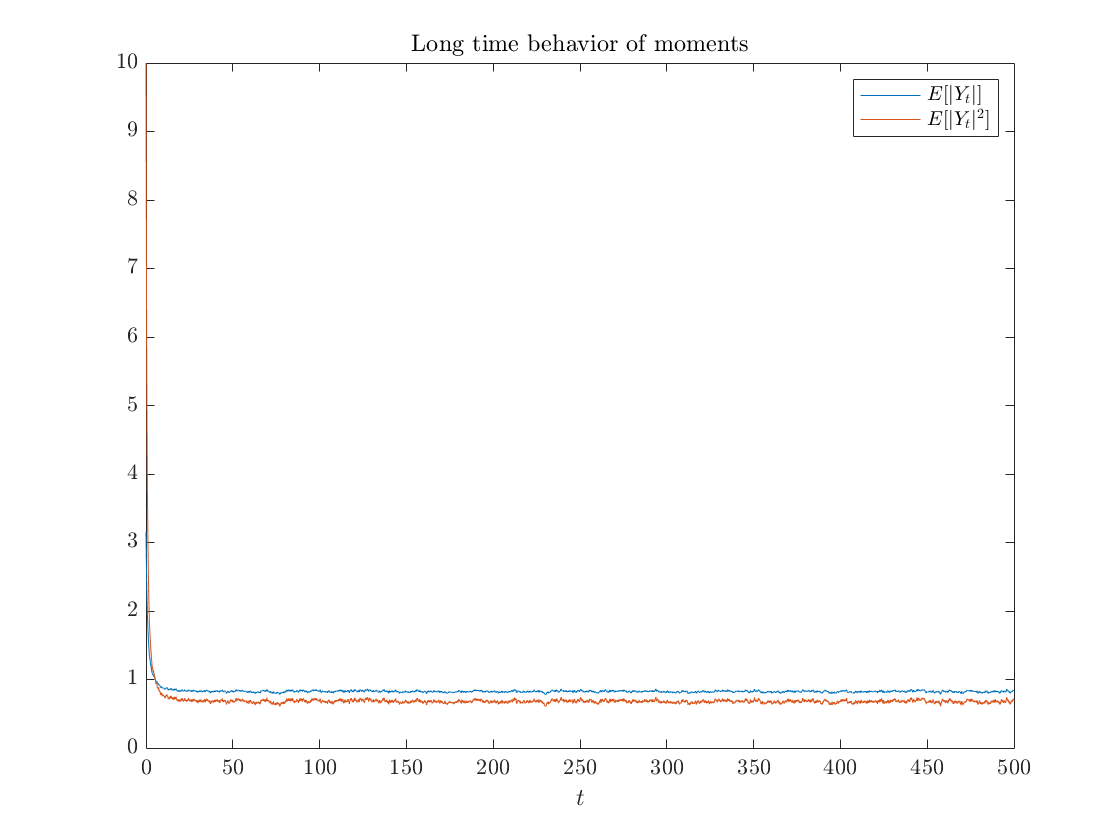}}
        \subfigure[]{
        \includegraphics[width=6.6cm,height = 5.4cm]{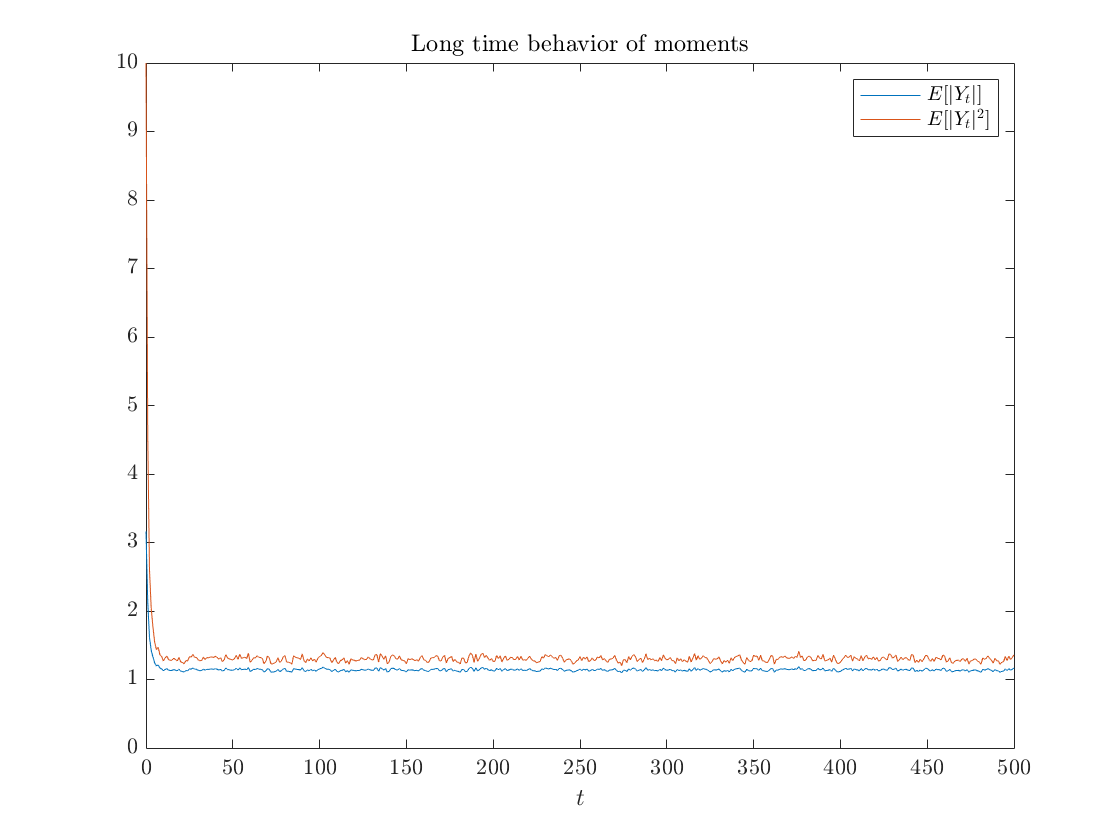}}
        \caption{\small Evolution of numerical moments: (a) $h=2^{-7}$; (b) $h=1$ }
        \label{fig:LV_longtime_moment}
\end{figure}

    \begin{figure}[H]
      \centering
      \includegraphics[scale=0.382]{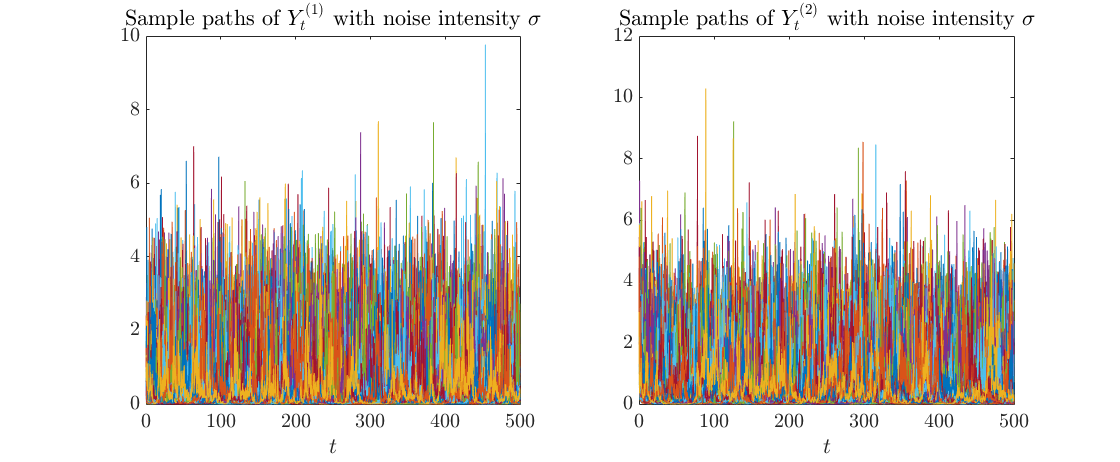}
      \caption{\small Permanent performance of linear-implicit Milstein method with $h=2^{-3}$
      \label{fig:LV_permanent_property}}
    \end{figure}

        \begin{figure}[H]
      \centering
      \includegraphics[scale=0.382]{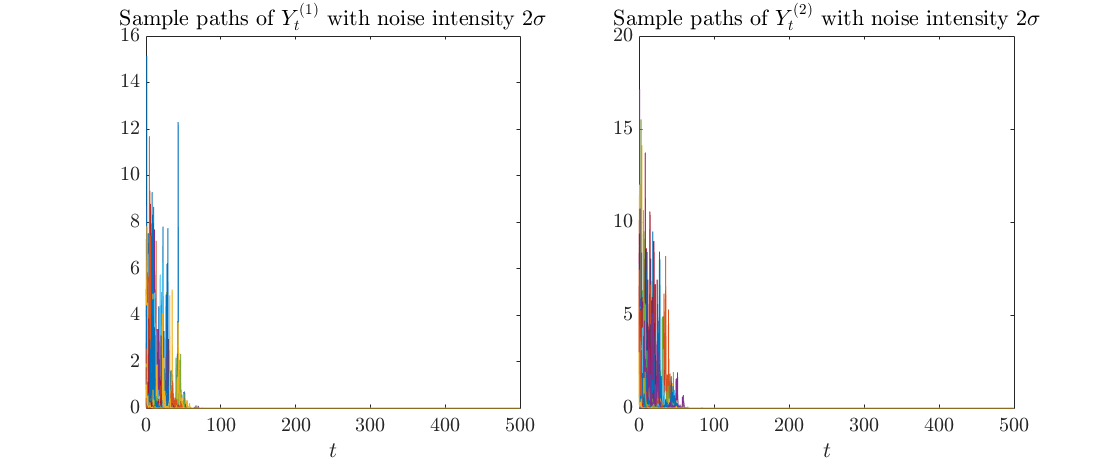}
      \caption{\small Extinct performance of linear-implicit Milstein method with $h=2^{-3}$
      \label{fig:LV_extinction_property}}
    \end{figure}
}
\section{Conclusion} 
\label{section:conclusion}
{\color{black}
In this paper we successfully reveal order-one strong convergence of two kinds of numerical methods for
several SDEs without globally monotone coefficients,
which fills the gap left by  \cite{Hutzenthaler2020}.}
This is accomplished by developing some
new perturbation estimates and some more careful estimates.
Numerical experiments are also provided to support the theoretical findings.
As an ongoing project, we propose and analyze higher order
(strong order $1$ and $1.5$) time-stepping schemes for general SDEs without globally monotone coefficients. 

\bibliographystyle{abbrv}

\bibliography{perturbation_order_one_reference}

\begin{thebibliography}{10}

\bibitem{alfonsi2013strong}
A.~Alfonsi.
\newblock {Strong order one convergence of a drift implicit Euler scheme:
  Application to the CIR process}.
\newblock {\em Statistics \& Probability Letters}, 83(2):602--607, 2013.

\bibitem{andersson2017mean}
A.~Andersson and R.~Kruse.
\newblock {Mean-square convergence of the BDF2-Maruyama and backward Euler
  schemes for SDE satisfying a global monotonicity condition}.
\newblock {\em BIT Numerical Mathematics}, 57(1):21--53, 2017.

\bibitem{bahar2004stochastic}
A.~Bahar and X.~Mao.
\newblock Stochastic delay population dynamics.
\newblock {\em International Journal of Pure and Applied Mathematics},
  11:377--400, 2004.

\bibitem{beyn2017stochastic}
W.-J. Beyn, E.~Isaak, and R.~Kruse.
\newblock Stochastic {C-stability and B-consistency of explicit and implicit
  Milstein-type} schemes.
\newblock {\em Journal of Scientific Computing}, 70(3):1042--1077, 2017.

\bibitem{brehier2023approximation}
C.-E. Brehier.
\newblock {Approximation of the invariant distribution for a class of ergodic
  SDEs with one-sided Lipschitz continuous drift coefficient using an explicit
  tamed Euler scheme}.
\newblock {\em ESAIM: Probability and Statistics}, 27:841--866, 2023.

\bibitem{cai2023positivity}
Y.~Cai, Q.~Guo, and X.~Mao.
\newblock {Positivity preserving truncated scheme for the stochastic
  Lotka--Volterra model with small moment convergence}.
\newblock {\em Calcolo}, 60(2):24, 2023.

\bibitem{cox2024local}
S.~Cox, M.~Hutzenthaler, and A.~Jentzen.
\newblock {\em {Local Lipschitz continuity in the initial value and strong
  completeness for nonlinear stochastic differential equations}}, volume 296.
\newblock American Mathematical Society, 2024.

\bibitem{cui2022density}
J.~Cui, J.~Hong, and D.~Sheng.
\newblock Density function of numerical solution of splitting {AVF} scheme for
  stochastic {Langevin} equation.
\newblock {\em Mathematics of Computation}, 91(337):2283--2333, 2022.

\bibitem{Fang2020}
W.~Fang and M.~B. Giles.
\newblock {Adaptive Euler--Maruyama method for {SDEs} with nonglobally
  {Lipschitz} drift}.
\newblock {\em The Annals of Applied Probability}, 30(2):526--560, 2020.

\bibitem{higham2000mean}
D.~J. Higham.
\newblock Mean-square and asymptotic stability of the stochastic theta method.
\newblock {\em SIAM Journal on Numerical Analysis}, 38(3):753--769, 2000.

\bibitem{higham2002strong}
D.~J. Higham, X.~Mao, and A.~M. Stuart.
\newblock Strong convergence of {Euler-type} methods for nonlinear stochastic
  differential equations.
\newblock {\em SIAM Journal on Numerical Analysis}, 40(3):1041--1063, 2002.

\bibitem{hong2021positivity}
J.~Hong, L.~Ji, X.~Wang, and J.~Zhang.
\newblock Positivity-preserving symplectic methods for the stochastic
  {Lotka--Volterra} predator-prey model.
\newblock {\em BIT}, 62:493--520, 2022.

\bibitem{hutzenthaler2011convergence}
M.~Hutzenthaler and A.~Jentzen.
\newblock Convergence of the stochastic {Euler} scheme for locally {Lipschitz}
  coefficients.
\newblock {\em Foundations of Computational Mathematics}, 11(6):657--706, 2011.

\bibitem{hutzenthaler2015numerical}
M.~Hutzenthaler and A.~Jentzen.
\newblock {\em {Numerical approximations of stochastic differential equations
  with non-globally Lipschitz continuous coefficients}}, volume 236.
\newblock American Mathematical Society, 2015.

\bibitem{Hutzenthaler2020}
M.~{Hutzenthaler} and A.~{Jentzen}.
\newblock On a perturbation theory and on strong convergence rates for
  stochastic ordinary and partial differential equations with nonglobally
  monotone coefficients.
\newblock {\em The Annals of Probability}, 48(1):53--93, 2020.

\bibitem{hutzenthaler2011strong}
M.~Hutzenthaler, A.~Jentzen, and P.~E. Kloeden.
\newblock Strong and weak divergence in finite time of {Euler's} method for
  stochastic differential equations with non-globally {Lipschitz} continuous
  coefficients.
\newblock {\em Proceedings of the Royal Society A: Mathematical, Physical and
  Engineering Sciences}, 467(2130):1563--1576, 2011.

\bibitem{hutzenthaler2012strong}
M.~Hutzenthaler, A.~Jentzen, and P.~E. Kloeden.
\newblock Strong convergence of an explicit numerical method for {SDEs} with
  nonglobally {Lipschitz} continuous coefficients.
\newblock {\em The Annals of Applied Probability}, 22(4):1611--1641, 2012.

\bibitem{hutzenthaler2018exponential}
M.~Hutzenthaler, A.~Jentzen, and X.~Wang.
\newblock Exponential integrability properties of numerical approximation
  processes for nonlinear stochastic differential equations.
\newblock {\em Mathematics of Computation}, 87(311):1353--1413, 2018.

\bibitem{Kelly2023}
C.~Kelly, G.~Lord, and F.~Sun.
\newblock Strong convergence of an adaptive time-stepping {Milstein} method for
  {SDEs} with monotone coefficients.
\newblock {\em BIT Numerical Mathematics}, 63:33, 2023.

\bibitem{1992Numerical}
P.~E. Kloeden and E.~Platen.
\newblock {\em Numerical Solution of Stochastic Differential Equations}.
\newblock Berlin: Springer-Verlag, 1992.

\bibitem{kumar2019milstein}
C.~Kumar and S.~Sabanis.
\newblock On {Milstein} approximations with varying coefficients: the case of
  super-linear diffusion coefficients.
\newblock {\em BIT Numerical Mathematics}, 59(4):929--968, 2019.

\bibitem{li2009population}
X.~Li and X.~Mao.
\newblock {Population dynamical behavior of non-autonomous Lotka-Volterra
  competitive system with random perturbation}.
\newblock {\em Discrete and Continuous Dynamical Systems-Series A},
  24(2):523--593, 2009.

\bibitem{Li2019explicit}
X.~Li, X.~Mao, and G.~Yin.
\newblock Explicit numerical approximations for stochastic differential
  equations in finite and infinite horizons: truncation methods, convergence in
  pth moment, and stability.
\newblock {\em IMA Journal of Numerical Analysis}, 39(2):847--892, 2019.

\bibitem{LI2023107260}
Y.~Li and W.~Cao.
\newblock {A positivity preserving Lamperti transformed Euler--Maruyama method
  for solving the stochastic Lotka--Volterra competition model}.
\newblock {\em Communications in Nonlinear Science and Numerical Simulation},
  122, 2023.

\bibitem{liu2013strong}
W.~Liu and X.~Mao.
\newblock Strong convergence of the stopped {Euler--Maruyama} method for
  nonlinear stochastic differential equations.
\newblock {\em Applied Mathematics and Computation}, 223:389--400, 2013.

\bibitem{mao2007stochastic}
X.~Mao.
\newblock {\em Stochastic Differential Equations and Applications}.
\newblock Elsevier, 2007.

\bibitem{mao2015truncated}
X.~Mao.
\newblock The truncated {Euler--Maruyama} method for stochastic differential
  equations.
\newblock {\em Journal of Computational and Applied Mathematics}, 290:370--384,
  2015.

\bibitem{mao2013strong}
X.~Mao and L.~Szpruch.
\newblock {Strong convergence rates for backward Euler--Maruyama method for
  non-linear dissipative-type stochastic differential equations with
  super-linear diffusion coefficients}.
\newblock {\em Stochastics}, 85(1):144--171, 2013.

\bibitem{mao2021positivity}
X.~Mao, F.~Wei, and T.~Wiriyakraikul.
\newblock Positivity preserving truncated {Euler--Maruyama} method for
  stochastic {Lotka--Volterra} competition model.
\newblock {\em Journal of Computational and Applied Mathematics}, 394:113566,
  2021.

\bibitem{mattingly2002ergodicity}
J.~C. Mattingly, A.~M. Stuart, and D.~J. Higham.
\newblock {Ergodicity for SDEs and approximations: locally Lipschitz vector
  fields and degenerate noise}.
\newblock {\em Stochastic processes and their applications}, 101(2):185--232,
  2002.

\bibitem{Milstein2004}
G.~N. {Milstein} and M.~V. {Tretyakov}.
\newblock {\em Stochastic Numerics for Mathematical Physics}.
\newblock Berlin: Springer, 2004.

\bibitem{neuenkirch2014first}
A.~Neuenkirch and L.~Szpruch.
\newblock First order strong approximations of scalar {SDEs} defined in a
  domain.
\newblock {\em Numerische Mathematik}, 128(1):103--136, 2014.

\bibitem{sabanis2016euler}
S.~Sabanis.
\newblock Euler approximations with varying coefficients: the case of
  superlinearly growing diffusion coefficients.
\newblock {\em The Annals of Applied Probability}, 26(4):2083--2105, 2016.

\bibitem{talay2002stochastic}
D.~Talay.
\newblock Stochastic {Hamiltonian} systems: exponential convergence to the
  invariant measure, and discretization by the implicit {Euler} scheme.
\newblock {\em Markov Process. Related Fields}, 8(2):163--198, 2002.

\bibitem{to2000nonlinear}
C.~W. To.
\newblock {\em Nonlinear random vibration: Analytical techniques and
  applications}.
\newblock Crc Press, 2000.

\bibitem{Tretyakov2013fundamental}
M.~V. Tretyakov and Z.~Zhang.
\newblock A fundamental mean-square convergence theorem for {SDE}s with locally
  {L}ipschitz coefficients and its applications.
\newblock {\em SIAM Journal on Numerical Analysis}, 51(6):3135--3162, 2013.

\bibitem{wang2023mean}
X.~Wang.
\newblock {Mean-square convergence rates of implicit Milstein type methods for
  SDEs with non-Lipschitz coefficients}.
\newblock {\em Advances in Computational Mathematics}, 49(3):37, 2023.

\bibitem{wang2013tamed}
X.~Wang and S.~Gan.
\newblock The tamed {Milstein} method for commutative stochastic differential
  equations with non-globally {Lipschitz} continuous coefficients.
\newblock {\em Journal of Difference Equations and Applications},
  19(3):466--490, 2013.

\bibitem{wang2020mean}
X.~Wang, J.~Wu, and B.~Dong.
\newblock {Mean-square convergence rates of stochastic theta methods for SDEs
  under a coupled monotonicity condition}.
\newblock {\em BIT}, 60(3):759--790, 2020.

\bibitem{zong2018convergence}
X.~Zong, F.~Wu, and G.~Xu.
\newblock {Convergence and stability of two classes of theta-Milstein schemes
  for stochastic differential equations}.
\newblock {\em Journal of Computational and Applied Mathematics}, 336:8--29,
  2018.

\end{thebibliography}

\end{document}